\documentclass{article}
\usepackage{amssymb,amsmath,amsthm}
\usepackage[mathscr]{euscript}

\usepackage{graphicx}
\input{epsf}

\usepackage[bbgreekl]{mathbbol}
\DeclareSymbolFontAlphabet{\mathbb}{AMSb}
\DeclareSymbolFontAlphabet{\mathbbb}{bbold}

\newcounter{sec}

\newcounter{punct}[sec]

\def\punct{\refstepcounter{punct}{\arabic{sec}.\arabic{punct}.  }}

\newtheorem{theorem}{Theorem}[sec]
\newtheorem{proposition}[theorem]{Proposition}

\newtheorem{lemma}[theorem]{Lemma}

\newtheorem{corollary}[theorem]{Corollary}
\newtheorem{observation}[theorem]{Observation}

\newtheorem{conjecture}[theorem]{Conjecture}

\def\COUNTERS{\addtocounter{sec}{1}
              \setcounter{punct}{0}
          \setcounter{equation}{0}
          \setcounter{theorem}{0}
            \setcounter{figure}{0}
          }
          
          \def\sm{\smallskip}

\begin{document}

\newcommand{\supp}{\mathop {\mathrm {supp}}\nolimits}
\newcommand{\rk}{\mathop {\mathrm {rk}}\nolimits}
\newcommand{\Aut}{\mathop {\mathrm {Aut}}\nolimits}
\newcommand{\Out}{\mathop {\mathrm {Out}}\nolimits}
\renewcommand{\Re}{\mathop {\mathrm {Re}}\nolimits}
\newcommand{\Inn}{\mathop {\mathrm {Inn}}\nolimits}
\newcommand{\Char}{\mathop {\mathrm {Char}}\nolimits}
\newcommand{\ch}{\cosh}
\newcommand{\sh}{\sinh}
\newcommand{\Sp}{\mathop {\mathrm {Sp}}\nolimits}
\newcommand{\SOS}{\mathop {\mathrm {SO^*}}\nolimits}
\newcommand{\Ams}{\mathop {\mathrm {Ams}}\nolimits}
\newcommand{\Gms}{\mathop {\mathrm {Gms}}\nolimits}

\def\0{\mathbf 0}

\def\ov{\overline}
\def\wh{\widehat}
\def\wt{\widetilde}
\def\pol{\twoheadrightarrow}

\renewcommand{\rk}{\mathop {\mathrm {rk}}\nolimits}
\renewcommand{\Aut}{\mathop {\mathrm {Aut}}\nolimits}
\newcommand{\Ob}{\mathop {\mathrm {Ob}}\nolimits}
\renewcommand{\Re}{\mathop {\mathrm {Re}}\nolimits}
\renewcommand{\Im}{\mathop {\mathrm {Im}}\nolimits}
\newcommand{\sgn}{\mathop {\mathrm {sgn}}\nolimits}
\newcommand{\ver}{\mathop {\mathrm {vert}}\nolimits}
\newcommand{\val}{\mathop {\mathrm {val}}\nolimits}
\newcommand{\edge}{\mathop {\mathrm {edge}}\nolimits}
\newcommand{\germ}{\mathop {\mathrm {germ}}\nolimits}
\newcommand{\PB}{{\mathop {\EuScript {PB}}\nolimits}}

\def\bfa{\mathbf a}
\def\bfb{\mathbf b}
\def\bfc{\mathbf c}
\def\bfd{\mathbf d}
\def\bfe{\mathbf e}
\def\bff{\mathbf f}
\def\bfg{\mathbf g}
\def\bfh{\mathbf h}
\def\bfi{\mathbf i}
\def\bfj{\mathbf j}
\def\bfk{\mathbf k}
\def\bfl{\mathbf l}
\def\bfm{\mathbf m}
\def\bfn{\mathbf n}
\def\bfo{\mathbf o}
\def\bfp{\mathbf p}
\def\bfq{\mathbf q}
\def\bfr{\mathbf r}
\def\bfs{\mathbf s}
\def\bft{\mathbf t}
\def\bfu{\mathbf u}
\def\bfv{\mathbf v}
\def\bfw{\mathbf w}
\def\bfx{\mathbf x}
\def\bfy{\mathbf y}
\def\bfz{\mathbf z}

\def\bfA{\mathbf A}
\def\bfB{\mathbf B}
\def\bfC{\mathbf C}
\def\bfD{\mathbf D}
\def\bfE{\mathbf E}
\def\bfF{\mathbf F}
\def\bfG{\mathbf G}
\def\bfH{\mathbf H}
\def\bfI{\mathbf I}
\def\bfJ{\mathbf J}
\def\bfK{\mathbf K}
\def\bfL{\mathbf L}
\def\bfM{\mathbf M}
\def\bfN{\mathbf N}
\def\bfO{\mathbf O}
\def\bfP{\mathbf P}
\def\bfQ{\mathbf Q}
\def\bfR{\mathbf R}
\def\bfS{\mathbf S}
\def\bfT{\mathbf T}
\def\bfU{\mathbf U}
\def\bfV{\mathbf V}
\def\bfW{\mathbf W}
\def\bfX{\mathbf X}
\def\bfY{\mathbf Y}
\def\bfZ{\mathbf Z}

\def\frA{\mathfrak A}
\def\frB{\mathfrak B}
\def\frD{\mathfrak D}
\def\frS{\mathfrak S}
\def\frL{\mathfrak L}
\def\frP{\mathfrak P}
\def\frQ{\mathfrak Q}
\def\frR{\mathfrak R}
\def\frT{\mathfrak T}
\def\frG{\mathfrak G}
\def\frg{\mathfrak g}
\def\frh{\mathfrak h}
\def\frf{\mathfrak f}
\def\frk{\mathfrak k}
\def\frl{\mathfrak l}
\def\frm{\mathfrak m}
\def\frn{\mathfrak n}
\def\fro{\mathfrak o}
\def\frp{\mathfrak p}
\def\frq{\mathfrak q}
\def\frr{\mathfrak r}
\def\frs{\mathfrak s}
\def\frt{\mathfrak t}
\def\fru{\mathfrak u}
\def\frv{\mathfrak v}
\def\frw{\mathfrak w}
\def\frx{\mathfrak x}
\def\fry{\mathfrak y}
\def\frz{\mathfrak z}

\def\bfw{\mathbf w}

\def\R {{\mathbb {R} }}
 \def\C {{\mathbb C }}
  \def\Z{{\mathbb Z}}
  \def\H{{\mathbb H}}
\def\K{{\mathbb K}}
\def\N{{\mathbb N}}
\def\Q{{\mathbb Q}}
\def\A{{\mathbb A}}
\def\O {{\mathbb O }}
\def\G {{\mathbb G }}

\def\T{\mathbb T}
\def\P{\mathbb P}
\def\SS{\mathbb S}

\def\G{\mathbb I}

\def\cD{\EuScript D}
\def\cL{\EuScript L}
\def\cK{\EuScript K}
\def\cM{\EuScript M}
\def\cN{\EuScript N}
\def\cP{\EuScript P}
\def\cQ{\EuScript Q}
\def\cR{\EuScript R}
\def\cT{\EuScript T}
\def\cW{\EuScript W}
\def\cY{\EuScript Y}
\def\cF{\EuScript F}
\def\cG{\EuScript G}
\def\cZ{\EuScript Z}
\def\cI{\EuScript I}
\def\cB{\EuScript B}
\def\cA{\EuScript A}
\def\cE{\EuScript E}
\def\cC{\EuScript C}

\def\bbA{\mathbb A}
\def\bbB{\mathbb B}
\def\bbD{\mathbb D}
\def\bbE{\mathbb E}
\def\bbF{\mathbb F}
\def\bbG{\mathbb G}
\def\bbI{\mathbb I}
\def\bbJ{\mathbb J}
\def\bbL{\mathbb L}
\def\bbM{\mathbb M}
\def\bbN{\mathbb N}
\def\bbO{\mathbb O}
\def\bbP{\mathbb P}
\def\bbQ{\mathbb Q}
\def\bbS{\mathbb S}
\def\bbT{\mathbb T}
\def\bbU{\mathbb U}
\def\bbV{\mathbb V}
\def\bbW{\mathbb W}
\def\bbX{\mathbb X}
\def\bbY{\mathbb Y}

\def\U{\mathbb{U}}

\def\kappa{\varkappa}
\def\epsilon{\varepsilon}
\def\phi{\varphi}
\def\le{\leqslant}
\def\ge{\geqslant}

\def\B{\mathrm B}

\def\la{\langle}
\def\ra{\rangle}
\def\tri{\triangleright}

\def\lambdA{{\boldsymbol{\lambda}}}
\def\alphA{{\boldsymbol{\alpha}}}
\def\betA{{\boldsymbol{\beta}}}
\def\mU{{\boldsymbol{\mu}}}
\def\nU{{\boldsymbol{\nu}}}
\def\omegA{{\boldsymbol{\omega}}}
\def\GammA{{\mathbbb{\Gamma}}}

\def\const{\mathrm{const}}
\def\rem{\mathrm{rem}}
\def\even{\mathrm{even}}
\def\SO{\mathrm{SO}}
\def\SL{\mathrm{SL}}
\def\PSL{\mathrm{PSL}}
\def\cont{\mathrm{cont}}
\def\Isom{\mathrm{Isom}}
\def\Isolated{\mathrm{Isolated}}
\def\coun{\mathrm{count}}
\def\S{\mathrm{S}}
\def\out{\mathrm{out}}

\def\Un{\operatorname{U}}
\def\GL{\operatorname{GL}}
\def\Mat{\operatorname{Mat}}
\def\End{\operatorname{End}}
\def\Mor{\operatorname{Mor}}
\def\Aut{\operatorname{Aut}}
\def\inv{\operatorname{inv}}
\def\red{\operatorname{red}}
\def\Ind{\operatorname{Ind}}
\def\dom{\operatorname{dom}}
\def\im{\operatorname{im}}
\def\md{\operatorname{mod\,}}
\def\indef{\operatorname{indef}}
\def\Gr{\operatorname{Gr}}
\def\Pol{\operatorname{Pol}}
\def\UU{\operatorname{U}}

\def\arr{\rightrightarrows}
\def\bs{\backslash}

\def\cH{\EuScript{H}}
\def\cO{\EuScript{O}}
\def\cQ{\EuScript{Q}}
\def\cL{\EuScript{L}}
\def\cX{\EuScript{X}}
\def\cJ{\EuScript{J}}

\def\Di{\Diamond}
\def\di{\diamond}

\def\fin{\mathrm{fin}}
\def\ThetA{\boldsymbol {\Theta}}

\def\0{\boldsymbol{0}}

\def\F{\,{\vphantom{F}}_2F_1}
\def\FF{\,{\vphantom{F}}_3F_2}
\def\H{\,\vphantom{H}^{\phantom{\star}}_2 H_2^\star}
\def\HH{\,\vphantom{H}^{\phantom{\star}}_3 H_3^\star}
\def\Ho{\,\vphantom{H}_2 H_2}




\def\abst
{We consider groups $\mathbb{I}$ of isometries  of   ultrametric Urysohn spaces $\mathbb{U}$. Such spaces $\mathbb{U}$ admit transparent
realizations as boundaries of certain $R$-trees and the groups $\mathbb{I}$ are groups of automorphisms of these $R$-trees.
 Denote by $\mathbb{I}[X]\subset \mathbb{I}$ stabilizers of finite subspaces
$X\subset \U$. Double cosets $\mathbb{I}[X]\cdot g\cdot \mathbb{I}[Y]$, where $g\in \mathbb{I}$,
 are enumerated by  ultrametrics on union of spaces $X\cup Y$.  We construct natural associative multiplications on double coset 
spaces $\mathbb{I}[X]\backslash \mathbb{I}/\mathbb{I}[X]$ and, more generally, multiplications
$\mathbb{I}[X]\backslash \mathbb{I}/\mathbb{I}[Y]\,\times\, \mathbb{I}[Y]\backslash \mathbb{I}/\mathbb{I}[Z]\to \mathbb{I}[X]\backslash \mathbb{I}/\mathbb{I}[Z]$.
These operations are a kind of  canonical amalgamations of ultrametric spaces. 
On the other hand, 
this product can be interpreted in terms of partial isomorphisms of certain $R$-trees 
(in particular, we come to an inverse category).
 This allows us
to classify all unitary representations of the groups $\mathbb{I}$ 
and to prove that groups $\mathbb{I}$ have type $I$.
We also describe a universal semigroup compactification of $\mathbb{I}$
 whose image in any unitary representation of $\mathbb{I}$ is compact.}

\begin{center}
\bf\Large
Groups of isometries of ultrametric Urysohn spaces and  
their unitary representations 

\medskip

\sc \large
Yu. A. Neretin%
\footnote{Supported by the grant FWF, Project P31591.}
\end{center}

{\small
We consider groups $\G$ of isometries  of   ultrametric Urysohn spaces $\U$. Such spaces $\U$ admit transparent
realizations as boundaries of certain $\R$-trees and the groups $\G$ are groups of automorphisms of these $\R$-trees. Denote by $\G[X]\subset \G$ stabilizers of finite subspaces
$X\subset \U$. Double cosets $\G[X]\cdot g\cdot \G[Y]$, where $g\in \G$, are enumerated by  ultrametrics on union of spaces $X\cup Y$.  We construct natural associative multiplications on double coset 
spaces $\G[X]\backslash \G/\G[X]$ and, more generally, multiplications
$\G[X]\backslash \G/\G[Y]\,\times\, \G[Y]\backslash \G/\G[Z]\to \G[X]\backslash \G/\G[Z]$.
These operations are a kind of  canonical amalgamations of ultrametric spaces. On the other hand, 
this product can be interpreted in terms of partial isomorphisms of certain $\R$-trees (in particular, we come to an inverse category).
 This allows us
to classify all unitary representations of the groups $\G$ and to prove that groups $\G$ have type $I$.
We also describe a universal semigroup compactification of $\G$ whose image in any unitary representation of $\G$ is compact.

}

\section{Introduction}

\COUNTERS

{\bf\punct Ultrametric spaces. Notation.%
\label{ss:notation}} Recall that an {\it ultrametric space} $X$ is a metric space, where a metric satisfies the inequality
$$
d(x,z)\le \max \bigl(d(x,y),\, d(y,z)\bigr).
$$

For ultrametric spaces we introduce the following notation:

\sm

--- $\Isolated(X)$ is the set of isolated points of $X$;

\sm

--- $\Isom(X)$ is the group of all isometries of a space $X$;
we use the term '{\it isometry of metric spaces}' for bijections preserving distances and
'{\it isometric embedding}' for injective maps preserving distances;

\sm

---   $\B^c(a,r)$ is
 a {\it closed ball} with center $a$ and radius
$r$, i.e., the set of all $x\in X$ such that $d(a,x)\le r$;

\sm

---  $\B^o(a,r)$ is the {\it open ball}, i.e.,
the set of $x$ satisfying 
$d(a,x)<r$;

\sm

--- a {\it ball} $\B^\epsilon(a,r)$, where $\epsilon=$'$c$' or '$o$',
is a closed or open ball;

\sm

--- $\S(a,r)$ is the sphere of radius $r$ with center $a$deal;

\sm

---   $\Lambda$ is the set 
 of all possible distances $d(x,y)>0$, we call  $\Lambda$
  by {\it spectrum} of $X$;

\sm

--- $\Lambda^\circ$ is the set of all $r$,
for which there exists a strictly increasing sequence $\lambda_k\in \Lambda$
convergent to $r$; 

\sm

---   $\Lambda^{+0}$  is $\Lambda$ if $0$ is a limit point of $\Lambda$, and $\Lambda^{+0}=\Lambda\cup 0$
otherwise.

\sm

Notice that both $\B^c(a,r)$, $\B^o(a,r)$
are {\it clopen sets} (closed  and open  simultaneously).
However, we use the terms 'open ball' and 'closed ball'.

  Clearly, for any point $b\in \B^c(a,r)$ we have $\B^c(b,r)=\B^c(a,r)$.
  The same holds for balls $\B^o(a,r)$. Two balls $\B_1$, $\B_2$ are disjoint, or
coincide,
or one ball contains another.

For two disjoint balls $\B_1$, $\B_2$ we define a distance
$$\delta(\B_1,\B_2):=d(x,y), \qquad \text{where $x\in \B_1$, $y\in \B_2$,}$$
the right hand side does not depend on the choice of $x$, $y$.    Consider a set
$\{\B_j\}$, whose elements are  
 several
 mutually disjoint balls. Then the metric  $\delta(\cdot,\cdot)$ defines a structure of  
 an ultrametric space on the set $\{\B_j\}$ . 

\sm 

{\it Perfect balls.}
It is known (see our Subsect. \ref{ss:trees-of})
  that for ultrametric spaces it is useful
to consider a ball as an abstract ideal object determined by 
  the corresponding subset in $X$ and radius.
We say that a ball $\B^\epsilon(a,r)$ is {\it perfect} if it does not coincide
(as a set)
with a ball $\B^{\epsilon'}(a,r')$ that looks as smaller.
For instance, a closed ball $\B^c(a,r)$
is perfect if the sphere $\S(a,r)$ is non-empty,
or, equivalently, $\B^c(a,r)
\supsetneqq
\B^o(a,r)$. An open ball $\B^o(a,r)$ is perfect
if for arbitrary $\delta>0$ we have a proper inclusion of sets
$\B^o(a,r)\supsetneqq \B^c(a,r-\delta)$.
A radius of a perfect open ball is contained in the set $\Lambda^\circ$
defined above.

 \sm

{\it We consider only  separable ultrametric spaces.}
 For  such spaces  spectra $\Lambda$
are at most countable.
For each $r>0$ the set of distinct balls  of radius $r$ is at most countable.

\sm

{\bf \punct Urysohn ultrametric  spaces.%
\label{ss:urysohn1}} Fix a countable or finite subset $\Lambda\subset (0,\infty)$.
The {\it Urysohn ultrametric space} $\U_\Lambda$ is a unique up to an isometry complete separable ultrametric space 
$\U_\Lambda$
with spectrum $\Lambda$ satisfying the following equivalent conditions:

\sm

{\it The Urysohn  property}. Let $X\subset Y$ be finite ultrametric spaces whose spectra are contained in $\Lambda$. Then each isometric embedding $X\to \U_\Lambda$ can be extended to an isometric embedding $Y\to \U_\Lambda$. See
Gao, Shao \cite{GS}, Definition 2.3.

\sm 

{\it  The Ismagilov property}. For each positive $\lambda\in \Lambda$ and each point $x$ there is a countable collection
of  points $\{x_j\}$ including $x$ such that $d(x_i,x_j)=\lambda$ if $i\ne j$. See Ismagilov \cite{Ism97}, Lemma 1.

\sm

The space $\U_\Lambda$ also satisfies the following properties:

\sm

{\it The universality}. Any separable ultrametric space whose spectrum is contained in $\Lambda$ admits an isometric embedding to $\U_\Lambda$. See Gao, Shao \cite{GS}, Proposition 2.7.

\sm

{\it The ultrahomogeneity}. Let  $A$ be a finite ultrametric space and $\phi$, $\psi$ be embeddings $A\to \U_\Lambda$. Then there exists an isometry $\theta$ of  $\U_\Lambda$ such that $\psi=\theta\circ \phi$.
See \cite{GS}, Proposition 2.7.

\sm

Clearly, universality $+$ ultrahomogeneity imply the Urysohn property.

See the survey by Gao, Shao \cite{GS} containing several constructions of spaces $\U_\Lambda$.

\sm

{\sc Remarks.}
1)
Recall that  the original {\it Urysohn space}%
\footnote{A posthumous Urysohn's work \cite{Ury} was published in 1927, it  
was prepared for a publication by P.~S.~Alexandroff. 
In 1951 the paper was republished in collected works of Urysohn. 
I do not any continuation of this work until
a publication of  Kat\v{e}tov \cite{Kat}, 1988. Between these dates there were discovered several
relatives of the Urysohn space as the Rado graph, the universal poset, the Gurarii space.}
 is the complete separable metric space $\U$ satisfying the following condition:
 
 \sm
 
 {\it The Urysohn   property}:  Let $X\subset Y$ be finite metric spaces.
  Then each isometric embedding $X\to \U$ can be extended to an isometric embedding $Y\to \U$.
  
  \sm
  
   The Urysohn space is unique and satisfies
 a collection of amazing properties, see, e.g., Uspenskii \cite{Usp}, Bogatyi \cite{Bog-U}, Vershik \cite{Ver}, Melleray \cite{Mel}.

\sm

2)
A universal ultrametric space was constructed by Vestfrid \cite{Ves}. He considered universality for all possible distances $\in\R$ and his construction gives a nonseparable space. However, his approach allows to produce spaces $\U_\Lambda$, see \cite{GS}.  Ismagilov \cite{Ism97} considered non-Archimedean non-locally compact
 normed fields as ultrametric Urysohn spaces. A general construction of spaces $\U_\Lambda$
 is contained in Nguyen \cite{Ngu}, see also \cite{GS}. 
 \hfill $\boxtimes$

\sm

{\bf \punct Constructions of ultrametric Urysohn spaces.%
\label{ss:urysohn-construction}} 

\sm 

{\it  A special case}, see Ismagilov \cite{Ism97}.
Let $\Sigma$ be a countable  subgroup in the additive group of $\R$
(for instance, $\Sigma=\Q$).
Let $\Bbbk$ be a countable field. 
Let $t$ be a formal variable. Consider the space $\K$ consisting of formal 
 Puiseux-type series 
\begin{equation}
\label{eq:puizo}
f(t)=\sum_j a_j t^{s_j}, \qquad \text{where $a_j\in \Bbbk$, $s_j\in\Sigma$, $s_1<s_2<\dots$, $s_j\to +\infty$.}
\end{equation}
Clearly, $\K$ is a field with respect to the usual operations over formal series. Fix $h$ such that $0<h<1$ and define a norm
in $\K$ by 
\begin{equation}
\label{eq:norm}
\|f\|:=h^{\min\bigl( \text{$s_j$ such that $a_j\ne 0$}\bigr)}.
\end{equation}
Then $d(f_1,f_2)=\|f_1-f_2\|$ is an ultrametric on $\K$, the spectrum $\Lambda$ of $\K$
consists of numbers $h^s$, where $s\in\Sigma$, and $\K$ is an Urysohn space $\U_\Lambda$.

\sm

{\it A general construction. Variant} 1.
  Let $\Sigma\subset \R$ be arbitrary countable of finite set. 
Choose $h$ as above. 
We fix a countable Abelian group  $\mathbbb{z}$  
and consider the space $\cZ$ of formal series \eqref{eq:puizo} with coefficients in $\mathbbb{z}$.
In this case there is no structure of a field or a ring on $\cZ$. However, we 
have a norm \eqref{eq:norm} and the corresponding distance. Again, we get
the Urysohn space $\U_\Lambda$ whose spectrum $\Lambda$ consists of numbers $h^s$, where $s\in\Sigma$.

\sm 

{\it A general construction. Variant} 2.
 It is more reasonable to reformulate the previous construction in the following form  
 (apparently this was originally done by Nguyen \cite{Ngu}).
 Fix a countable set $Z$ with a distinguished element $z^*$.
The space $\U_\Lambda$ is the space of all functions $\omegA:\Lambda\to Z$
such that $\omegA(\lambda)=z^*$ for all $\lambda$ except a finite set or a sequence convergent to 0.
The distance between $\omegA_1$ and $\omegA_2$ is the maximal $\lambda$ for which $\omegA_1(\lambda)\ne \omegA_2(\lambda)$. 

Below we set $Z=\Z$, $z^*=0$.

\sm

{\bf\punct The groups $\boldsymbol{\G_\Lambda}$.%
\label{ss:isometries}}
Denote by $\G_\Lambda$ the group of all isometries of an Urysohn space $\U_\Lambda$.

 Recall that a topological group is called {\it Polish} (see, e.g., \cite{Kech}) if 
its underlying topological space is Polish, i.e., is homeomorphic to a complete separable metric space.
Consider the group $\Isom(A)$ of isometries of a complete separable metric space $A$. We say that a sequence $g_j\in \Isom(A)$
converges to $g$,  if for each $a\in A$ the sequence
$g_j(a)$ converges to $g(a)$. This convergence determines a Polish topology on 
$\Isom(A)$, see, e.g., \cite{Kech}, I.9.B. So the groups $\G_\Lambda$ are Polish.
 
 For a finite collection $\{\B^c(a_1, r_1)$, \dots, $\B^c(a_n, r_n)\}$ of disjoint balls
 in $\U_\Lambda$
 consider the subgroup $\G\{\B^c(a_j,r_j)\} \subset \G_\Lambda$
 consisting of isometries sending these balls to themselves. Such subgroups are open and form
 a fundamental system of neighbourhoods of  unit in $\G_\Lambda$  (this is a reformulation
 of pointwise convergence).
 
 \sm
 
 {\sc Remark.} The statement: '{\it Any homomorphism between Polish groups is continuous}'
 is compatible with the Zermelo--Fraenkel system plus the Axiom of dependent choice, see, e.g., \cite{Wri}, see also
 the Pettis theorem \cite{Kech}, Theorem 9.9.
 Informally, this implies that any group (including $\G_\Lambda$) admits at most one explicit Polish topology. \hfill $\boxtimes$

\sm  

{\bf \punct Amalgams of ultrametric spaces.%
\label{ss:amalgama1}} 
The following statement was obtained by Bogatyi \cite{Bog}, proof of Theorem 2.2.

\begin{theorem}
\label{th:amalgama}
Let $X$, $Y$ be sets, $Z=X\cap Y$. Let $d_X$ and $d_Y$ be ultrametrics on $X$ and $Y$
respectively, let $d_X=d_Y$ on $Z$. For $x\in X\setminus Z$, $y\in Y\setminus Z$
 we define a distance by
 $$
 d(x,y):=\min_{z\in Z} \max \bigl(d_X(x,z),d_Y(z,y)\bigr).
 $$
Then we get an ultrametric on $X\cup Y$.
\end{theorem}

In Subsect. \ref{ss:amalgamas2} we present a simple proof of this theorem.

We call this construction by {\it amalgamation} of ultrametric spaces.
For our purposes it is sufficient to consider finite ultrametric spaces.

\sm 

We define the following  {\it category} $\cK_\Lambda$. {\it Objects} are finite 
ultrametric spaces $X$ whose spectra are contained in $\Lambda$. {\it Morphisms}
 $\frp:X\to Y$ are
triples $\frp=(P, p_+,p_-)$, where

\sm

--- $P$ is an   ultrametric space;

\sm

--- $p_+:X\to P$, $p_-:Y\to P$ are isometric embeddings and
$$p_+(X)\cup p_-(Y)=P.$$

\sm 

Two morphisms $\frp=(P, p_+,p_-)$, $\frp'=(P', p_+',p_-')$ are {\it equal}
if there is an isometry
 $\pi:P\to P'$ such that $p_\pm'=\pi\circ p_\pm$.

\sm

For two morphisms $\frp: X\to Y$, $\frq:Y\to Z$ we define their {\it product} 
$\frr=\frq\diamond\frp:X\to Z$
in the following way (see Fig.\ref{fig:1}.a). We consider the disjoint union $P\sqcup Q$ and for each $y\in Y$ we identify
points $p_-(y)\in P$ and $q_+(y)\in Q$. Denote the resulting set by $P\Join Q$.
 By Theorem \ref{th:amalgama} we have a canonical ultrametric on $P\Join Q$. By construction,
 we have maps $p_+:X\to P\Join Q$, $q_-:Z\to P\Join Q$. We set
 $R:=p_+(X)\cup q_-(Z)\subset P\Join Q$, and $r_+=p_+$, $r_-=q_-$. 
 
 \begin{figure}
$${\mathrm a)}\qquad \epsfbox{kvadraty.1}$$

$${\mathrm a)}\qquad \epsfbox{kvadraty.2}$$
\caption{a) To Subsects. \ref{ss:amalgama1}, \ref{ss:PB-variant}. To the definition of the product in the category $\cK$.
\newline
b) To Subsect. \ref{ss:PB-variant}. 
 } 
 \label{fig:1}
\end{figure} 
 
 \begin{lemma}
 \label{l:ass}
 This product is associative, i.e., for any finite ultrametric spaces $X$, $Y$, $X$, whose spectra are contained in $\Lambda$, and  for any morphisms $\frp:X\to Y$, $\frq:Y\to Z$, $\frr:Z\to U$ we have
 $$
 (\frp\diamond \frq)\diamond \frr=  \frp\diamond (\frq\diamond \frr).
 $$ 
 \end{lemma}
 
 In Section \ref{s:amalgams} we present an equivalent description of our category
 in terms of partial isomorphisms of trees. This reformulation automatically implies the associativity.

\sm

We  have an {\it involution} $\frp\mapsto \frp^*$ on this category.
Namely, we define $\frp^*:Y\to X$ leaving $P$ the same and renaming $p_+$ to
$(p^*)_-$ and $p_-$ to $(p^*)_+$. Clearly,
$$(\frq\diamond \frp)^*=\frp^* \diamond \frq^*.$$

\sm

For each $X$ we have a {\it unit endomorphism} $1_X$, namely $P=X$, $p_\pm:X\to X$ are
identical maps.
 
\sm

We need in some standard notation related to categories.
For arbitrary category $\cL$ we denote by 

\sm 

--- $\Ob(\cL)$
the set (class) of objects of $\cL$;

\sm

---  ${\Mor}_\cL(X,Y)$ the set of all morphisms   $X\to Y$;

\sm

--- ${\End}_\cL(X)$ the semigroups of morphisms $X\to X$;

\sm

--- ${\Aut}_\cL(X)$ the group of automorphisms of $X$.

\sm

The topic of our interest is the category $\cK_\Lambda$,
in this case  we usually  omit a subscript $\cK_\Lambda$. 

Clearly, 
$\Aut_{\cK_\Lambda}(X)\simeq\Isom(X).$

 \sm 
 
{\bf \punct Reduction of unitary representations%
\footnote{The  subject of the present paper are spaces $\U_\Lambda$ and the groups $\G_\Lambda$. They are far
from traditional topics of representations theory, it seems that specialists in Polish groups, general topology, and descriptive set theory
are more familiar with objects of such kind than representation theorists. For this reason, I am trying to minimize a representation
theoretic background necessary for understanding of the paper. For a minimal collection of definitions, see, e.g., \cite{BHV}, Appendix A.} of groups $\boldsymbol{\G_\Lambda}$
 to representations
of the category $\boldsymbol{\cK_\Lambda}$.%
\label{ss:double-cosets}}
Let $X\subset \U_\Lambda$ be a finite subspace. Denote by
$\G[X]\subset \G_\Lambda$ the point-wise stabilizer of $X$.
Let $X$, $Y$ be finite subsets in $\U_\Lambda$, i.e., objects of
the category $\cK(\Lambda)$. For $g\in \G_\Lambda$
set $P:=Y\cup Xg$. 
By construction, we have an embedding $p_+:X\to P$, namely the map $g\Bigr|_X$,
and the identical embedding $p_-:Y\to P$.
So we have a morphism of the category  $\cK_\Lambda$. 
Clearly, for $h_1\in \G[X]$, $h_2\in \G[Y]$, isometries
$h_2gh_1$ and $g$ determine the same morphism.
So we have a bijection of the double coset space
with the set of morphisms:
\begin{equation}
\G[Y] \backslash \G_\Lambda/\G[X]\simeq \Mor(X,Y).
\label{eq:coset-identification}
\end{equation}

 Let $\rho$ be a unitary representation
of $\G_\Lambda$ in a Hilbert space%
\footnote{A {\it unitary representation} of a topological
group $G$ is a continuous homomorphism to the group $\mathrm{U}(H)$ of unitary operators
of a separable Hilbert space; the group $\mathrm{U}(H)$ is 
 equipped with the strong operator topology, see, e.g., \cite{Dix}, 13.1.1-13.1.2, \cite{BHV}, Appendix A. The strong operator topology on $\mathrm{U}(H)$
coincides with the weak operator topology.
The group $\mathrm{U}(H)$ is  Polish, see Remark at the end of Subsect. \ref{ss:isometries}.%
} 
 $H$.
Denote by $H[X]\subset H$ the subspace consisting of $\G[X]$-fixed vectors,
by $\Pi[X]$ the operator of orthogonal projection to $H[X]$.
For any $g\in \G_\Lambda$ we define the operator
$\wh\rho_{X,Y}(g):H[X]\to H[Y]$ 
by 
$$
\wh\rho_{X,Y}(g):=\Pi[Y]\, \rho(g)\Bigr|_{H[X]}.
$$
It is easy to see that $\wh\rho_{X,Y}(g)$ depends only on the double coset%
\footnote{Indeed, let $h_1\in \G[X]$, $h_2\in \G[Y]$, let $\xi\in H[X]$, $\upsilon\in H[Y]$.
Then
\newline
 $\la \wh\rho_{X,Y}(h_2 g h_1)\xi, \upsilon\ra_{H[Y]}=  \la \Pi[Y]\rho(h_2 g h_1)\xi, \upsilon\ra_{H}=\la \rho(h_2 g h_1)\xi,\Pi[Y] \upsilon\ra_{H}=
 \la \rho(h_2 g h_1)\xi,\upsilon\ra_{H} =$
 \newline  
 $=
 \la \rho(g)(\rho(h_1)\xi),\rho(h_2^{-1}) \upsilon\ra_{H}
=\la \rho(g)\xi, \upsilon\ra_{H}=\la\Pi[Y] \rho(g)\xi, \upsilon\ra_{H}=\la\wh \rho_{X,Y}(g)\xi, \upsilon\ra_{H[Y]}$.
\newline
So all matrix elements of the operators $\wh\rho_{X,Y}(h_2 g h_1)$ and $\wh\rho_{X,Y}(g)$ coincide.%
}
$$\frg:=\G[Y]\cdot g\cdot \G[X]\,\in \G[Y]\backslash \G_\Lambda/\G[X].$$
By \eqref{eq:coset-identification}, we can regard $\wh\rho_{X,Y}(g)$
as a function on $\Mor(X,Y)$.

\begin{theorem}
\label{th:multiplicativity}
{\rm a)}
The maps $\frp\mapsto \wh\rho_{X,Y}(\frp)$ determine a representation%
\footnote{For basic definitions of representations of categories, see, e.g., \cite{Ner-book}, Sect. 2.5, Sect. 2.8.
However, we avoid a usage of such definitions.} of
the category $\cK_\Lambda$, i.e.,
 a 
functor from the category
$\cK_\Lambda$ to the category of Hilbert spaces and  linear
operators.
This means that for any finite subsets $X$, $Y$, $Z\subset \U_\Lambda$
and any $\frp\in\Mor(X,Y)$, $\frq\in\Mor(Y, Z)$ we have
$$
\wh\rho_{Y,Z}(\frq)\, \wh\rho_{X,Y}(\frp)=\wh\rho_{X,Z}(\frq\diamond\frp).
$$

{\rm b)} The operators $\wh\rho_{X,Y}(\cdot)$ are contractive,
i.e.,
$$
\|\wh\rho_{X,Y}(\frp)\|\le 1.
$$

\sm

{\rm c)} The representation $\wh\rho$ is a {\rm $*$-representation 
of the category $\cK_\Lambda$},
i.e., for any finite subsets $X$, $Y\subset \U_\Lambda$ and any 
$\frp\in \Mor(X,Y)$ we have
$$
\wh\rho_{X,Y}(\frp)^*=\wh\rho_{Y,X}(\frp^*).
$$

{\rm d)} Each map $\wh \rho_{X,Y}$ is continuous as a map
from the quotient topological   space%
\footnote{These sets are countable, however they are not discrete,
see the description of the quotient topology in Lemma \ref{l:topology}.}
 $\G[Y] \backslash \G_\Lambda/\G[X]$
to the space of operators $H[Y]\to H[X]$ equipped with the weak operator
topology.

\sm

{\rm e)} The representation $\wh\rho$ is {\it non-degenerate} in the following sense%
\footnote{In our case this is equivalent to the following property: for any $X$ the closed subspace generated by all vectors
$\wh\rho (\frp)\xi$, where $\xi$ ranges in $H[X]$ and $\frp$ ranges in $\End(X)$, is the whole $H[X]$. In  this form, this
corresponds to the definition of non-degenerate representations of $C^*$-algebras, see, e.g., \cite{Dix}, 2.2.6. 
}: for all
$X$ we have
$$
\rho_{X,X}(1_X)=1_{H[X]}.
$$

\end{theorem}

The nontrivial statement is the first claim of the theorem, other  statements are obvious.

 This theorem (it is proved in Section \ref{s:multiplicativity}) reduces a classification of unitary representations of $\G_\Lambda$
 to an  easy problem (which is solved in Section \ref{s:classification}).
 
 \begin{theorem}
 The construction above provides a one-to-one correspondence
 between  unitary representation 
 of the group $\G_\Lambda$ and nondegenerate
 $*$-rep\-re\-sen\-ta\-tions $\tau$ of the category $\cK_\Lambda$. 
 \end{theorem}

 \sm
 
 {\bf \punct Classification of unitary representations 
 of $\boldsymbol{\G_\Lambda}$.%
 \label{ss:classification}}
We say that a pair $(\epsilon, r)$, where $\epsilon=$'$c$', '$o$'
and $r\in \R$
 is  $\Lambda$-perfect (or simply perfect), if balls $\B^\epsilon(x,r)$ are perfect
 in $\U_\Lambda$.
 Let us formulate the definition more carefully. $\Lambda$-{\it Perfect
 pairs} are pairs of the form:
 
  \sm
 
 1) $(c,r)$ with $r\in \Lambda^{+0}$;
 
 \sm
 
 2) $(o,r)$ with $r\in \Lambda^\circ$. If the set
 $\Lambda$ is unbounded, then we also allow the pair $(o,\infty)$,
 i.e., we consider the whole space as an open ball of infinite radius.
 
\sm

For a perfect pair $(\epsilon, r)$ denote by $\Xi_{\epsilon,r}$
the set of all balls $\B^\epsilon(x,r)$. Since the space $\U_\Lambda$
is separable, the space $\Xi_{\epsilon,r}$ is countable.
The group $\G_\Lambda$ acts on 
$\Xi_{\epsilon,r}$ and therefore it acts in the Hilbert
space%
\footnote{For a countable or finite set $\Omega$ we denote by $\ell^2(\Omega)$
the Hilbert space of functions $f$ on $\Omega$ satisfying the condition
$\|f\|^2:=\sum_{\omega\in\Omega} |f(\omega)|^2<\infty$.}
 $\ell^2(\Xi_{\epsilon,r})$.

A shortest but non-explicit form of classification theorem
is:

\newtheorem*{theorema1}{Theorem \ref{th:classification-main}$^\circ$}
\newtheorem*{theorema3}{Theorem \ref{th:classification-main}$^{\circ\circ}$}

\begin{theorema1}
Any unitary representation of the group $\G_\Lambda$
is a subrepresentation in some tensor product
$$
\bigotimes_{j=1}^N \ell^2(\Xi_{\epsilon_j,r_j}),
$$
where pairs $(\epsilon_j,r_j)$ are $\Lambda$-perfect 
{\rm (}and not necessary distinct{\rm)}.
\end{theorema1}

Our next purpose is to  formulate the classification theorem in a precise form.
Denote by $\Xi^\bullet=\Xi^\bullet(\Lambda)$ the space of all subsets 
$L\subset \U_\Lambda$
that can be represented as unions of finite collections 
of disjoint balls, $L=\cup_{j=1}^N \B^{\epsilon_j}(x_j,r_j)$. Denote by $\Xi(\Lambda)$ the set of all finite ordered collections
of disjoint ball $\{\B^{\epsilon_j}(x_j,r_j)\}_{j\le N}$. 

 Let us describe orbits of $\G_\Lambda$ on $\Xi^\bullet(\Lambda)$ and $\Xi(\Lambda)$.

 \sm
 
Fix a finite subset $X\subset \U_\Lambda$.
For each $x\in X$ denote
$$
d_x=\min_{y\in X} d(x,y),
$$
if $X$ consists of one point $x$, then we assume $d_x=\infty$.
Let us for each $x\in X$ fix a $\Lambda$-perfect pair  $(\epsilon_x, r_x)$,
where $r_x>0$ satisfies the conditions:
\begin{align*}
r_x<d_x\qquad \text{if $\epsilon_x=c$};\\
r_x\le d_x\qquad \text{if $\epsilon_x=o$}.
\end{align*} 
This  means that for each point $x\in X$ we draw  a ball $\B^{\epsilon_x}(x,r_x)$ around $x$ that does not contain points
$y\in X\setminus x$.
We call such data  $X$ by a {\it labeling} of the ultrametric space $X$.

Consider a
collection of balls 
\begin{equation}
\bigl\{\B^{\epsilon_x}(x, r_x)\}_{x\in X}\in \Xi.
\label{eq:collection}
\end{equation}
Denote by $\G\bigl\{\B^{\epsilon_x}(x, r_x)\}$
its stabilizer in $\G_\Lambda$ (by the definition, it consists of transformations sending each ball to itself).
Clearly, each $\G_\Lambda$-orbit on $\Xi$ has a representative of  
this form.
Two collections $\bigl\{\B^{\epsilon_x}(x, r_x)\}_{x\in X}$
and $\bigl\{\B^{\epsilon_{x'}}(x', r_{x'})\}_{x'\in X'}$
are $\G_\Lambda$-equivalent iff there is an isometry $\iota:X\to X'$
sending the labeling to the labeling, i.e.,
 $\epsilon_{\iota(x)}=\epsilon_{x'}$, 
$r_{\iota(x)}=r_{x'}$. 

Similarly, consider
\begin{equation}
\cup_{x\in X} \B^{\epsilon_x}(x, r_x)\in\Xi^\bullet.
\label{eq:union}
\end{equation}
and denote the stabilizer of this set in $\G_\Lambda$ by
$\G^\bullet\bigl\{\B^{\epsilon_x}(x, r_x)\}$.
Remarks of the previous paragraph remain valid in this case.

The group $\G^\bullet\bigl\{\B^{\epsilon_x}(x, r_x)\}/
\G\bigl\{\B^{\epsilon_x}(x, r_x)\}$ is finite. It is the subgroup
in $\Isom(X)$ consisting of transformations $X$  preserving the labeling
$\{\epsilon_x,r_x\}$.

\begin{theorem}
\label{th:classification-main}
{\rm a)}
Each irreducible representation of $\G_\Lambda$ 
 is induced
from an irreducible representation $\nu$ of some group $\G^\bullet\{\B_j\}$, which
is trivial on the subgroup $\G\{\B_j\}$.

\sm

{\rm b)} Two irreducible representations obtained in this way
 are equivalent if and only if the corresponding configurations of balls are $\G_\Lambda$-equivalent
  and  the corresponding representations of $\G^\bullet\{\B_j\}/\G\{\B_j\}$ are equivalent.
 
 \sm
 
 {\rm c)} The groups $\G_\Lambda$ have type $I$.
 
 \sm
 
 {\rm d)} Any unitary representation of $\G_\Lambda$ is a direct integral
 of irreducible representations. 
\end{theorem}

A definition of induced representations appropriate for our case, 
is given below in Subsect. \ref{ss:induced}.

\sm 

We also formulate the third version of the classification theorem:

 \begin{theorema3}
  \label{th:classification3}
 Each irreducible representation of $\G_\Lambda$ is a subrepresentation
 in a quasiregular representation in $\ell^2$
 on some homogeneous space
  $\G\bigl\{\B_j\}\bigl\backslash\G_\Lambda$.
 \end{theorema3}

Theorems \ref{th:classification-main}$^\circ$,  \ref{th:classification-main},  \ref{th:classification-main}$^{\circ\circ}$ are equivalent,
this equivalence can be established by 
  Mackey's machinery of induced representations, see \cite{Mack0}, \cite{Mack}, Chapter 3. In fact, we prove the statement in the strong
  form (Theorem \ref{th:classification-main}),
   two other statements are easy corollaries.
  
  \sm

  {\bf \punct Closures of $\boldsymbol{\G_\Lambda}$ in unitary representations.%
  \label{ss:compactifications}}
  Let $H$ be a Hilbert space. Consider the set $\bfB(H)$
  of all operators in $H$ with norm $\le 1$. The weak operator topology
  restricted to this set is metrizable and the set $\bfB(H)$
  is compact with respect to the weak operator topology.
  On the other hand $\bfB(H)$ is a semigroup with respect to multiplication
  of operators, and this operation is separately continuous. See, e.g.,
  \cite{Ner-book}, Sect. 1.4.
  
  Let $G$ be a topological group and $\rho$ be its unitary representation
   in $H$. Denote by ${\rho(G)}$ the set of all operators $\rho(g)$
   and by   $\ov {\rho(G)}$ the closure of this set in the weak operator
   topology. It is easy to see that $\ov {\rho(G)}\subset \bfB(H)$ 
   is a compact semigroup with separately continuous product.  See, e.g.,
  \cite{Ner-book}, Sect. 1.4.
  
  There is also a question about universalization of such compactifications.
  Consider a family $\rho_\alpha$, where $\alpha$ ranges in some set $A$,  of unitary representations of $G$ in Hilbert spaces $H_\alpha$ (for instance, we
  can consider the set of all unitary representations or the set of all 
  irreducible representations, there exist other natural situations). Then we have 
  a family of embeddings
  $\rho_\alpha:G\to \bfB(H_\alpha)$ and therefore we have the diagonal embedding of
  $G$ to the direct product  $\prod_{\alpha\in A}\bfB(H_\alpha)$. By the Tychonoff theorem the direct product is compact and therefore  the closure of $G$
  is a compact semigroup. As usual, a usage of the Axiom of choice (the Tychonoff theorem depends on it)
   can produce non-tame
  objects (as Bohr compactification of $\R$, see, e.g., \cite{Dix}, 16.1.1), but such phenomena happen mainly for Abelian groups.
  
  Problems  of description of such semigroups for various  $G$ arose many times by many reasons.
  Cases of Abelian groups, Lie groups, infinite-dimensional groups looks as seriously different.
  
  Question about closures of powers of a measure preserving
  transformation (this is related to the  closure of the group $\Z$ acting in a Hilbert space)
   is highly non-trivial and was widely discussed in ergodic theory, see, e.g., \cite{King},
  \cite{Jan},
  \cite{Sol}. Melleray and Tsankov described structure of the closure of powers of  generic unitary operators, see \cite{MT}. 
  
   Howe and Moore \cite{HM} considered the case of Lie groups (over real and $p$-adic fields).
  For  irreducible unitary representations  of semisimple Lie groups with finite center we  always get the one-point compactification.
    So a compactification is a trivial object, but this is an important statement about behaviour
   of matrix elements of unitary representations. 
  For general Lie groups the picture is similar, but closures of centers can produce arbitrary compactifications of centers.

  Olshanski
  \cite{Olsh-kiado}, \cite{Olsh-semigroups} noticed that for
  infinite-dimensional groups $G$ compactifications must be seriously different
  from the group themselves, and (this is more important), that 
  such compactifications admit constructive descriptions. In fact, before this, boundary elements of groups appeared 
  as a tool in several representation-theoretical works.
  See also \cite{Ner-book}, Sect. 1.4, Notes to Sect. 4.2, Sect. 8.1-8.5.
  
   Ben Yakov, Ibarluc\'\i a, Tsankov \cite{BIT} considered  the problem for 
  oligomorphic groups%
  \footnote{A totally disconnected group $G$ is {\it oligomorphic} if it can be realized
  as a closed subgroup of the group $S_\infty$ (see below Subsect. \ref{ss:symmetric})
   of all permutations of a countable set $\Omega$ in such a way that $G$ has a finite number of orbits
  on any finite product $\Omega\times\dots\times\Omega$. According Tsankov \cite{Tsa}, any irreducible representation of an oligomorphic
  group is induced from an open subgroup. The same conclusion holds for inverse limits of oligomorphic groups.
   For our groups $\G_\Lambda$ all irreducible representations are induced from finite-dimensional representations of open subgroups.
If $\Lambda$ is finite, then $\G_\Lambda$ is oligomorphic, in this case it is a multi-stage wreath product of infinite symmetric groups,
$(((S_\infty\wr \dots\wr S_\infty) \wr S_\infty)\wr S_\infty)\infty$. If $\Lambda$ is a sequence convergent to 0, then the group $G_\Lambda$
is a countable-stage wreath product and is an inverse limit of oligomorphic groups. Apparently, all other groups $\G_\Lambda$ are not
inverse limits of
oligomorphic groups, see Conjecture \ref{conj:1}.}
   and inverse limits of oligomorphic groups and get
  an abstract description of such compactifications in terms of partial bijections,
   which for many groups is explicit (for instance, for the group
   of automorphisms of the Rado graph).

   In Section \ref{s:university} we  construct the  universal compactifications (for all unitary representations) $\GammA_\Lambda$
   of  groups $\G_\Lambda$. The semigroups $\GammA_\Lambda$  are described in terms of some kind of partial isomorphisms
   of  $\R$-trees. The problem  a priori is
   related to products of double cosets, since projectors $\Pi[X]$ are contained in 
    $\wt{\rho(\G_\Lambda)}$, see Lemma \ref{l:weakly}. Therefore all operators
   $
   \Pi[Y]\rho(g)\Pi[X]
   $
   are contained in $\ov{\rho(\G_\Lambda)}$. So the whole category $\cK_\Lambda$ can be regarded as a part
   of the semigroup $\GammA_\Lambda$.
   
   Our construction has obvious parallels with \cite{BIT}, on the other
   hand partial isomorphisms of trees arise in the context of \cite{Olsh-new},
   see comments in \cite{Ner-book}, Sect. 8.5, Notes.
   
   \sm

 {\bf \punct General comments on multiplication of double cosets.}
 Multiplications of double cosets and  multiplicativity theorems 
 (as Theorem \ref{th:multiplicativity}) are a usual phenomenon for 
 infinite-dimensional%
 \footnote{We use the term 'infinite-dimensional groups' as commonly accepted.
  Formally, dimensions of topological spaces $\G_\Lambda$ are 0.}
   groups. They arise for various types
 of  groups as infinite analogues of symmetric groups,
 see, e.g., Olshanski \cite{Olsh-sym}, \cite{Olsh-semigroups},
 Neretin
 \cite{Ner-field}; for infinite-dimensional real classical groups, see Olshanski
 \cite{Olsh-GB}, \cite{Olsh-semigroups}, Neretin \cite{Ner-book}, Sect. 9.3-9.4.
 \cite{Ner-coll}, \cite{Ner-faa}; for groups of measurable transformations
 \cite{Ner-bist}, \cite{Ner-matching}, \cite{Ner-book}, 10.3-10.4. There are constructions in this spirit
 for infinite-dimensional groups over finite fields \cite{Olsh-semigroups},
 \cite{Ner-finite} and infinite-dimensional $p$-adic groups (but these cases are less understood). There are
 also some 'exotic' constructions related to groups acting on trees, see, e.g., \cite{Olsh-new},
  \cite{Ner} (and the present paper); some constructions can be extracted from Tsankov \cite{Tsa}.
   Quite often, multiplications of double cosets produce unexpected algebraic
  structures, which looks drastically non-similar to the group themselves.
  
   Initially, such multiplications were introduced as a
  tool for classification of unitary representations. See various classification results
  in  \cite{Olsh-kiado},
  \cite{Olsh-new}, \cite{Olsh-rank}, \cite{Olsh-sym}, \cite{Ner-bist},
Subsect. 3.2, 3.4
  \cite{Ner-book}, Sect. 8.1-8.4, \cite{Ner-field}, Sect. 8,
   \cite{Ner-finite}.
 Apparently, the first construction of this type was discovered by Ismagilov
 \cite{Ism67} in 1967. He considered double cosets $\SL_2(\bbO)\backslash\SL_2(\K)/\SL_2(\bbO)$, where $\K$ is a 
  field of formal Laurent series over a countable field $\Bbbk$ and $\bbO\subset \K$ is the ring of integer elements,
  i.e., formal power series;
  cf. Subsect. \ref{ss:urysohn-construction} above. 
  
  Notice that for all cases considered in works by references double cosets can be interpreted as elements of the closure
  of the group (as in the previous subsection). On the other hand, there are infinite-dimensional
  groups having interesting representation theory (the group of diffeomorphisms of circle and loop groups),
   for which double coset machinery does not work,
  however the problem of closure in representations exists, see \cite{Ner-book}, Sect.9.4.
 
 \sm
 
 {\bf\punct Further structure of the paper.}
 Section \ref{s:trees} contains preliminaries on $\R$-trees corresponding to ultrametric spaces and on embeddings
 of $\R$-trees and ultrametric spaces to Hilbert spaces. 
 
 In Section \ref{s:amalgams}, we interpret amalgamate  product 
 as a product of partial isomorphisms of trees.

  Section \ref{s:urysohn} contains
 initial facts on the Urysohn ultrametric spaces $\U_\Lambda$.
 
  In Section \ref{s:multiplicativity},
 we reduce representations of the groups $\G_\Lambda$
  to $*$-representations of the category $\cK_\Lambda$ (Theorem \ref{th:multiplicativity}). The main tools are Lemma \ref{l:weakly}
  and Proposition \ref{pr:density}.
  The section also contains preliminaries on the infinite symmetric groups
  and preliminaries on induced representations of totally disconnected groups.
  
   In Section \ref{s:classification}, we 
  classify $*$-representations of the category $\cK_\Lambda$. In fact, the problem a priory is reduced to classifications
  of representations of semigroups $\End_{\cK_\Lambda}(X)$. They
are inverse semigroups (in the sense of V.~V.~Wagner and G.~B.~Preston),  this simplifies the situation.  
  
In Section \ref{s:university}, we describe the universal compactifications
$\GammA_\Lambda$ of the groups $\G_\Lambda$.
 
  In Section \ref{s:affine}, we discuss affine isometric actions of the groups $\G_\Lambda$
  on Hilbert spaces and  show that a group $\G_\Lambda$ has Kazhdan
  property (T) if and only if $\Lambda$ has a maximal element.  
  
  An $\R$-tree related to a  non-Archimedean normed field $\K$ described in Subsect. \ref{ss:urysohn-construction},
  admits a group of symmetries, which is larger than the group of  isometries of the ultrametric space $\K$. Apparently, this
  larger group also
  admits a classification of unitary representations (and it is more interesting, since it has irreducible representations
   that are not induced from open subgroups).
  In Section \ref{s:extension} we formulate a conjecture about its representations.
  
 \sm
 
 {\bf\punct Notation.} Here we list some notation that are introduced below:
 
 \sm
 
 --- $\ver(J)$ is the set of vertices of a tree or an $\R$-tree, see Subsect. \ref{ss:R-trees};
 
 \sm
 
 --- $\val(v)$ is the valence of a vertex of a tree or an $\R$-tree, see Subsect. \ref{ss:R-trees};
 
 \sm
 
 --- $\germ(v)$ is the set of germs of edges in a vertex $v$,  see Subsect. \ref{ss:R-trees};
  $\germ^\downarrow(v)$ denotes the set of germs looking down, See Subsect. \ref{ss:trees-of};
 
 \sm
 
 --- $[\xi,\eta]$, where $\xi$, $\eta$ are point of an $\R$-tree, is a segment, see Subsect. \ref{ss:R-trees};
 
 \sm
 
 --- $[[\xi_1,\dots,\xi_n]]$ is a convex hull in an $\R$-tree, see Subsect. \ref{ss:R-trees};
 
 \sm
 
 --- $T(X)$ is the $\R$-tree corresponding to an ultrametric space $X$, see Subsect. \ref{ss:trees-of};
 $\ov{T(X)}$ is the completion of $T(X)$, see Subsect. \ref{ss:completion}; 
  $\wt T(X)$, see Subsect. \ref{ss:woolly};
 
 \sm
 
 --- $\frh(\xi)$ is the height function on an $\R$-tree of an ultrametric space, see Subsect. \ref{ss:trees-of}
 
 \sm
 
 --- $L_x$ are 'ascending paths' in $T(X)$, see Subsect. \ref{ss:trees-of};
 
 \sm
 
 --- $\Ind_H^G(\nu)$ is an induced representation, see Subsect. \ref{ss:induced};
 
 \sm
 
 --- $\ell^2(\Omega,V)$ is the space $\ell^2$ consisting of functions on a countable set $\Omega$
 with values in a Hilbert space $V$, see Subsect. \ref{ss:induced}.
 
 \sm
 
 --- $S_\infty$ is the group of all permutations of $\N$; $S_\infty(\Omega)$ is the group of all permutations
 of a countable set $\Omega$, see Subsect. \ref{ss:symmetric};
 
 \sm
 
 --- $G\ltimes N$ is a semidirect product of groups, $N$ is a normal subgroup, $G$ is the quotient group.

\section{R-trees and ultrametric spaces. Preliminaries\label{s:trees}}

In this section, we discuss a standard construction
of  $\R$-trees of  ultrametric spaces
(Subsect. \ref{ss:R-trees}--\ref{ss:trees-of}, only these topics  are 
used in the proof of the classification theorem in Sect. \ref{s:amalgams}, \ref{s:multiplicativity}, \ref{s:classification}), embeddings
of $\R$-trees and ultrametric spaces to Hilbert spaces
(Subsect. \ref{ss:tree-hilbert}--\ref{ss:ultra-hilbert}, they are used in Sect. \ref{s:affine} and \ref{s:extension}),
 structure of completions of $\R$-trees 
of ultrametric spaces (Subsect. \ref{ss:given}--\ref{ss:completion}, they are used in discussion of homomorphisms
between different groups $\G_\Lambda$ in Sect. \ref{s:urysohn}), and an extension of the notion of 'subtree'
for $\R$-trees ('woolly subtrees', Subsect. \ref{ss:woolly}, it is used in the construction of the universal compactifications
of $\G_\Lambda$, Sect. \ref{s:university}).


\COUNTERS

{\bf\punct $\R$-trees.%
\label{ss:R-trees}}  See, e.g., Tits \cite{Tit},
Shalen \cite{Sha}, Chiswell \cite{Chi}, Bestvina \cite{Bes}.
Consider a metric space $T$.
An {\it arc} connecting two point $x$, $y$ is an image of a topological embedding
of a segment $[0,a]$ to $T$. 
 A {\it geodesic segment} in $T$ is an image of an
isometric embedding of a segment $[0,a]$ to $T$. 
A metric space in an {\it $\R$-tree} (or, a {\it real tree}) if for any two points $x$, $y\in T$ there is a unique
arc connecting $x$ with $y$ and this arc is a geodesic segment,
below we denote this arc by $[x,y]$.

There are many other equivalent definitions of $\R$-trees, we prefer the following constructive
way, which is valid for separable $\R$-trees.

Recall that a {\it tree} $J$ is a connected graph without cycles.
By $\ver (J)$ we denote the set of vertices of $J$; we consider only trees with countable or finite number
of vertices.
 We regard edges of a tree as segments.
Assigning a positive length to each edge we get a metric space, such spaces are called {\it metric trees}.

\sm

Let $J_1$, $J_2$, \dots be metric trees. Consider a chain of isometric embeddings
$$
J_1\stackrel{\iota_1}{\longrightarrow} J_2 \stackrel{\iota_2}{\longrightarrow} J_3 \stackrel{\iota_3}{\longrightarrow} \dots
$$
such that $\iota_k\bigl(\ver(J_k)\bigr)\subset \ver(J_{k+1})$ and consider its inductive limit $J_\infty$ equipped with the natural metric. This means that we have an increasing family of trees
$J_1\subset J_2\subset \dots$
 and  take its union. A point of $J_\infty$ is a {\it vertex}  if it is a vertex of some $J_m$
 and therefore for all $J_n$ for $n>m$. Clearly, $J_\infty$ is a separable metric space.
 Generally, a metric space $J_\infty$ is not complete. 
 
 Its completion $\ov J_\infty$ also
 is an $\R$-tree, see \cite{MNO}, Th.~1.11. It is easy to show that any separable $\R$-tree $T$
  is an intermediate space, $J_\infty\subset T\subset \ov J_\infty$ for some inductive limit
 $J_\infty$.
 
 \sm
 
 Fix $v$ be a point of  $J_\infty$. Then we have the following variants for $J_\infty\setminus v$:

 \sm
 
 1) if $v\notin \ver(J_\infty)$, then $J_\infty\setminus v$ splits into two connected components.

\sm

 2) if $v\in \ver(J_\infty)$, then $J_\infty\setminus v$ can consist of $\val(v)=1$, 2, \dots, $\infty$ connected
 components $C_1$, $C_2$, \dots 

\sm 

In the second case,  we say that $\val(v)$ is the {\it valence} of the vertex $v$. It is a limit of 
valences of $v$ in $J_m$ as $m$ tend to $\infty$. If $\val(v)=2$, then we say, that a vertex $v$ is
 {\it artificial},
it is not a vertex from the point of view of intrinsic geometry of an $\R$-tree. However for us is important
 to allow such vertices.

\begin{figure}
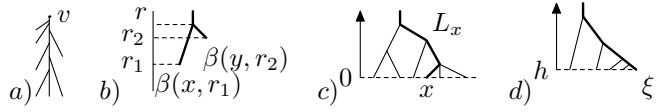

$$a)\epsfbox{trees.1}\quad b) \epsfbox{trees.2} \quad c) \epsfbox{trees.3} \quad d) \epsfbox{trees.4}
$$
\caption{a) To Subsect. \ref{ss:R-trees}. A germ of edge at a vertex $v$.
\newline
b) To Subsect. \ref{ss:trees-of}. A distance between two ideal balls.
\newline
c) To Subsect. \ref{ss:trees-of}. An ascending path.
\newline
d) To proof of Proposition \ref{pr:completion-tree}.
}
\label{fig:2}
\end{figure}

 Let $\psi_1$, $\psi_2$ be geodesic segments
$\psi_1:[0,a]\to  C_k\cup v$, $\psi_2:[0,b]\to  C_k\cup v$ such that $\psi_{1,2}(0)=v$. Then $\psi_1$, $\psi_2$ coincide on some domain 
$[0,\epsilon]$, where
$\epsilon>0$ depends on $\psi_1$, $\psi_2$. Generally, infimum of all possible $\epsilon$
for a given vertex $v$ is 0. So we have  {\it germs%
\footnote{Consider some family $\cF$ of maps $\psi_\alpha$ of segments $[0,a]$ to some space $T$. We say that $\psi_\alpha$, $\psi_\beta$
are equivalent if they coincide on some subsegment $[0,\epsilon]$, where $\epsilon>0$. A {\it germ} at 0 is a class of equivalent elements of $\cF$.
For separable $\R$-trees, the set of germs at a vertex is finite or countable.} of edges} at a vertex $v$, see Fig.\ref{fig:2}.a,
 but generally an $\R$-tree does not have edges. Denote the set of germs of edges at a vertex $v$ by $\germ(v)$. 

\sm

In the case $J_\infty \subset T \subset \ov J_\infty$ and $v\in T\setminus J_\infty$ the space
$T\setminus v$ is connected. 

\sm

{\sc Remark.} Points of the completion $\ov J_\infty$ arise in the following way. Consider a subsequence 
$J_{k_n}$ of our sequence of trees and a sequence of vertices $v_n\in \ver(J_{k_n})$ such that for each geodesic segment 
$[v_n,v_{n+1}]$ we have $[v_n,v_{n+1}]\setminus v_n\subset J_{k_{n+1}}\setminus J_{k_n}$.
Let $\sum d(v_n,v_{n+1})<\infty$. Then the sequence $v_n$ is fundamental and its limit is not 
contained in any $J_m$.

If vertices of valence $\ge 3$ are dense in $J_\infty$, then
the 'boundary'%
\footnote{We avoid a usage of the term 'boundary' in this context,  since 'boundary of
a tree'  has another meaning.} $\ov J_\infty\setminus J_\infty$ is dense in $\ov J_\infty$.
\hfill $\boxtimes$.

\sm

Consider a  collection of points $\{z_\alpha\}$ in an $\R$-tree $T$. Then 
$\cup_{\beta,\gamma} [z_\beta,z_\gamma]\subset T$ is an $\R$-tree, we call it 
the {\it hull} of $\{z_\alpha\}$  and denote by $[[z_1, z_2,\dots]]$.

\sm 

{\sc Remark.} There are several reasons for appearance of $\R$-trees in mathematics, in particular the following:

\sm

a) Consider the space of all actions of a discrete group $\Gamma$ on a hyperbolic space $L$
(in particular, on Lobachevsky plane).
It admits a natural compactification by actions of $\Gamma$ on $\R$-trees, see, e.g., \cite{Bes}, this idea arises to the construction
of the Thurston boundary of the Teichm\"uller space. Also the space of actions
of a given discrete group $\Gamma$ on simplicial trees admits a natural compactification by actions of $\Gamma$
on $\R$-trees.

\sm 

b) Recall that a {\it Bruhat-Tits tree} is a tree such that each vertex is contained in $(p+1)$ edges.
Groups $\SL(2)$ over $p$-adic numbers act in a natural way on such trees.
It is well-known that Brihat-Tits trees are right $p$-adic counterparts of rank 1 Riemannian symmetric spaces.
If we start with a non-Archimedean field with a nondiscrete valuation group, then we come to certain $\R$-trees,
see, e.g., \cite{Chi}. In particular, we come to such tree for the field $\K$ described in Subsect. \ref{ss:urysohn-construction}.

\sm 

c) $\R$-trees arise also in the context of branching stochastic processes.
\hfill $\boxtimes$

\sm

{\sc Remark.} The definition of $\R$-tree was proposed by Tits \cite{Tit} in 1977, in 1976 
Chiswell considered actions of groups on some $\R$-trees. It seems that the first example of an
$\R$-tree
appeared in one of  posthumous publications of Urysohn \cite{Ury2}, 1927, he constructed a nonseparable $\R$-tree such
that all points of the tree are vertices of continual valence, see comments in \cite{Ber}.
\hfill $\boxtimes$

\sm 

{\bf\punct $\R$-trees of ultrametric spaces.%
\label{ss:trees-of}}
Fix a countable or finite subset  $\Lambda$ in the ray $(0,\infty)$. 
  Consider a separable ultrametric space $X$, whose spectrum is contained 
  in $\Lambda$.
 We embed $X$ into a canonically defined $\R$-tree
$T(X)$ in the following way (see Berestovskii \cite{Ber0}, Huges \cite{Hug}).

We say that an {\it ideal closed ball} is a symbol $\beta(x,r)$, where $x\in X$, $r\ge 0$. We say
that two ideal balls $\beta(x_1,r_1)$, $\beta(x_2,r_2)$ coincide if their radii are equal,
$r_1=r_2=r$, and $d(x,y)\le r$. Otherwise they are different%
\footnote{It can happened that actual balls $\B^c(x,r_1)$, $\B^c(x,r_2)$ with $r_1\ne r_2$
 coincide as subsets in
$X$.}%
. If $r_1<r_2$, we write 
$\beta(x,r_1)\Subset \beta(x,r_2)$ and say '$\beta(x,r_2)$ {\it is over} $\beta(x,r_1)$' or 
'$\beta(x,r_1)$ {\it is under} $\beta(x,r_2)$'. More formally,
$\beta(x,r_1)\Subset \beta(y,r_2)$
if $r_1<r_2$ and $d(x,y)\le r_2$).

Points of the metric space $T(X)$
are ideal balls. The distance between $\beta(x,r_1)$ and $\beta(x,r_2)$ is $|r_1-r_2|/2$.
Let $\beta(x,r_1)$, $\beta(y,r_2)$ be {\it distinct}, i.e., $d(x, y)>\max (r_1,r_2)$.
Denote $r:=d(x,y)$, this is the radius of $\Subset$-minimal ideal ball $\Subset$-containing
$\beta(x,r_1)$, $\beta(y,r_2)$.
Then distance between $\beta(x,r_1)$, $\beta(y,r_2)$ is
\begin{equation}
\tfrac12 (r-r_1)+\tfrac12 (r-r_2),
\label{eq:rrr}
\end{equation}
see Fig.\ref{fig:2}.b.

Finally,  we get an $\R$-tree $T(X)$ containing $X$ (points of $X$ correspond to ideal balls
of zero radius) and the metric on $T(X)$ extends the metric on $X$.

\sm 

This $\R$-tree is equipped with an additional structure ({\it height}), namely the map
$\frh: T(X)\to [0,\infty)$ given by $\frh\bigl(\beta(x,r)\bigr)=r$. So $X$ consists of points of
height zero.

\sm 

{\sc Ascending paths.}
For each $x\in X$ we define an {\it ascending path} $L_x:r\mapsto \beta(x,r)$. Changing a parameter $r\to 2r$,
we get a geodesic path  $L_x(2r)$. By the definition,
$\frh(L_x(r))=r$.  For two ascending paths $L_x(r)$, $L_y(r)$ 
we have
$$
L_x(r)=L_y(r), \qquad \text{for $r\ge d(x,y)$.}
$$
More generally, for each $w=\beta(x,r)\in T(X)$ we have a unique  {\it ascending path} starting at
$w$,
\begin{equation}
t\mapsto L_{w}(t)=\beta(x, r+t),
\label{eq:t-mapsto}
\end{equation}
see Fig.\ref{fig:2}.c.
A point $v$ is over $w$ iff it $v$ is contained  in the  path $L_w$.
By the construction,
$$
T(X)=\cup_{x\in X} L_x.
$$
It is natural to add to $T(X)$ a point $\infty$ of infinite height, it corresponds to 'the ball of infinite
radius' containing the whole space $X$. All ascending  paths  lead to $\infty$.

The following important remark is obvious.

\begin{lemma}
\label{l:max}
For any $\xi$, $\eta\in T(X)$ the segment $[\xi,\eta]$ has a unique highest point, say $\zeta$.
Moreover, the segment $[\xi,\zeta]$ is a part of the path $L_\xi$ and the segment $[\eta,\zeta]$
is a part of the path $L_\eta$.
\end{lemma}

\sm 

{\sc Vertices, artificial vertices, and germs.}
For any vertex $v$ of $T(X)$ of valence $\ge 3$ we have $\frh(v)\in \Lambda$
and the corresponding  ball $\B^c(y, \frh(v))$ is perfect.
 This ball consists of all $x\in X$, for which $L_x$ passes through $v$ (or equivalently
of all $x$ lying under $v$).

 The set $\germ(v)$ contains one germ going up (the germ of the ascending path $L_v$),
other germs go down, they are enumerated by disjoint balls $\B^o(x_j,\frh(v))$
that are contained in $\B^c(x,\frh(v))$. We denote the set of germs going down
by
$$
\germ^\downarrow(v).
$$

Below {\it we assume that all points of any level $\lambda\in \Lambda$
are vertices of the $\R$-tree $T(X)$}. Of course, some vertices can be artificial.

\sm 

{\sc Remark.} Let $X$ be an ultrametric space, $\ov X$ be its completions. Balls
of nonzero radii in these spaces
are the same. So $T(X)\setminus X=T(\ov X)\setminus \ov X$.
\hfill $\boxtimes$

\sm

{\sc Naked segments.}
Let $\xi$, $\eta\in T(X)$, let $\xi$ be over $\eta$. 
We say that the {\it segment} $[\xi,\eta]$ is {\it naked}
if it does not contain actual vertices may be except $\xi$ and $\eta$.
Any naked segment is contained in a unique maximal naked segment.
A maximal naked segment corresponds to a family of ideal
balls $\beta(x, r)$, where $x$ is fixed and $r$ ranges in a segment
$[\rho_1,\rho_2]$ such that 

--- $\B^c(x,\rho_1)=\B^c(x,\rho_2)$;

--- $\rho_2=\infty$ or for
any $\delta>0$ we have
$
\B^c(x,\rho_2)\subsetneqq \B^c(x,\rho_2+\delta)$,

--- $\rho_1=0$ (in this case $x\in X$ is an isolated point)
 or  for
any $\delta>0$ we have
  $\B^c(x,\rho_1)\supsetneqq \B^c(x,\rho_1-\delta)$.
  
  \sm
  
 A {\it naked interval} is a naked segment without ends.
  
  Denote by $W$ the closure of the set of actual vertices
  Then the set $T_{>0}(X)\setminus W$ splits into a disjoint union 
  of naked intervals. 

\sm

{\bf \punct Embeddings of $\R$-trees to affine Hilbert spaces.%
\label{ss:tree-hilbert}}
We start with two obvious facts that are used in several places of the paper.

\sm

1) Let $H$, $H'$ be Hilbert spaces, let $\cX$ be a set. Consider maps
$\psi:\cX\to H$, $\psi:\cX\to H'$ such that for all $x$, $y\in \cX$ we have
$$\la \psi(x),\psi(y)\ra_H=\la \psi'(x),\psi'(y)\ra_{H'},$$
and in both cases linear spans of images of $\cX$ are dense in the Hilbert spaces.
Then there is a unique unitary operator $U:H\to H'$ such that
$U\psi(x)=\psi'(x)$ for all $x$.

\sm 

 2)  Let $H$, $H'$ be real Hilbert spaces, let $\cX$ be a set. Consider
$\phi:\cX\to H$, $\phi:\cX\to H'$ such that for all $x$, $y\in \cX$ we have
$$\|\phi(x)-\phi(y)\|_H=\| \phi'(x)-\phi'(y)\|_{H'},$$
and in both cases there is no proper closed affine subspace containing the image of $\cX$.
Then there is a unique (bijective) affine isometry  $R:H\to H'$ such that
$R\phi(x)=\phi'(x)$ for all $x$.

\begin{theorem}
\label{th:psi}
{\rm a)} Any  separable $\R$-tree $T$ admits an embedding $\psi:x\mapsto \psi_x$ 
to an affine real Hilbert space $\cH=\cH(T)$
such that:

\sm

\,\, $1^*$. For each pair of points $x$, $y\in T$ we have
\begin{equation}
\|\psi_x-\psi_y\|^2=d(x,y).
\label{eq:psi}
\end{equation}

\,\, $2^*$. The image of $\psi$ is not contained in a proper affine subspace.

\sm 

{\rm b)} An embedding $\psi$ is unique in the following sense. If $\psi'$ is another such embedding
to a Hilbert space $\cH'$, then there is an isometry $Q:\cH\to \cH'$ such that
$\psi'_x=Q(\psi_x)$. 

\sm

{\rm c)} If segments $[x_1,x_2]$, $[y_1,y_2]$ are disjoint or have a unique common point, when
lines 
$\psi_{x_1}\psi_{x_2}$, $\psi_{y_1}\psi_{y_2}$ are orthogonal%
\footnote{Two affine subspaces in a Hilbert space are {\it orthogonal} if the corresponding
linear subspaces are orthogonal.} in $\cH$. 

\sm

{\rm d)} For a subtree $J\subset T$ denote by $\cH(J)$ the minimal closed affine subspace
containing the $\psi$-image of $J$.
If subtrees $J_1$, $J_2\subset T$ are disjoint or have one common point, then 
affine subspaces $\cH(J_1)$, $\cH(J_2)$  are orthogonal.
\end{theorem}

For sets of vertices of simplicial trees the statement was obtained by Olshanski \cite{Olsh-new} and Watatani, \cite{Wat}, they considered  the sets of vertices of  simplicial trees.
The case of $\R$-trees  was mentioned in 
\cite{Ner-hie}. To be self-closed, we  present a geometric proof.

\sm

\begin{figure}
$$\epsfbox{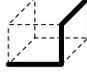}$$
\caption{To the proof of Theorem \ref{th:psi}.}
\label{fig:3}
\end{figure}

{\sc Proof.} First, consider a finite metric tree $J$ with $n$ vertices. Let us show that
there is an embedding $\psi^J$ of $\ver(J)$ to an $(n-1)$-dimensional Euclidean space,
satisfying our condition \eqref{eq:psi}. Let $p_1$, \dots, $p_{n-1}$ be edges of $J$, let $l_j$ be their lengths.
  We choose pairwise orthogonal
vectors $\vec f_j$ such that $\|f_j\|^2=l_j$. Fix $v\in\ver(J)$ and choose a point $\psi^J_v$
in an arbitrary way. For any $w\in\ver(J)\setminus v$ consider the shortest way
$p_{w,1}$, $p_{w,2}$, \dots from $v$ to $w$.  We set
$$
\psi^J_w:=\psi^J_v+\sum_{j} \vec f_{w,j}.  
$$
By the  Pythagoras 
theorem, we get the desired embedding, see Fig. \ref{fig:3}. 

Let $u\in \ver(J)$ have valence 1 or 2. Then the set $\ver(J)\setminus u$ remains to be a set 
of vertices of a metric tree, denote it by $J^{\backslash u}$.  Restricting the embedding $\psi^J$ to $J^{\backslash u}$
we get an embedding of $\ver(J^{\backslash u})$ to our Euclidean space satisfying the condition \eqref{eq:psi}
(but now the image is contained in an $(n-2)$-dimensional affine subspace).

\sm

Next, consider an  arbitrary separable $\R$-tree $T$ and choose a dense sequence  $s_1$, $s_3$, $s_5, \dots\in T$ containing all vertices.
For each $j$ consider the subtree $[[s_1,s_3,\dots,s_{2j-1}]]$ and denote by $s_{2j}$ the point of $[[s_1,s_3,\dots,s_{2j-1}]]$ nearest to
$s_{2j+1}$. For each $k$ we denote by $J_k$ the subtree $[[s_1,s_2,\dots,s_k]]\subset T$, by the construction
all vertices of this subtree of valence $\ne 2$ are contained in the set $\{s_1,\dots,s_k\}$. 
Adding artificial vertices of valence 2 we can assume that
$\ver(J_k):= \{s_1,\dots,s_k\}$. For each $k$ consider the embedding $\psi^{(k)}$ of $\ver(J_k)$ to the corresponding
Euclidean space $H(J_k)$. Then we have a chain of canonical embeddings
$$
\dots\xrightarrow{\iota_{k-1}}  H(J_k)\xrightarrow{\iota_k} H(J_{k+1}) \xrightarrow{\iota_{k+1}} \dots,
$$   
each embedding sends each  point $\psi^{(k)}_{s_\beta}\in H(J_k)$ with $\beta\le k$ to 
the point $\psi^{(k+1)}_{s_\beta}\in H(J_{k+1})$. 

Consider the union of spaces $ H(J_k)$ and denote its completion by $\cH(T)$, it is equipped with
the embedding of the set
 $\{s_j\}_{j=1}^\infty$. We extend this embedding to the whole $\R$-tree $T$ by continuity.
\hfill $\square$

\begin{corollary}
The group $\Isom(T)$ of isometries of an $\R$-tree $T$ acts on $\cH(T)$ by
affine isometric transformations.
\end{corollary}

\begin{corollary}
\label{cor:cEs}
 Let $T$ be an $\R$-tree. Let $0<s<1$. Then

\sm 

{\rm a)} There exists a Hilbert space $\cE_s(T)$ and a system of vectors
$\xi_x\in\cF(T)$, where $x$ ranges in $T$, such that

\sm

-- for each pair $x$, $y\in T$ we have
\begin{equation}
\la \xi_x,\xi_y\ra = s^{d(x,y)}.
\end{equation}

-- the linear span of vectors $\xi_x$ is dense in $\cE_s(T)$.

\sm

{\rm b)} A system $\xi_s$ is unique in the following sense. If $\xi'_x\in \cF_s'$
is another such system, then there is a unique unitary operator $U:\cF_s(T)\to \cF'_s(T)$
sending $\xi_x$ to $\xi'_x$.
\end{corollary}

 Indeed, according Schoenberg \cite{Sch}, 1937, for any $\lambda>0$
 there is an embedding  $\Phi_\lambda:h\mapsto \Phi_h$
of a separable Hilbert space $H$ to a Hilbert space $\cF(H)$  such that%
\footnote{In fact, this is a boson
Fock space, see, e.g., \cite{Ner-hie}, Subsect. 4.3}
$$
\la \Phi_h,\Phi_{h'}\ra_{\cF(H)}=\exp\bigl\{ -\lambda\|h-h'\|^2_H\bigr\}.
$$
We take the composition of embeddings $\xi:T\to \cH(T)$ and $\Phi_{-\ln s}:\cH(T)\to \cF(\cH(T))$.
\hfill $\square$

\begin{corollary}
\label{cor:olsh}
The group $\Isom(T)$ acts in $\cE_s(T)$ by unitary operators.
\end{corollary}

{\sc Remark.}
De facto this construction was used in old works by Ismagilov \cite{Ism67}, \cite{Ism73} (in the second paper he considered spherical functions on a certain $\R$-tree, however he did not observe a presence of trees), also construction was used explicitly by Olshanski \cite{Olsh-new} for non locally finite Bruhat-Tits 
trees. He also showed that such representations are infinite-dimensional limits
of representations of complementary series of groups of automorphisms of the usual Bruhat-Tits trees.
\hfill $\boxtimes$

 \sm
 
 {\bf \punct Embeddings of ultrametric spaces to affine Hilbert spaces.%
 \label{ss:ultra-hilbert}}

\begin{theorem}
\label{th:embedding}
Any separable ultrametric space $X$ admits an isometric embedding to
a Hilbert space. 
\end{theorem}

{\sc Proof.} We change  ultrametric metric on $X$ to $d'(x,y)=(d(x,y))^2$ and take
the corresponding $\R$-tree $T(X)$. Consider the embedding $\psi:T(X)\to \cH(T(X))$
described in Theorem \ref{th:psi}. We restrict it to $X\subset T(X)$
and get
$$
\qquad\qquad\qquad\qquad\qquad
\|\psi_x-\psi_y\|^2=d'(x,y)=(d(x,y))^2. 
\qquad\qquad\qquad\qquad
\square
$$

{\sc Remark.} This statement was formulated and proved in different ways by Vestfrid, Timan \cite{VT1}, \cite{VT2}
and Lemin \cite{Lem}. Earlier, in 1967, the statement was proved but not formulated by Ismagilov
\cite{Ism67}, Lemma 2%
\footnote{The conclusion of this lemma is precisely Menger's condition \cite{Men} of isometric embeddability of a metric space
to a Hilbert space.}. In 1997 Ismagilov \cite{Ism97} again proved (and formulated) statement. See, also,
Berestovskii \cite{Ber0} for additional details. \hfill $\boxtimes$

\sm

For an ultrametric space $X$
denote by $\cE_s(X)$ the subspace in $\cE_s(T(X))$ generated by vectors $\xi_x$, where $x$ ranges in $X$.
By construction, we have a unitary representation of the group $\Isom(X)$ in  $\cE_s(X)$.


\sm









\sm 

{\bf \punct Spaces of balls of a given radius.%
\label{ss:given}} Let $X$ be a separable ultrametric space with spectrum
$\Lambda$.
Denote by $(X)_h$ the set of all closed balls in $X$ of radius%
\footnote{Notice that for a fixed $h$ the sets of
geometrical balls $\B^c(x,h)$ and of ideal balls $\beta(x,h)$ coincide.}
 $h$.
There is a natural ultrametric on $(X)_h$, namely
$$
\delta\bigl(\B^c(x,h),\B^c(y,h)\bigr):=
\begin{cases}
0 \qquad &\text{if $\B^c(x,h)=\B^c(y,h)$};\\
d(x,y) &\qquad\text{otherwise}.
\end{cases}
$$
Clearly, we get a countable discrete ultrametric space with spectrum $\Lambda\cap(h,\infty)$.

\sm

Consider another ultrametric on $(X)_h$:
\begin{equation}
d_h\bigl(\B^c(x,h),\B^c(y,h)\bigr):=
\begin{cases}
0, \qquad &\text{if $\B^c(x,h)=\B^c(y,h)$};\\
d(x,y)-h &\qquad\text{otherwise}.
\end{cases}
\end{equation}
We get an ultrametric space with spectrum consisting
of (positive) numbers of the form $\lambda-h$, where $\lambda$ ranges in $\Lambda\cap (h,\infty)$.
We denote this ultrametric space by $[X]_h$.

By the construction, $[X]_h$ as a metric space coincides
with a level subset of $T(X)$, i.e., the set of all points of $T(X)$
of height $h$.

Denote by $T_{>h}(X)$ (resp. $T_{\ge h}(X)$)
the subset in $T(X)$ consisting of points $v$ such that $\frh(v)>h$ (resp. $\frh(v)\ge h$).
The
following statement is straightforward.

\begin{proposition}
Consider the space $[X]_h$ equipped with the metric $d_h$. Then we have an isometry
$T([X]_h)=T_{\ge h}(X)$.
\end{proposition}

Also, we have
$\bigl[[X]_{h_1}]_{h_2}=[X]_{h_1+h_2}$.

\sm

Generally, spaces $[X]_h$ are not complete. Let us describe their completions
$\ov{[X]}_h$.
We say that an {\it $h$-nested sequence of  balls} is a chain
$$
\B^c(a_1,r_1)\supset \B^c(a_2,r_2)\supset \B^c(a_3,r_3)\supset\dots,
$$
where $r_1>r_2>r_3>\dots$ and $\lim r_j=h$.
  For any $h$-nested sequence
 of balls we assign a   point of the completion, say $z$.
 Let $\B^c(a_j,r_j)$ and  $\B^c(a'_k,r'_l)$ be
 two  $h$-nested sequences.  There are two variants:
 
\sm  
 
 1) For each $j$, $k$   one of balls $\B^c(a_j,r_j)$,
 $\B^c(a'_k,r'_k)$ contains another.  Then we assume that these sequences are equivalent and the corresponding ideal points
 $z$, $z'$
 coincide. 

 \sm 
 
2) Some balls
 $\B^c(a_{j_0},r_{j_0})$,
 $\B^c(a'_{k_0},r'_{k_0})$ are  disjoint. Then balls
  $\B^c(a_{j},r_{j})$,
 $\B^c(a'_{k},r'_{k})$ are disjoint for all $j\ge j_0$, $k\ge k_0$
 and the corresponding distance 
 $\delta=\delta\bigl(\B^c(a_{j},r_{j}),\B^c(a'_{k},r'_{k})\bigr)$
 does not depend on $j\ge j_0$, $k\ge k_0$.
 We assume that $d_h(z,z')=\delta-h$.
 
 \sm
 
For a point $\B^c[a,h]\in [X]_h$
 we can assign a chain $\B^c[a,h+\tfrac 1j]$.
 
 \begin{proposition}
 The space of $h$-nested sequences of balls defined up to the equivalence
equipped with the metric $d_h$ 
  is the completion of
 $[X]_h$. For any point $z\in \ov{[X]}_h\setminus [X]_h$
 we can choose an $h$-nested chain consisting of perfect closed balls.
 \end{proposition}
 
This follows from 
Proposition \ref{pr:completion-tree} from the next subsection. 
 
 \sm

{\bf\punct Completions  of the trees $\boldsymbol{T(X)}$.%
\label{ss:completion}}  Denote by $\ov{T(X)}$ the completion of $T(X)$.
Denote by $\Omega_h(X)$ the set of all $w\in \ov{T(X)}$ such that $\frh(w)=h$.

\begin{proposition}
\label{pr:completion-tree}
{\rm a)} $[X]_h$ is dense in $\Omega_h(X)$ and so the level set $\Omega_h(X)$ is the completion of $[X]_h$
with respect to the metric $d_h$.

\sm

{\rm b)} Any sequence 
\begin{equation}
\beta(x_1,r_1)\Supset \beta(x_2,r_2)\Supset \dots,\qquad \text{where $r_j\to h+0$,} 
\label{eq:Supset}
\end{equation}
has a limit in $\Omega_h(X)$ and any element of $\Omega_h(X)$
is a limit of such chain.

\sm 

{\rm c)} For any point $\xi\in\Omega_h$ there is an ascending path  $L_\xi(t)$ such that 
$\frh(L_\xi(t))=h+t$.

\sm 

{\rm d)} For any point in $\Omega_h(X)\setminus [X]_h$ there is a chain \eqref{eq:Supset}
such that  all  balls $\B^c(x_j,r_j)$ are perfect,
$$
\B^c(x_1,r_1)\varsupsetneq \B^c(x_2,r_2)\varsupsetneq \dots
$$
\end{proposition}

{\sc Proof.} 
c) Let $\xi\in \Omega_h(X)$. Let $\eta\in T(X)$, $\frh(\eta)>h$. Since $\ov{T(X)}$ is an $\R$-tree,
there is a segment $[\xi,\eta]$. Chose an interior point $\zeta_1$ of  $[\xi,\eta]$, denote $\delta:=d(\xi,\zeta_1)$.
Since $\zeta_1$ is contained in $T(X)$, by Lemma \ref{l:max} there is a unique highest point on $[\eta,\zeta]$, say $\upsilon$,
and $\upsilon\in L_{\zeta_1}$. Next, take a point $\zeta_2\in [\xi,\zeta_1]$ such that $d(\xi,\zeta_2)=\delta/2$.
Then $[\upsilon,\zeta_2]\supset [\upsilon,\zeta_1]$ and by Lemma \ref{l:max}, $\upsilon$, $\zeta_1\in L_{\zeta_2}$, etc.
Further, we choose $\zeta_3\in [\xi,\zeta_2]$ such that $d(\xi,\zeta_3)=\delta/4$, etc. In this way, we get 
a sequence $L_{\zeta_k}$ of ascending paths, such that $L_{\zeta_{k+1}}$ includes $L_{\zeta_k}$. We take a union of these paths,
by construction, the sequence $\zeta_{k}$ converges to $\xi$, and we add $\xi$ to the union.

The statement b) follows from c). 

\sm

d), a) Let $\xi\in \Omega_h$. Assume that some segment $[\xi,\kappa]\subset L_\xi$ does not contain actual vertices. 
Let $\eta\in T(X)$ be an element of a small neighbourhood of $\xi$. Consider
the highest element $\upsilon$ of the segment $[\xi, \eta]$. We have
$d(\upsilon,\xi)\le d(\eta,\xi)$. Therefore $\upsilon$ can not be upper than
$\xi$. Therefore, $\upsilon=\xi$ and $\eta$ is under $\xi$. Hence, $\xi\in T(X)$, and so $\xi\in [X]_h$.
Also, $\xi$ is an isolated point of $\Omega$.

If there is a sequence of actual vertices   $v_j\in L_\xi$ convergent to $\xi$, then 
it corresponds to a sequence of perfect balls. This proves d).
This also proves a), since under each $v_j$ there is at least one point 
 $\zeta_j\in[X]_h$ and $d(\zeta_j,\xi)=2 d(v_j,\xi)$, see Fig.\ref{fig:2}.d. So the sequence
 $\zeta_j$ converges to $\xi$. 
 \hfill $\square$

\sm

\sm

{\bf \punct Dilative systems of balls and woolly subtrees.%
\label{ss:woolly}}
Let $X$ be an ultrametric space.
We say that a non-empty  family $\frS$ of perfect balls is {\it dilative}
if  $\B\in \frS$, $\B'\supset \B$ imply  $\B'\in\frS$. We consider
isolated points of $X$ as closed perfect balls and exclude other points of $X$ from considerations.

\sm

{\sc Remarks.}
a) To be accurate, we explain the meaning of the condition of
non-emptiness.
Recall that if the space $X$ is unbounded, then we consider 
the whole space as an open perfect ball $\B^o(a,\infty)$ of radius 
$\infty$. So non-emptiness for the case of an unbounded space means that a dilative family
contains this ball. If $X$ is bounded, then $X$ itself is a perfect
ball (closed or open), and this ball must be contained in any dilative family.

\sm

b) Our definition dilative system is motivated by Theorem \ref{th:extension} below,  the 
topic of our interest there are balls in $X$ and not points in $X$, so we exclude non-isolated points
of $X$ from considerations. We can vary a definition in the following way: we assume that points of $X$ are arbitrary prefect balls and
add the condition: for any nested sequence $\B^{\epsilon_j}(x_j,r_j)\in \frS$  such that $r_j\to 0$ 
the point of their intersection is contained in $\frS$. Then we will get the same set of dilative families.
\hfill $\boxtimes$ 

\sm

Let us reformulate this definition in terms of $\R$-trees $T(X)$, see Fig.\ref{fig:4}.a).
Denote by $\Isolated(X)$ the set of all isolated points of $X$.
Consider the set
$$
\wt T(X):=
\Bigl(T_{>0}(X)\bigcup \Isolated(X)\Bigr)\coprod
 \Bigl(\coprod\limits_{v\in \ver T(X)}\germ^\downarrow(v)\Bigr).
$$
We say that a non-empty subset $\bfS$ in $\wt T(X)$ is a {\it woolly subtree} in $T(X)$
if it satisfies the following conditions:

\sm

1) If $\xi\in T_{>0}(X)\cup \Isolated(X)$ is contained
in   $\bf S$, then all points of the ascending path $L_\xi$ are contained in
$\bfS$. Also for each vertex $v\in L_\xi$
the element of $\germ^\downarrow(v)$ directed to $\xi$ is contained
in $\bf S$.

\sm

2) If $\gamma\in \germ^\downarrow(v)$ is contained in $\bfS$, then $v\in \bfS$.

\sm

3) Let $[\xi,\eta]\subset T(X)$ be a maximal naked segment, let $\xi$ be over $\eta$.
 If an  interior point $\zeta$ of $[\xi,\eta]$,
is contained in $\bfS$, then $\eta\in \bfS$.
If $\frh(\xi)\in\Lambda$ and the unique element of $\germ^\downarrow(\xi)$
is contained in $\bfS$, then $\eta\in\bfS$. 

\sm

\sm

\begin{proposition}
There is a canonical one-to-one correspondence between 
dilative families of perfect balls in $X$ and woolly subtrees
in $T(X)$. 
\end{proposition}  

Let us describe this correspondence.
Let $\frS$ be a dilative family of balls.
We map $\frS\to \wt T(X)$ in the following way.

\sm

--- For each perfect ball $\B^c(x,r)\in \frS$
 we assign the point $\beta^c(x,r)$.
 
 \sm
 
 --- Let $\B^o(y,\rho)\in\frS$.
There are two variants.  First, the corresponding closed ball $\B^c(y,\rho)$
can be non-perfect (in this case $\B^c(y,\rho)=\B^o(y,\rho)$). Then we
draw the point $\beta(y,\rho)$ in the tree. Second,
let  $\B^c(y,\rho)$ be perfect. Then the point $v:=\beta(y,\rho)$
is an actual vertex of $T(X)$, and $\B^o(y,\rho)$
corresponds to an element of $\germ^\downarrow(v)$. We add this germ to the subtree. 

\sm

--- If a lower vertex of a maximal naked segment
$[\xi,\eta]$ is contained in the image of this map, we add the whole segment
to the subtree.

\sm 

 In this way we come to a woolly subtree.

\section{Amalgamate product\label{s:amalgams}}

\COUNTERS

We describe morphisms $X\to Y$ of the category $\cK_\Lambda$ in terms of partial
isomorphisms of trees $T(X)\to T(Y)$, products of morphisms correspond to products of
partial isomorphisms.

\sm

{\bf \punct Amalgams of ultrametric spaces.%
\label{ss:amalgamas2}}
Let $A$, $X$ be ultrametric spaces, let $f:A\to X$ be an isometric embedding. Then it induces
an isometric embedding of the corresponding $\R$-trees $\wt f: T(A)\to T(X)$. 

Let
$f:A\to X$, $A:g\to Y$ be isometric embeddings of ultrametric spaces. We extend them to isometric embeddings $ \wt f:T(A)\to T(X)$, $ \wt g:T(A)\to T(Y)$
and glue trees $T(X)$ and $T(Y)$ identifying points $ \wt f(u)\in T(X)$ 
and $\wt g(u)\in T(Y)$, where $u$ ranges in $T(A)$.
We get a new $\R$-tree, say $\cT$, equipped with a height function $\frh$. In particular,
we glue sets $X$, $Y$ identifying points $f(a)$ and $g(a)$, where $a$ ranges in $A$, and get 
a new set, say $Z$; now $A$ is identified with subset of $Z$. The tree $\cT$ determines a metric on $Z$. For two point $x\in X$ and $y\in Y$
we take geodesic path $[x,y]$ and the distance  is the length of this path
(or, equivalently the double height  of the upper point  of this path).

This construction defines the same operation on ultrametric spaces 
as was described in Subsect. \ref{ss:amalgama1}:

\begin{lemma} 
\label{l:minimax}
For $x\in X$, $y\in Y$ the distance  is
\begin{equation}
d(x,y)=\min_{a\in A} \max\bigl(d(x,a),d(a,y)\bigr).
\label{eq:min-max}
\end{equation}
\end{lemma}

\begin{figure}
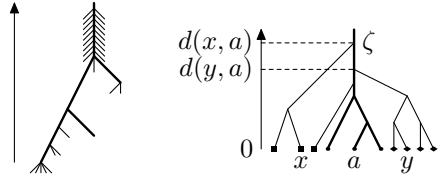

$$\epsfbox{trees2.1}\qquad\epsfbox{trees.6}$$
\caption{
a) To Subsect. \ref{ss:woolly}. A woolly subtree (a tree itself is not drawn).
\newline
b) To proof of Lemma \ref{l:minimax}. Fat points correspond to elements of $A$, squares to elements of $Y$, rhombi to elements of$X$.}
\label{fig:4}
\end{figure}

{\sc Proof.} The set $T(Z)\setminus T(A)$ splits into connected components, and each component
is contained either in $T(X)\setminus T(A)$ or in $T(Y)\setminus T(A)$. So a path $[x,y]$
contains a point  $v\in T(A)$.
Let $\zeta$ be the upper point of $[x,y]$. It  is upper than $v$ and therefore $\zeta\in T(A)$, i.e., $\zeta\in L_a$ for some $a\in A$. We have
 $d(x,y)=\frh(\zeta)$. 
By construction, $d(a,x)$, $d(a,y)\le \frh(\zeta)$.
 Consider the following 3 elements of $\germ^\downarrow(\zeta)$, the germ $x^\circ$ of the way $[\zeta,x]$,
the germ $y^\circ$ of the way $[\zeta,y]$ and a germ $a^\circ$ of the way $[\zeta,a]$. 
Since $\zeta\in [x,y]$ then   $x^\circ\nsim y^\circ$. Therefore at least one of germs $x^\circ$ or $y^\circ$
differs from $a^\circ$ and therefore at least one of distances $d(x,a)$, $d(y,a)$ is $\frh(\zeta)$, see Fig. \ref{fig:4}b.
 This also shows that $\min$ in \eqref{eq:min-max}
is reached.   
 \hfill $\square$

\sm 
 
{\bf \punct Partial bijections.%
\label{ss:partial}}
Let $L$, $M$ be sets. A {\it partial bijection} $\phi:L\to M$ is a bijection
of some subset $A\subset L$ to some subset $B\subset L$.
We say that $\dom \phi=A$ is the {\it domain of definition} of $\phi$, $\im \phi=B$
is the {\it image} of $\phi$. 

Let $\phi:L\to M$, $\psi:M\to N$
 be partial bijections. Their {\it product} $\psi \phi$ is a partial bijection
 $L\to N$ defined in the following way:
 
 \sm
 
 --- $\dom \psi \phi=\phi^{-1} (\im \phi\cap \dom \psi)$;
 
 \sm
 
 --- the image of $a\in \dom \psi \phi$ is  $\psi(\phi(a))$.
 
 \sm
 
So we have a {\it category $\PB$} whose objects are sets and morphisms are partial bijections.

This category is equipped with an {\it involution}. Let $\phi:L\to M$ be a partial bijection.
Then we have a partial bijection $\phi^*:M\to L$ with $\dom \phi^*=\im \phi$,
$\im\phi^*=\dom\phi$, and for $m\in \im\phi$ we assume $\phi^*(m)=\phi^{-1}(m)$.

\sm

A partial bijection $\phi:L\to L$ is {\it self-adjoint} ($\phi^*=\phi$) iff
$\dom \phi= \im\phi$ and the map $\phi$ is identical. In particular
any self-adjoint of ${\End\vphantom{q}}_{\PB}\,(L)$ is an {\it idempotent}, $\phi^2=\phi$ (and vice versa).
 
\sm 

{\bf\punct A variant of definition of the category $\boldsymbol{\PB}$.\label{ss:PB-variant}}
 Let us define another category $\wt\PB$, which in fact is equivalent to the category of partial bijections.
 {\it Objects} of  $\wt \PB$ are sets. 
A {\it morphism} $L\to M$ is a triple 
$\frp=(P,p_+,p_-)$,
where $P$ is a set, $p_+:L\to P$, $p_-:M\to P$ are embeddings such that
$p_+(L)\cup p_-(M)=P$. Two morphisms $\frp=(P,p_+,p_-)$,
 $\frp'=(P',p'_+,p'_-)\in \Mor_{\wt{\PB{\vphantom{\PB'}}}}(L,M)$ are equal if there is a bijection $\pi:P\to P'$ 
 such that $p'_\pm=\pi\circ p_\pm$.  
 
 We define the product  of $\frp=(P,p_+,p_-)\in \Mor_{\wt{\PB{\vphantom{\PB'}}}}(L,M)$, 
$\frq=(Q,q_+,q_-)\in\Mor_{\wt{\PB{\vphantom{\PB'}}}} (M, N)$. We consider the disjoint union $P\sqcup Q$,
identify pairs of points $p_-(m)$ with $q_+(m)$, where $m$ ranges in $M$,
and get new set $P\Join Q$. This set is equipped with embeddings 
$p_+:L\to P\Join Q$, $p_-:N\to P\Join Q$.
 We define $R\subset P\Join Q$ as $p_+(L)\cup q_-(N)$, and assume
 $r_+:=p_+$, $r_-:=p_-$. We get
  $\frr=(R,r_+,r_-)\in \Mor_{\wt{\PB{\vphantom{\PB'}}}}(L,M)$.
 
 \sm 
 
 Let us show that this multiplication  is reduced to the multiplication of partial bijections.
  For $\frp=(P,p_+,p_-)\in  \Mor_{\wt\PB}(L,M)$ we define a partial bijection
 $\phi:L\to M$ assuming
 $$
 \dom\phi= p_+^{-1}\bigl(p_+(L)\cap p_-(M)\bigr)
 \quad\text{and\quad $\phi(l)=(p_-)^{-1} (p_+(l))$ for $l\in \dom\phi$,}
 $$
 See Fig.\ref{fig:1}.a-b.
 Conversely, let $\phi:L\to M$ be a partial bijection. Then 
 we consider $L\sqcup M$, for each $l\in L$ we identify points 
$l\in L$ and $\phi(l)\in M$ and denote the result by $P$ or by  $L\sqcup_{\sim\phi} M$. We have identical
embeddings $p_+:L\to P$, $p_-:M\to P$.

 Let $\phi:L\to M$, $\psi:M\to N$ be partial bijections.
Then we get sets $P=L\sqcup_{\sim\phi} M$,
$Q=M\sqcup_{\sim \psi} N$. To evaluate the product $\frq\frp$,
we glue them by the set $M$ and get the set
$L\sqcup M\sqcup N$ with identified points $l\sim \phi(l)$,
$m\sim\psi(m)$ for $l\in \dom\phi$, $m\in\dom\psi$. 
Therefore for $l\in \dom \psi\phi$ we have $l\sim  \psi\phi(l)\in N$.
So the product $\frq\frp$ described above is 
the set 
$R=L\sqcup_{\sim\psi\phi} N$ equipped with canonical embeddings
$L\to R$, $M\to R$. So our product corresponds to a product
of partial bijections.

\sm    

{\bf \punct Partial isomorphisms of trees and amalgamation.%
\label{ss:partial-trees}}
Consider a finite ultrametric space $X$ whose spectrum is contained in $\Lambda$
and the corresponding $\R$-tree $T(X)$. In fact it is a usual tree,
points of $X$ are vertices of this tree of valence 1.
 Recall that we regard as vertices all points of $T(X)$ whose
 height is contained in $\Lambda$. We say that an {\it ascending subtree} is a convex
hull of $\infty$ and a finite non-empty collection of vertices  of $T(X)$ and
of  points of $X$.
Equivalently, we consider a finite collection of vertices $w_\alpha$ and a finite collection of points
$x_\beta\in X$ and take the union of ascending paths $\cup_\alpha L_{w_\alpha}\cup \cup_\beta L_{x_\beta}$.
 
  \sm

We define a category $\cK^*_\Lambda$. Its objects are
finite ultrametric spaces whose spectra are contained 
in $\Lambda$,
Let $X$, $Y\in\Ob(\cK^*_\Lambda)$ be such spaces. 
A morphism $X\to Y$ is a  height preserving isomorphism  of an ascending subtree 
$U\subset T(X)$ to an ascending subtree $V\subset  T(Y)$.
We multiply morphisms as partial bijections.

\begin{lemma}
The categories $\cK^*_\Lambda$ and $\cK_\Lambda$ are isomorphic.
\end{lemma}

{\sc Proof.} We define the third category $\wt \cK_\Lambda$. Its objects
are trees $T(X)$, where $X$ are finite ultrametric spaces. Morphisms
$T(X)\to T(Y)$ are triples $\wt\frp=(T(P),\wt p_+,\wt p_-)$, where $T(P)$
is a finite ultrametric space, and $\wt p_+:T(X)\to T(P)$, $\wt p_-:T(Y)\to T(P)$
are height preserving embeddings. We multiply morphisms as morphisms of $\wt\PB$.

Clearly, the category $\wt\cK_\Lambda$ is isomorphic to $\cK^*_\Lambda$.
On the other hand, consider  a morphism $\frp=(P,p_+,p_-):X\to Y$ of the category $\cK_\Lambda$.
Recall that $P$ is an ultrametric space equipped with isometric embeddings
$p_+:X\to P$, $p_-:Y\to P$  such that $P=p_+(X)\cup  p_-(Y)$. 
Then we have embeddings $\wt p_+:T(X)\to T(P)$, $\wt p_-:T(Y)\to T(P)$
and $\wt p_+(T(X))\cup\wt p_-(T(Y))=T(P)$. So we get a morphism, say $\wt\frp$,
of the category $\wt \cK_\Lambda$. By the construction of Subsect. \ref{ss:amalgamas2},
the product in $\cK_\Lambda$ corresponds to the product in $\wt \cK_\Lambda$.
\hfill $\square$

\section{Ultrametric Urysohn spaces\label{s:urysohn}}

\COUNTERS

\begin{figure}
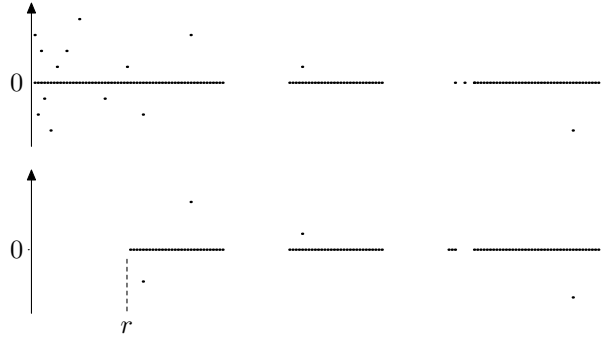

$$
\epsfbox{nguen.1}
$$
$$
\epsfbox{nguen.2}
$$
\caption{To Subsect. \ref{ss:nguen} and \ref{ss:glossary}. A graph of a function $\omegA(\lambda)$ and the function $\nU(\lambda)$ corresponding the ball $\B_\nU(r)$
containing
$\omegA$.}
\label{fig:5}
\end{figure}

{\bf \punct Ultrametric Urysohn spaces.%
\label{ss:nguen}} Fix a countable $\Lambda\subset (0,\infty)$. Following Nguyen \cite{Ngu} and Gao, Shao \cite{GS},
we define the  {\it ultrametric Urysohn} space $\U_\Lambda$ as the space
of functions $\omegA:\Lambda\to \Z$ whose support%
\footnote{The support of a function  is the set, where it is not zero.} is finite or is a sequence convergent to $0$, see Fig. \ref{fig:5}.
The distance is defined by
\begin{equation} 
 d(\omegA,\omegA')=\max\bigl\{ \text{$\lambda\in \Lambda$ such that 
 $\omegA(\lambda)\ne\omegA(\lambda')$}\bigr\}.
 \label{eq:d-omega}
 \end{equation}
 
 \begin{lemma}
We really get the space $\U_\Lambda$.
\end{lemma} 

{\sc Proof.}
The completeness and the separability of our space are obvious. Also obvious is the Ismagilov property,
see Subsect. \ref{ss:urysohn1}. Indeed, let $\omegA\in\U_\Lambda$, $\lambda^\circ \in\Lambda$. For each $n$ we define
a function $\omegA_n$ such that $\omegA_n(\lambda)=\omegA(\lambda)$ for $\lambda\ne \lambda^\circ$ and $\omegA_n(\lambda^\circ)=n$.
we get a family of functions including $\omegA$ such that $d(\omegA_n,\omegA_m)=\lambda^\circ$ for $n\ne m$. 
 
 Let us derive the Urysohn  property from the  Ismagilov property.
 Let the spectrum of a finite ultrametric space $Y$ with metric $\rho$
 be contained in $\Lambda$. Let 
$y\in Y$, $X=Y\setminus y$. Consider an isometric embedding $\psi$ of $X$ to our space.
Let $x_1$, \dots, $x_p$ be points of $X$ nearest to $y$, let $\rho(y,x_j)=s$.  Let $z_k\in Y$ be the remaining points. Choose in the collection $\{x_j\}$
a maximal subset $\{x_i^*\}$ such that distances between these points are $s$.
By Ismagilov property, there is a point $\theta \in \U_\Lambda$ such that $\theta\ne \psi(x_i^*)$ and $d(\theta,\psi(x_i^*))=s$.
The ultrametricity implies that all remaining distances $d(\theta,\psi(x_j))$ also are $s$.

Since $d(\psi(x_j),\psi(z_k))>d(\theta,\psi(x_j))$, we have 
$$d(\theta, \psi(z_k))=d(\psi(x_j),\psi(z_k))=d(x_j,z_k)=d(y,z_k).$$
We set $\psi(y):=\theta$. \hfill $\square$

\sm

{\sc Remark.} To be complete let us verify ultrahomogeneity. Let $X$, $Y$ be isomorphic finite
subspaces in $\U_\Lambda$, let $\psi:X\to Y$ be an isometry. Choose dense sequences $p_j$ and $q_j\in \U_\Lambda$. By Urysohn property we
can extend $\psi$ to an isometric embedding. 
$X\cup p_1\to \U_\Lambda$. Let keep the notation $\psi$ for this extended isometry. Next,
we add a point $q_1$ to $Y\cup \psi(p_1)$ and extend $\psi^{-1}$ to  $Y\cup \psi(p_1)\cup q_1$. Then we take $p_2$, etc. \hfill $\boxtimes$ 

\begin{lemma}
Elements of $\Mor (X,Y)$ are in one-to-one correspondence with
double cosets $\G[X]\backslash \G_\Lambda/ \G[Y]$.
\end{lemma}

{\sc Proof.} This correspondence to one side was described in Subsect. \ref{ss:amalgama1}. Let $\frp=(P,p_+,p_)$ be as above.
Then we have an isometric embedding $(p_-)^{-1}: p_-(Y)\to Y$. By the Urysohn property it can be extended to an isometric
embedding $r:P\to \U_\Lambda$. By the universality, there is $g\in \G_\Lambda$ such that
$g\bigr|_{X}=r^{-1}\circ p_+$. Then $\frp$ corresponds to the double coset $\G[X]\,g\, \G[Y]$.
\hfill $\square$

\sm 

{\bf\punct 
 A  glossary for Urysohn spaces and their  trees.%
 \label{ss:glossary}}

\sm
 
1) Closed balls (and ideal balls $\beta(\cdot)$)  of radius $r>0$ correspond to  functions
$\nU: (r,\infty)\cap \Lambda\to \Z$ with finite support, we denote such balls  by%
\footnote{$r$ is not necessary in notation since it is contained a definition of a function $\nU$.} $\B_\nU^c(r)$. Such a ball contains  points
$\omegA:\Lambda\to \Z$ satisfying the condition $\omegA\Bigr|_{(r,\infty)\cap \Lambda}=\nU$.

\sm

2) A closed ball
$\B_\nU^c(r)$ of radius $r$ contains a ball $\B_{\nU'}^c(r')$ of smaller radius $r'$ if 
 $\nU=\nU'\Bigr|_{(r,\infty)}$. 
 
 \sm
 
 3) The ascending path  $L_\omegA$ consists of all restrictions
of $\omegA$ to 'rays' $\Lambda\cap (r,\infty)$, where $r\in(0,\infty)$. 

\sm 


\sm

4) Open balls $\B^o$  of radius $r>0$ correspond to  functions
$\mU: [r,\infty)\cap \Lambda\to \Z$ with finite support, we denote such balls  by $\B^o_\mU(r)$. This  ball contains  points
$\omegA:\Lambda\to \Z$ such that $\omegA\Bigr|_{[r,\infty)\cap \Lambda}=\mU$. The corresponding closed 
ball corresponds to the function 
\begin{equation}
\nU:=\mU\Bigr|_{(r,\infty)},
\label{eq:munu}
\end{equation}
 i.e., we forget a value of $\mU$ at $r$.
Conversely, for $r\in\Lambda$ a closed ball $\B^c_\nU(r)$ is a disjoint union of open balls $\B^o_\mU(r)$
satisfying \eqref{eq:munu}. So elements of $\germ^\downarrow(\nU)$ are enumerated by values of $\mU$ at $r$

\sm

5) The completion $\ov{T(\U_\Lambda)}$ consists of pairs $(r,\omegA)$,
such that $r\ge 0$ and $\omegA:\Lambda\cap (r,\infty)\to \Z$
is a function, whose support is finite or is a decreasing sequence
convergent to $r$.

\sm 

\begin{observation} Any ball $\B^c(x,r)\subset \U_\Lambda$ is isometric to $\U_{\Lambda\cap (0,r]}$, 
any ball $\B^o(x,r)$ is isometric to $\U_{\Lambda\cap (0,r)}$.
\end{observation}

Recall that by $\Xi_{\epsilon,r}$ we  denote the set of all balls of the type $\B^\epsilon(\cdot,r)$,
see Subsect. \ref{ss:classification}.

\begin{observation}
\label{obs:canonical}
All balls $\B^\epsilon(a,r)$ of a given type $\epsilon$ and a given radius $r$ are canonically isometric in the following sense. 
For each pair $\B_1$, $\B_2\in \Xi_{\epsilon,r}$ we can define an isometry
{\rm(a 'shift')} $\gamma^{\B_1}_{\B_2}:\B_1\to \B_2$
such that for any $\B$, $\B_1$, $\B_2$, $\B_3\in \Xi_{\epsilon,r}$ we have
$$
\gamma^{\B_2}_{\B_3} \circ \gamma^{\B_1}_{\B_2}=\gamma^{\B_1}_{\B_3}, \qquad \gamma^{\B}_{\B}=1. 
$$
\end{observation}

To be definite, consider the case of closed balls. 
Let $\B_1$, $\B_2$ correspond to (finitely supported) functions $\nU_1$, $\nU_2:(r,\infty)\cap \Lambda\to\Z$.
Let $\omegA\in \B_1$, i.e., $\omegA$  is an extension of $\nU_1$ to the set $(0,r]\cap\Lambda$.
We define the function $\gamma^{\B_1}_{\B_2}\omegA:\Lambda\to \Z$ by
$$
\gamma^{\B_1}_{\B_2}\omegA(\lambda)
=\begin{cases}
\nU_2(\lambda), \qquad \text{if $\lambda>r$};
\\
\omegA(\lambda), \qquad \text{if $\lambda\le r$}.
\end{cases}
$$

{\sc Remark.} Of course, shifts  are not canonical operations from the point of view of intrinsic geometry
of the space $\U_\Lambda$. However, in the model with Puiseux-type series (see Subsect. \ref{ss:urysohn-construction})
$\gamma$ are really shifts of balls.
\hfill$\boxtimes$ 
 
\begin{observation}
For any $v\in\ver(T(\U_\Lambda))$ denote by $\G_\Lambda[v]$  the stabilizer of $v$ in $\G_\Lambda$.
Then the natural homomorphism $\G_\Lambda[v]$ to the group of permutations of $\germ^\downarrow(v)$
is surjective.
\end{observation}

{\sc Proof.} Let $\frh(v)=r$. Let $\sigma$ be a permutation of the set $\Z$.
For $\omegA\in \U_\Lambda$ we assume
$$
(\sigma \omegA)(\lambda)=\begin{cases}
\omegA(\lambda), \qquad \text{if $\lambda\ne r$;}
\\
\sigma(\omegA(r)), \qquad \text{if $\lambda\ne r$,}
\end{cases}
$$ 
and we realized arbitrary permutation of  $\germ^\downarrow(v)$.
\hfill $\square$

\sm

{\bf \punct Countable Urysohn spaces.} A {\it countable Urysohn space} $\U^{\coun}_\Lambda$
 is a countable space, for which the Urysohn  property holds.
 Such space exists and is unique, see Nguyen \cite{Ngu}, Gao, Shao \cite{GS}.
 It is realized as the space of all finitely supported functions $\omegA:\Lambda\to\Z$
 with the same distance \eqref{eq:d-omega}. Clearly, $\U^{\coun}_\Lambda$ is a dense subset
 in $\U_\Lambda$ and 
  $$T_{>0}\bigl(\U^{\coun}_\Lambda\bigr)=T_{>0}\bigl(\U_\Lambda\bigr).$$
  
Denote by $\G^\coun_\Lambda$ the group of isometries of this space. The ultrahomogeneity implies
that $\G^\coun_\Lambda$   is dense
in $\G_\Lambda$.

Next, take $h>0$ and consider the shit $h+\Lambda$ of the spectrum $\Lambda$
and the corresponding Urysohn space $\U_{h+\Lambda}$. This space is discrete,
however, we have a one-one-to one correspondence between $\U^{\coun}_\Lambda$ and  $\U_{h+\Lambda}$,
for each function $\omega(\lambda)\in \U^{\coun}_\Lambda$ we assign the function
$\mu(\lambda)=\omega(\lambda+h)$, which is contained in $\U_{h+\Lambda}$. Since functions
$\mu(\lambda)\in \U_{h+\Lambda}$ have finite support, we get a bijection.
This bijection also identifies the corresponding groups,
$$
\G^{\coun}_\Lambda\simeq \G_{h+\Lambda}.
$$

  
  \sm 
  
  {\bf \punct Level sets in $T\boldsymbol{\bigl(\U_\Lambda\bigr)}$.}
  We use notation of Subsect. \ref{ss:given}-\ref{ss:completion}. The following statements
  is obvious.
  
  \begin{proposition}
   For each $h>0$ we have
\begin{align*}   
   &\bigl[\U_\Lambda\bigr]_h=\U^\coun_{[\Lambda]_h};
   \\
   &\Omega_h=\ov{\bigl[\U_\Lambda\bigr]}_h=\U_{[\Lambda]_h}
   \\
    &T_{\ge h}(\U_\Lambda)=T\bigl(\U^\coun_{[\Lambda]_h})\bigr).
   \end{align*}
   \end{proposition}
   
   \begin{proposition}
{\rm a)} The  natural homomorphisms 
$$\G_\Lambda\to \G^\coun_{[\Lambda]_h}, \qquad
\G^\coun_\Lambda\to \G^\coun_{[\Lambda]_h}$$
are surjective.

\sm

{\rm b)} For each $h>0$ the group $\G_\Lambda$ is a wreath product
$$
\G_\Lambda=\G^\coun_{[\Lambda]_h}\ltimes \prod_{\B\in \Xi_{c,h}}
\Isom(\B).
$$
  \end{proposition}

 {\sc Proof.} It is sufficient to verify the surjectivity. An element of $\G^\coun_{[\Lambda]_h}$ determines an isometry
 of the set $[X]_h$ whose elements enumerate balls. For each such a ball $\B$ we take the shift $\gamma^{\B}_{g(\B)}$
and get an element of $\G_\Lambda$.
\hfill $\boxtimes$

\sm

 {\bf \punct Remarks. Homeomorphisms of the spaces $\boldsymbol{\U_\Lambda}$ and isomorphisms of the groups $\boldsymbol{\G_\Lambda}$.}

\sm

a) Let $X$ be an ultrametric space. Let $u:(0,\infty)\to (0,\infty)$ be a
 monotone map.  
Then the function 
$\wt d(x,y)=u\bigl(d(x,y)\bigr)$ is an ultrametric on $X$.

\sm

b) Let $\Lambda$, $\Lambda'\subset (0,\infty)$ be countable sets. Let $0$ be or not  be
an isolated point of $\Lambda$ and $\Lambda'$ simultaneously. 
 Let $u:\Lambda\to \Lambda'$ be
an order preserving bijection. According our construction it identifies the spaces $\U_\Lambda$ and 
$\U_{\Lambda'}$, this identification is a homeomorphism. Also
we get an isomorphism of the isometry groups, $\G_\Lambda\simeq \G_{\Lambda'}$. In particular,
if both $\Lambda$, $\Lambda'$ are dense in $(0,\infty)$, then the corresponding isometry groups are isomorphic.

\sm 

c) For a countable  $\Lambda\subset (0,\infty)$ consider the group $\Aut(\Lambda)$
 of order preserving bijections $\Lambda\to\Lambda$. Then $u\in \Aut(\Lambda)$ determines an automorphism
 of the group $\G_\Lambda$ and we get an action of $\Aut(\Lambda)$ by automorphisms of $\G_\Lambda$. Also we get a semi-direct product 
 $\Aut(\Lambda)\ltimes \G_\Lambda$.
 \hfill $\boxtimes$

 \section{Multiplicativity\label{s:multiplicativity}}
 
 \COUNTERS
 
 The purpose of this section is the proof of the multiplicativity
 theorem \ref{th:multiplicativity}, it is occupies Subsect. \ref{ss:countable-product}-\ref{ss:density}. 
  The last Subsect. \ref{ss:correspondence} contains some immediate corollaries of multiplicativity.
  In Subsection \ref{ss:induced} we recall a definition 
 of a representation of a topological group induced from an open subgroup.
  Subsect. \ref{ss:symmetric} contains preliminaries on unitary representations of
infinite symmetric group (the Lieberman theorem and Olshanski's comments on it). In fact, for further proofs we need only Lemma 
\ref{l:olsh}.

\sm

{\bf \punct Induced representations.%
\label{ss:induced}} Recall a definition of induced representations in a degree of generality, which
is sufficient for our purposes.
 See, e.g., Mackey \cite{Mack0}, \cite{Mack}, Chapter 3, Kirillov, \cite{Kir}, Subsect 13.1-13.2.
  Let $G$ be a separable topological group, $H$ its open subgroup. Then the homogeneous space
$X=H\backslash G$ is countable, let $x_0$ be a point whose stabilizer is $H$.
 Let $V$ be a Hilbert space, let $\UU(V)$
be the group of unitary operators in $V$. We consider the space
$\ell^2(X,V)$ consisting of $V$-valued functions $f$ on $X$ satisfying
the condition
$$\|f\|^2:=\sum_{x\in X} \|f(x)\|^2_V<\infty.$$
It is equipped  with the inner product
$$
\la f_1,f_2\ra:=\sum_{x\in X} f_1(x) \ov{f_2(x)}.
$$
Consider unitary representations $T(g)$ of the group $G$ in $\ell^2(X,V)$ given by
formulas of the type
$$
T(g)f(x)=A(x,g) f(x,g),
$$
where $A:G\times X\to \UU(V)$ is a continuous function (a '{\it cocycle}'). The equality
$$T(g_1)T(g_2)=T(g_1g_2)$$
 is equivalent to the condition
\begin{equation}
A(x, g_1g_2)=A(x,g_1)A(xg_1, g_2).
\label{eq:AAA}
\end{equation}
Fix a function $\phi:X\to \UU(V)$.
If we perform a gauge  transformation $f(x)\to \phi(x)f(x)$, then $A(x,g)$
transforms as 
\begin{equation}
A(x,g)\mapsto A'(x,g)= \phi(x)^{-1} A(x,g) \phi(xg).
\end{equation}
We say that cocycles $A(g,x)$ and $A'(g,x)$ are {\it equivalent}.
In this case the corresponding representations are equivalent and the gauge transformation
provides this equivalence.

For $h_1$, $h_2\in H$ we have
$$
A(x_0, h_1h_2)=A(x_0, h_1)A(x_0,h_2),
$$
i.e., $\tau(h):=A(x_0,h)$ is a unitary representation of the group $H$ in $V$. 

\begin{proposition}
\label{pr:induced}
For each unitary representation $\tau$ there exists and is  unique   up to the equivalence
a cocycle $A(x,g)$
such that $A(x_0,h)=\tau(h)$ for $h\in H$.
\end{proposition}

Then the representation $T$ of $G$ is called a representation  {\it induced} from a representation $\tau$ of the subgroup $H$,
 notation is  
$$
T=\Ind_H^G(\tau).
$$ 




To be explicit, we  write a (noncanonical) formula for a  cocycle $A(x,g)$ corresponding $\tau$.
Choose a section $s:X\to G$ such that $x_0 s(x)=x$, $s(x_0)=1$. We can
choose unitary operators $A(x_0, s(x))$ in an arbitrary way, to be definite we
set $A(x_0, s(x))=1$. Substituting to \eqref{eq:AAA}:
 $x\to x_0$, $g_1\to s(x)$, $g_2\to s(x)^{-1}$, we get  $A(x,s(x)^{-1})=1$. For each $x\in X$ we decompose  $g\in G$ as
$$g=s(x)^{-1} h s(xg),$$
 where $h=h(x,g)=s(x)g s(xg)^{-1}\in H$. Then
\begin{multline*}
A(x,g)=A\bigl(x,s(x)^{-1} h s(xg)\bigr)=\\= A(x, s(x)^{-1}) A(x s(x)^{-1},h) A\bigl(x s(x)^{-1}h, s(xg)\bigr)
=1\cdot A(x_0,h)\cdot A(x_0,s(xg))=
\\=
 1\cdot \tau(h)\cdot 1= \tau\bigl(s(x)g s(xg)^{-1} \bigr).
\end{multline*}

\sm

{\sc Remarks.}  The representation of $G$ induced from the trivial one-dimensional representation
of $H$ is the quasi-regular representation of $G$ in $\ell^2(H\backslash G)$.
\hfill $\boxtimes$

\sm

{\sc Remark.}
There are two 'invariant' constructions of induced representations, in terms of $V$-valued functions
on $G$, see, e.g., \cite{Mack0}, \cite{Mack}, Sect. 3.2, and 
 in terms of skew products, see, e.g.,  \cite{Kir}, Subsect. 13.4.
On the other hand, see transparent explanations for finite groups (which remain to be valid in our case)  in Serre \cite{Ser}, Sect. 7. 
\hfill $\boxtimes$

\sm

{\bf \punct Induction in stages.} Let $G$ be a separable topological group, 
$N\supset H$ be open subgroups in $G$, let $\tau$ be a unitary representation
of $H$.
Then we have the following induction in stages (see, e.g., \cite{Mack}, Sect. 3.2, \cite{Kir}, 13.1) 
\begin{equation}
\Ind_N^G\bigl(\Ind_H^N(\tau)\bigr)\simeq\Ind_H^G(\tau).
\label{eq:stages}
\end{equation}

\begin{lemma}
\label{l:stages}
Let $G\supset N\supset H$ be as above, let $H$ be a normal subgroup in $N$ of finite index.
Then the quasiregular representation of $G$ in $\ell^2(H\backslash G)$ is equivalent 
to
$$
\bigoplus_{\nu\in \wh{N/H }} \dim (\nu)\, \Ind_N^G(\nu),
$$
where $\nu$ ranges in the set $\wh{N/H}$ of irreducible representations  of $N/H$, we understand such $\nu$ as a representations of $N$ trivial
on $H$.
\end{lemma}
 
{\sc Proof.} The statement is corollary of the following remarks.

\sm

1) A representation of a finite group $K$ induced from the trivial representation is the regular representation 
of $K$. Each irreducible representation $\nu$ of $K$ enters to the regular representation $\dim \nu$ times.

\sm

2) The representation  of $N$ induced from the trivial representation of the normal subgroup
$H$ is the regular representation of $N/H$ considered as a representation of $N$. 

\sm

It remains to apply 
the induction in stages.
\hfill $\square$

\sm 

{\bf \punct The infinite symmetric group.%
\label{ss:symmetric}} Denote by $S_\infty$ the group
of all permutations of the set $\N$. We also use notation $ S_\infty (\Omega)$
for the group of all permutations of a countable set $\Omega$. This group
is Polish, stabilizers of finite subsets are open subgroups and such subgroups form a
fundamental system of neighbourhoods of unit. 

Denote by $\cB_\infty$ the semigroup
of infinite $(0,1)$-matrices, i.e., matrices, which have $\le 1$ unit in each column 
and $\le 1$ unit in each row, all other matrix elements are 0. Equivalently, $\cB_\infty$ is the semigroup
of all partial bijections of $\N$. So $\cB_\infty\supset  S_\infty$.
We equip $\cB_\infty$ with the topology of element-wise convergence, 
$\sigma^{(k)}\to \sigma$ if we have convergences (in fact a stabilization) of all
matrix elements $\sigma^{(k)}_{ij}\to \sigma_{ij}$. This convergence is equivalent to the weak operator convergence
of the corresponding operators
in the Hilbert space $\ell^2(\N)$.  The semigroup $\cB_\infty$ is  compact and the product is
separately continuous, the group $ S_\infty$ is dense in $\cB_\infty$.
See, e.g. \cite{Ner-book}, Sect. 8.1.

The classification of all unitary representation of $ S_\infty$ was obtained by Lieberman \cite{Lie}
(see also an exposition given by Olshanski \cite{Olsh-kiado} and by the author \cite{Ner-book}, Sect. 8.1-8.2).
Denote by $S_n\times  S_{\infty-n}$ the subgroup in $S_\infty$ consisting of transformations sending the set $\{1,2,\dots,n\}\subset\N$
to itself.

\begin{theorem}
\label{th:lieberman}
{\rm a)} Any irreducible unitary representation of $ S_\infty$ is induced from
an irreducible representation of a subgroup $S_n\times  S_{\infty-n}$ that is trivial
on $ S_{\infty-n}$.

\sm 

{\rm b)} The group $ S_\infty$ has type%
\footnote{Recall that a {\it von Neumann algebra} is an algebra $\frA$ of operators in a Hilbert space $H$,
which is closed with respect to the weak operator topology. It has {\it type $I$} if there is a decomposition
of $H$ into a direct integral such that $\frA$ consists of all bounded operators preserving this decomposition.
A topological group $G$ has {\it type $I$} if for any its unitary representation $\rho$
the von Neumann algebra generated by all operators $\rho(g)$ has type $I$. See books by Mackey \cite{Mack}, 1, 3.5, Appendix 3.A,
Dixmier \cite{Dix}, 4, 5, 8, 9,   Addendum A, Kirillov \cite{Kir}, 4.5, 8.4, 11.2. In particular, they explain a fundamental
importance of this notion.  
 }
 $I$ and any its unitary representation is a direct sum of
irreducible representations.
\end{theorem}

Olshanski \cite{Olsh-kiado} proved the following statement, which in fact is equivalent to the Lieberman
theorem.

\begin{theorem}
Any unitary representation of $ S_\infty$ admits a weakly continuous extension to the semigroup
$\cB_\infty$.
\end{theorem}

We need a simple addition to this statement: the matrix $0\in\cB_\infty$ acts
in any unitary representations $\rho$ of $ S_\infty$ as the orthogonal projector $\Pi$ to the subspace of $ S_\infty$-fixed vectors,
see \cite{Olsh-kiado}, \cite{Ner-book}, Corollary 8.1.5. 

In particular, this implies the following statement:

\begin{lemma}
\label{l:olsh}
Consider the group $S_\infty(\Z)$ and its elements
\begin{equation}
\theta_j:n\mapsto n+j.
\label{eq:theta-j}
\end{equation}
Then $\rho(\theta_j)$ weakly converges to $\Pi$ as $j$ tends to $\infty$.
\end{lemma}

\sm

{\bf\punct The countable product of infinite symmetric groups.
\label{ss:countable-product}}

\begin{lemma}
\label{l:Theta}
Consider a direct product 
$$S_\infty(\Z)^\infty:=S_\infty(\Z)\times S_\infty(\Z)\times \dots$$
 of a countable number of copies $S^{(k)}_\infty(\Z)$  of the
group $S_\infty(\Z)$. In each copy we take a sequence $\theta_j^{(k)}$ as in
\eqref{eq:theta-j}. Consider the sequence
\begin{equation}
\Theta_j:=\bigl(\theta_j^{(1)}, \theta_j^{(2)}, \theta_j^{(3)},\dots\bigr).
\label{eq:Theta_j}
\end{equation}
Then for any unitary representation $\rho$ of the group $S_\infty(\Z)^\infty$
the sequence $\rho(\Theta_j)$ converges to  the orthogonal projector
$\Pi$ to the subspace of $S_\infty(\Z)^\infty$-fixed vectors.
\end{lemma}

The statement is a corollary of the following proposition

\begin{proposition}
\label{pr:tensor}
Any unitary representation $\rho$ of the group $(S_\infty)^\infty$ is a direct sum of irreducible representations.
Each irreducible unitary representation $\tau$ of this group is a tensor product,
$$
\tau(g_1,g_2,g_3,\dots)=\mu(g_1)\otimes \mu(g_2)\otimes \mu(g_3)\otimes\dots
$$
of irreducible representations of factors $S_\infty^{(k)}$. Moreover, all but a finite
number of factors
 are trivial one-dimensional representations.
\end{proposition}

The group $( S_\infty)^\infty$ is an inverse limit of oligomorphic groups $( S_\infty)^n$,
and our statement is a special case of Tsankov \cite{Tsa} (Theorem 1.3, Theorem 4.2) classification theorem.

On the other hand,  Proposition \ref{pr:Ism} formulated below  reduces the statement to Theorem \ref{th:lieberman}.b. 

\sm

{\bf \punct Stabilizers of finite subsets.%
\label{ss:stabilizers}} Consider a finite subset $A$ in the Urysohn space
$\U_\Lambda$.  Denote by 
$$\G[A]\subset \G_\Lambda$$
 the point-wise stabilizer of $A$ in
$\G_\Lambda$. Consider the subtree $T(A)\subset T(\U_\Lambda)$. For each $x\in\U_\Lambda\setminus A$
we consider the ascending path  $L_x$. Let $v(x)$ be the lowest point of 
$L_x\cap T(A)$. Conversely, consider a vertex%
\footnote{Recall that for $T(\U_\Lambda)$ all points of any height $\lambda\in \Lambda$ are vertices of infinite valence.
} $v$ of $T(A)$ and consider the set $D_v$ consisting
of all $x\in \U_\Lambda\setminus A$, for which $v=v(x)$. The  complement to $A$ splits into a disjoint
union 
$$
\U_\Lambda\setminus A= \coprod_{v\in \ver(T(A))} D_v.
$$
Any element  $g\in \G[A]$ preserves  this splitting.

Let $x\in D_v$, $y\in D_w$ and $v\ne w$. Then a value $d(x,y)$ depends only on $v$ and $w$, it equals to the maximal height of points of
the segment $[v,w]$. This implies that
$$
\G[A]=\prod_{v\in \ver(T(A))} \Isom(D_v).
$$ 

Fix a vertex $v\in T(A)$. It determines a closed ball, say $\B^c(a,\frh(v))$, where $a\in A$.
This ball is a union of a countable family of open balls $\B^o(\cdot,\frh(v)$
enumerated by the set $\germ^\downarrow(v)$. Finite number
of such balls contains points of  $A$. Denote the remaining balls by $\B^o_v\bigl(p_j,\frh(v)\bigr)$,
where $j$ ranges in $\Z$, so 
$D_v=\cup_{j\in \Z} \B^o_v\bigl(p_j,\frh(v)\bigr)$. For us it will be convenient to distinguish copies of $\Z$
corresponding to different $v$, so we will write $\Z^{(v)}$. In this notation,
$$
D_v=\coprod_{j\in \Z^{(v)}} \B^o_v\bigl(p_j,\frh(v)\bigr).
$$

{\sc Remark.}
Notice that $D_v$ as an ultrametric space is non-canonically isometric to
$\B^c(a,\frh(v))\simeq \U_{\Lambda\cap (0,\frh(v)]}$.
\hfill $\boxtimes$

\sm


According  Observation \ref{obs:canonical}, balls  $\B^o_v\bigl(p_j,\frh(v)\bigr)$ are canonically isometric.
So, $\Isom(D_v)$ is a wreath product,
$$
S_\infty(\Z^{(v)})\ltimes \prod_{j\in \Z} \Isom(\B^o_v\bigl(p_j,\frh(v)\bigr), \qquad \Isom(\B^o_v\bigl(p_j,\frh(v)\bigr)\simeq\G_{\Lambda\cap(0,\frh(v))}.
$$
The group $S_\infty(\Z^{(v)})$ acts by  permutations of balls $\B^o_v\bigl(p_j,\frh(v)\bigr)$, each factor $\Isom(\B^o_v\bigl(p_j,\frh(v)\bigr)$  
acts by isometries on the corresponding $\B^o_v\bigl(p_j,\frh(v)\bigr)$ and trivially outside it.

So we have an epimorphism
\begin{equation}
\kappa_A:
\G[A]\to \prod_{v\in \ver T(A)} S_\infty(\Z^{(v)}).
\label{eq:ker-kappa}
\end{equation}
Its kernel 
$$K[A]:=\ker\kappa_A
$$ 
consists of elements of $\G[A]$, which are identical on all sets $\germ^\downarrow(v)$
for all vertices of $T(A)$. 
On the other hand, we have an embedding 
$$
\prod_{v\in \ver T(A)} S_\infty(\Z^{(v)})\to \G[A].
$$

{\bf\punct Reformulation of the definition of the product in the category 
$\boldsymbol{\cK_\Lambda}$.}

\begin{proposition}
\label{pr:action-balls}
Consider a  double coset $\G[Y]\,g\,\G[X]$, let $\pi:T(X)\to T(Y)$ be the corresponding
partial isomorphism of trees. Let $\xi=\beta(x,r)\in T(X)$, 
$\eta=\beta(y,r)\in T(Y)$. Then $\pi$ sends $\xi\mapsto\eta$
if and only if $g(\B^c(x,r))=\B^c(y,r)$.
\end{proposition}

{\sc Proof.} We consider the tree $T(g(X)\cup Y)$. A point $\beta(z,r)$ of this tree
is contained in $T(g(X))\cap T(Y)$ if  the ball $\B^c(z,r)$ contains a  point
of $g(X)$ and a point of $T(Y)$.
\hfill $\square$




\sm

Fix a finite subset $Y\subset \U_\Lambda$.
Consider its pointwise stabilizer $\G[Y]$
 of $Y$, such stabilizers were described  in the previous subsection. Let $v$ ranges in  $\ver(T(Y))$.
In each $S_\infty(\Z^{(v)})$, we consider the sequence $\theta_j^v:n\mapsto n+j$.
For each $j$ consider an element 
$$
\Theta_j^Y:=\{\theta^v_j\}_{v\in \ver(T(Y))} \in \prod_{v\in \ver T(Y)} S_\infty(\Z^{(v)})\subset  \G[Y]. 
$$

\begin{proposition}
\label{pr:Theta-move}
Let $g_1$, $g_2\in \G_\Lambda$. Then for any finite $X$, $Y$, $Z\subset \U_\Lambda$
the sequence
$$
\G[X]\cdot g_1 \Theta_j^Y g_2 \cdot \G[Z]\in \G[X]\setminus \G_\Lambda/\G[Z]
$$
is eventually constant and its limit coincides with the product of double cosets 
$$
\bigl(\G[X]\, g_1\,\G[Y]\bigr)\diamond
 \bigl(\G[Y] \, g_2\, \G[Z]\bigr)
$$
defined in Subsect. {\rm\ref{ss:double-cosets}}.
\end{proposition}

{\sc Proof.}
Without loss of generality, we can assume $g_1=g_2=1$. Otherwise, we pass 
to sets
 $X':=Xg_1$, $Z':=Zg_2^{-1}$.
So it is sufficient to examine a sequence 
$$
\G[X]\cdot  \Theta_j^Y  \cdot \G[Z].
$$
 Denote by $\pi_j$ the corresponding partial bijections
$T(X)\to T(Z)$. To examine them, we use  Proposition \ref{pr:action-balls}.
 Denote by $\pi^\circ$ the partial isomorphism $T(X)\to T(Z)$
whose domain is $T(X)\cap T(Y)\cap T(Z)$ and which 
is identical on this domain. We must prove that for sufficiently
large $j$ we have $\pi_j=\pi^\circ$.

\begin{figure}
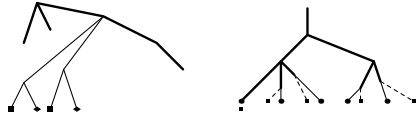

$$
\epsfbox{trees.7}\qquad \epsfbox{trees.8}
$$
\caption{To the proof of Proposition \ref{pr:Theta-move}. We mark $T(X)\cap T(Y)\cap T(Z)$ by thick segments,
points of $X$ by squares, points of $Z$ by rhombi. We show the tree $T(X\cup Y\cup Z)$ and $T(X\Theta_j^Y\cup Y\cup Z)$.
}
\label{fig:6}
\end{figure}

Let $L$, $M$ be sets, let $\phi$, $\psi:L\to M$ be partial bijections.
We say that $\psi$ is an {\it expansion} of $\phi$ if $\dom \psi\supset \dom\phi$
and $\phi(l)=m$ implies $\psi(l)=m$. 

Since $\Theta_j^Y$ is identical on $Y$, partial bijections
 $\pi_j$ are expansions of $\pi^\circ$. Let
  $\xi=\beta(x,r)\in T(X)$, $\eta=\beta(z,r)\in T(Z)$
   be not in $T(X)\cap T(Y)\cap T(Z)$. Since $\Theta_j\in \G[Y]$, a coincidence
    $\Theta_j^Y(\B^c(x,r))=\B^c(z,r)$
   is possible only if $\B^c(x,r)$, $\B^c(y,r)$ are contained 
   in one set $D_v$ defined in the previous subsection. But
   the sets $D_v\cap X$, $D_v\cap Z$ are finite, and for sufficiently
   large $j$ we have $\theta_j^{(v)}(\B^c(x,r))\ne \B^c(z,r)$ for any 
   $x\in D_v\cap X$, $z\in D_v\cap Z$.
   \hfill $\square$



\sm

{\bf \punct  Approximation of  projectors $\boldsymbol{\Pi[A]}$.}

\begin{lemma}
\label{l:weakly}
For any unitary representation $\rho$ of the group $\G_\Lambda$ in a Hilbert space $H$ 
the sequence $\rho(\Theta_j^A)$ weakly converges to the operator of orthogonal projection $\Pi=\Pi[A]$
to the subspace of $\G[A]$-fixed vectors.
\end{lemma}

{\sc Proof.}
We use the description of $\G[A]$ given in Subsect. \ref{ss:stabilizers}.
 Denote by $V$ the subspace of $\prod_{v\in \ver T(A)} S_\infty(\Z^{(v)})$-fixed vectors.
 By Lemma \ref{l:Theta}, the sequence  $\rho(\Theta_j^A)$ weakly converges to the operator $\Pi$ of orthogonal projection 
 to $V$. It is sufficient to prove the following lemma:
 
 \begin{lemma}
 \label{l:fixed-fixed}
 Any $\xi\in V$ is fixed by the whole group $\G[A]$.
 \end{lemma}
 
 {\sc Proof of Lemma \ref{l:fixed-fixed}.} 
 Let $\xi\in V$. It is sufficient to show that $\xi$ is fixed 
by the subgroup $K[A]=\ker\kappa_A$. Denote by $K^{fin}[A]\subset K[A]$ the subgroup
consisting of transformations that are non-trivial only on finite
number of open balls $\B^o_v(p_j, \frh(v))$. This subgroup is dense
in $K[A]$ and so its sufficient to show that $\xi$ is fixed by 
$K^{fin}[A]$.

Thus, let $q\in K^{fin}[A]$. It is sufficient to prove that  
$$\la \rho(q)\xi,\eta\ra=\la \xi, \eta\ra \qquad\text{for each $\eta\in H$.}$$
 The right hand side equals to 
$$\la \xi, \eta\ra =\la \Pi \xi,\eta\ra=\la \xi,\Pi \eta\ra.$$

 The sequence $\Theta_j^{-1} q \Theta_j$
converges to the unit in $\G_\Lambda$ (since for each point $x\in\U_\Lambda$
the sequence $x\Theta_j^{-1} q \Theta_j$ is eventually constant).
Therefore unitary operators $\rho(\Theta_j^{-1} h \Theta_j)$ strongly converge
to the unit operator.

 We have
$$
\rho(q)= \rho(q)\rho( \Theta_j) \xi=\rho(\Theta_j)\rho(\Theta_j^{-1} q \Theta_j)\xi.
$$
So
\begin{multline}
\la \rho(q)\xi,\eta\ra-\la \xi, \eta\ra=\\
\bigl\la \rho(\Theta_j)\rho(\Theta_j^{-1} q \Theta_j)\xi,\eta\bigr\ra-
\bigl\la \xi,\Pi \eta\bigr\ra
= \bigl\la \rho(\Theta_j^{-1} q \Theta_j)\xi,\rho(\Theta_j)^{-1} \eta\bigr\ra-
\bigl\la \xi,\Pi \eta\bigr\ra
=\\=
\bigl\la\bigl( \rho(\Theta_j^{-1} q \Theta_j)-1\bigr)\xi, \rho(\Theta_j)^{-1} \eta\bigr\ra+\bigl\la \xi, (\rho(\Theta_j)^{-1}-\Pi)\eta\bigr\ra. 
\label{eq:long}
\end{multline}
By the strong convergence, we have
 $\|\bigl( \rho(\Theta_j^{-1} q \Theta_j)-1\bigr)\xi\|\to 0$
and therefore the first summand in right hand side of \eqref{eq:long} tends to 0. 
The sequence $(\rho(\Theta_j)-\Pi)$ tends to 0 weakly and therefore the second summand tends to 0.
 So the right hand side tends to 0. But the left hand side does not depend on $j$.
 So our expression equals  0.
\hfill $\square$

\sm

{\bf \punct  Multiplicativity.}

 \sm 
 
 {\sc Proof of Theorem \ref{th:multiplicativity}.}
By Proposition \ref{pr:Theta-move}, it is sufficient to prove the identity
$$
\lim_{j\to\infty} \Pi[X]\rho(g_1\Theta_j^Y g_2)\Pi[Z]=
\Pi[X]\rho(g_1) \Pi[Y] \rho(g_2) \Pi[Z],
$$ 
where $\lim$ denotes a limit in the weak operator topology. We have
\begin{equation*}
 \lim_{j\to\infty} \Pi[X]\rho(g_1)\rho(\Theta_j^Y)\rho( g_2)P[Z]=
 \Pi[X]\rho(g_1)\bigl( \lim_{j\to\infty}\rho(\Theta_j^Y)\bigr)\rho( g_2)P[Z],
\end{equation*}
and by Lemma \ref{l:Theta}, we get the desired expression.
\hfill $\square$

\sm 

{\sc Remark.}
Our way of proof of Theorem \ref{th:multiplicativity} based on Proposition \ref{pr:Theta-move}
and Lemma \ref{l:weakly} follows a scheme working 
in numerous situations, see \cite{Ner-book}, Sec. 8.5,  \cite{Ner-field}, Subsect. 3.5, \cite{Ner-faa}, Lemma 3.1, \cite{Ner-finite},
Theorem 1.9. On the other hand, our Theorem \ref{th:multiplicativity}
itself makes no sense without Proposition \ref{pr:density} formulated in the next subsection.
\hfill $\boxtimes$

\sm

{\bf \punct Density of fixed vectors.%
\label{ss:density}}
Let $\rho$ be a unitary representation of $\G_\Lambda$
in a Hilbert space $H$.

\begin{proposition}
\label{pr:density}
Let $X_1\subset X_2\subset X_3\subset \dots$ be finite subsets in $\U_\Lambda$ such that
$\cup X_j$ is dense in $\U_\Lambda$. Consider a unitary representation of $\G_\Lambda$ 
in a Hilbert space $H$. Denote by $H[X_j]$ the space of $\G[X_j]$-fixed vectors. Then
the subspace $\cup_{j=1}^\infty H[X_j]$ is dense in $H$.
\end{proposition}

This is an immediate corollary of the following statement, see \cite{Ner-book}, Proposition VIII.1.3.

\begin{proposition}
\label{pr:Ism}
Let $\rho$ be a unitary representation of a topological group  $G$ in a Hilbert space $W$. Let $G\supset G_1\supset G_2\supset \dots$ be a family of subgroups such that each neighbourhood of the unit in $G$ contains some subgroup
$G_n$. Denote by $W_k\subset W$ the subspace of $G_k$-fixed vectors. Then
$\cup_k W_k$ is dense in $W$.
\end{proposition}

To apply the latter proposition, we recall that  stabilizers $\G\{\B^c(a_j,r_j\}$ of finite collections $\{\B^c(a_j,r_j)\}_{j=1,\dots,m}$
 of closed balls form
a fundamental system of neighbourhoods of unit in $\G_\Lambda$. For sufficiently large
$N$ each ball contains a point of $X_N$ and therefore $\G[X_N]\subset \G\{\B^c(a_k,r_k)\}$.
\hfill $\square$

\sm

{\bf \punct Correspondence of subrepresentations.%
\label{ss:correspondence}}

\begin{proposition}
\label{pr:red}
{\rm a)} Let $\rho$ be a unitary representation of $\G_\Lambda$ in a Hilbert space $H$. Let $X\subset \U_\Lambda$ be a finite set
and the subspace $H[X]$ be cyclic%
\footnote{Let $\rho$ be a unitary representation of a group $G$ in a Hilbert
space $H$. Consider a collection of vectors $\{\xi_\alpha\}\in H$.
Its {\it cyclic span} is the closure of the space of finite linear combinations of all
vectors of form
$\rho(g)\xi_\alpha$. This is a minimal subrepresentation containing all vectors $\xi_\alpha$.
A system of vectors is {\it cyclic} if its cyclic span is the whole space.
}.
 For a subrepresentation $L\subset H$ consider the subspace $L\cap H[X]$.
The map $L\to L\cap H[X]$ is a bijection between subrepresentations of $\rho$ in $H$
and subrepresentations of the representation of $\End(X)$ in $H[X]$.

\sm

{\rm b)} Let $\rho_1$, $\rho_2$ be unitary representations of $\G_\Lambda$ in Hilbert spaces $H_1$, $H_2$.
Let the corresponding representations of $\End[X]$ be equivalent and 
$H_1[X]$, $H_2[X]$ be cyclic. Then $\rho_1$ and $\rho_2$ are equivalent.

\end{proposition}

This statement is a kind of a 'general nonsense' but the author does not know a reference that  precisely covers our case. 

\sm

{\sc Proof.}
a) First, we notice that the projection of $L$ to $H[X]$ coincides with
the intersection $L\cap H[X]$. Indeed, let $\eta\in L$. Then $\Pi[X]\eta$ is a limit of the sequence  of vectors $\rho(\Theta_j) \eta$
in the weak topology of the Hilbert space $H$.
Since $\rho(\Theta_j)\eta\in L$, then $\Pi[X]\eta\in L$.

Next, let us assume that $\Pi L=0$, i.e., $L\bot H[X]$. Let $\zeta\in H[X]$, $\eta\in L$, and $g\in\G_\Lambda$.
 Then 
$$
\la \rho(g) \zeta,\eta\ra=\la  \zeta,\rho(g^{-1})\eta\ra=0.
$$
Therefore,  $L$ is orthogonal to the $\G_\Lambda$-cyclic span of $H[X]$, and hence $L=0$.
So the map $L\mapsto L\cap H[X]$ is injective.

Further, let $V\subset H[X]$ be an $\End[X]$-invariant subspace.  Consider its $\G_\Lambda$-cyclic span $M$. Let $\eta\in V$.
Then
$$\Pi[X]\rho(g)\eta=\wh\rho_{X,X} (g)\eta\in V$$
and $\Pi[X]M=V$.

\sm 

b) Let us show that a representation $\wh\rho_{X,X}$ uniquely determines $\rho$ if $H[X]$ is cyclic.
Let $\eta$, $\zeta\in H[X]$.
Then 
$$
\la\rho(g_1)\eta,\rho(g_2)\zeta\ra_H=\la\eta,\rho(g_1^{-1}g_2)\zeta\ra_H=
\la\eta,\wh\rho_{X,X}(g_1^{-1}g_2)\zeta\ra_{H[X]},
$$
and the representation $\wh\rho_{X,X}$ determines the 'Gram matrix' of vectors $\rho(g)\eta$, where $\eta$ ranges in $H[X]$ and $g$ ranges in $\G_\Lambda$.
So the system of vectors $\{\rho(g)\eta\}$ in a Hilbert space is uniquely determined up to a unitary isometry.

Next, for $\eta$, $\zeta\in H[X]$, we have
\begin{multline*}
\la \rho(q)(\rho(g)\eta), \rho(p)\zeta\ra_H=
\la \eta, \rho(g^{-1}q^{-1}p)\zeta\ra_H=\\=\la \eta, \rho(g^{-1}q^{-1}p)\zeta\ra_H=\la \eta, \wh\rho_{X,X}(g^{-1}q^{-1}p)\zeta\ra_{H[X]}= 
\end{multline*}
Since the subspace $H[X]$ is cyclic,
a position of any $\rho(q)$-image of any vector $\rho(g)\eta$ is uniquely determined.
\hfill $\square$

\sm 

\begin{corollary}
\label{cor:splitting}
{\rm a)}
Denote by $X_{(1)}^1$ the one-point ultrametric space. For each $n>1$ we choose
a sequence $X_{(n)}^1$, $X_{(n)}^2$, \dots, which contains once each $n$-point ultrametric
space whose spectrum is contained in $\Lambda$ {\rm(}they are defined up to  isometries{\rm)}.
Then  any unitary representation $\rho$ of $\G_\Lambda$ splits into the following direct sum
of representations:
\begin{equation}
\label{eq:bigoplus}
V[X_{(1)}^1]\bigoplus
\Bigl(\bigoplus_{n=2}^ \infty \bigoplus_{\alpha=1}^\infty V[X_{(n)}^\alpha]\Bigr),
\end{equation}
where

\sm

---  the representation of $\G_\Lambda$ in $V[X_{(n)}^\alpha]$ has no nonzero vectors
fixed by subgroups $\G[X_{(m)}^\beta]$ for $m< n$  and no vectors fixed by subgroups 
$\G[X_{(n)}^\beta]$ with $\beta<\alpha$;

\sm

--- the subspace $H[X_{(n)}^\alpha]$ is $\G_\Lambda$-cyclic in $ V[X_{(n)}^\alpha]$. 

\sm

{\rm b)}
The summands of this splitting are invariant with respect to any $\G_\Lambda$-intertwining operator.
\end{corollary}

Notice that
$$
H[Xg]=\rho(g)H[X],
$$
therefore the definition of $V[X_{(n)}^\alpha]$ does not depend on the choice of the embedding $X_{(n)}^\alpha\to \U_\Lambda$.

\sm

{\sc Proof.} a) We take $H[X_{(1)}^1]$ and its $\G_\Lambda$-cyclic span $V[X_{(1)}^1]$. We split off this summand
and pass to $R_{(1)}=V[X_{(1)}^1]^\bot$. In this space there is no $\G[X_{(1)}^1]$-invariant vectors. 
Next consider the subspace of $\G[X_{(2)}^1]$-invariant vectors in $R_{(1)}$, consider its cyclic span
$V[X_{(2)}^1]$, and pass to $R_{(1)}\ominus V[X_{(2)}^1]$.
Repeating these arguments, 
we come to a decomposition
$$
V[X_{(1)}^1]\bigoplus\Bigl( \bigoplus_{\alpha=1}^\infty V[X_{(2)}^\alpha]\Bigr)\bigoplus R_{(2)},
$$ 
 where $R_{(2)}$ does not contain vectors fixed by  subgroups $\G[X_{(1)}^1]$ and $\G[X_{(2)}^\alpha]$.
 Repeating the same procedure we split off the sum \eqref{eq:bigoplus}. In the rest $R_{(\infty)}$, we have a representation
 of $\G_\Lambda$ that has no $\G[X_{(n)}^\alpha]$-fixed vectors for any $\alpha$, $n$. By Proposition \ref{pr:density},
 $R_{(\infty)}=0$.
 
 \sm
 
 b)
Next, consider an $\G_\Lambda$-intertwining operator in $H$.  By Lemma \ref{l:weakly}, projectors $\Pi[\cdot]$ are contained  
in the weak closure $\ov{\rho(\G_\Lambda)}$ and therefore their images are invariant with respect to all intertwining operators.
Therefore their cyclic spans also are invariant.
 \hfill $\square$

\section{Classification of representations\label{s:classification}}

\COUNTERS

Here we obtain the classification of unitary representations of the groups
$\G_\Lambda$. As we have seen in the previous section such representations are reduced to representations of 
semigroups $\End(X)$. Below, for an irreducible representation of $\G_\Lambda$  we consider a
minimal $X$ such that the subspace $H[X]$ of $\G[X]$-fixed vectors is nonzero.  
In this case we have an action of a certain quotient $\End^0(X)$ of $\End[X]$, this quotient 
has a  simple structure.
 We classify $*$-representations of $\End^0(X)$
(Lemma \ref{l:representations-semigroup}) and show that all such representations really
arise from unitary representations of $\G_\Lambda$ (Proposition \ref{pr:trivial}).
This gives us Theorem \ref{th:classification-main} about classification
of unitary representations of $\G_\Lambda$. Finally, we derive from this statement
weaker forms of the classification  (Theorems 
\ref{th:classification-main}$^\circ$, \ref{th:classification-main}$^{\circ\circ}$).

\sm

{\bf \punct Some remarks.}

\begin{observation}
The category $\cK_\Lambda$ is an {\it inverse category}.
\end{observation}

 This precisely means that for any morphism
$\frp:X\to Y$ we have
\begin{equation}
\frp\diamond\frp^*\diamond\frp=\frp.
\label{eq:ppp}
\end{equation}
This identity holds for any partial bijection. In Subsect. \ref{ss:partial-trees}  we realized our category $\cK_\Lambda$
as a subcategory of the category of partial bijections, so   identity \eqref{eq:ppp}
is obvious.

\sm 

{\sc Remark.}
The notion of inverse category (see \cite{Kas}, \cite{CL}) is a natural extension of the notion of inverse semigroup
(see, e.g., \cite{ClPr}). Below we do not use any general statements on this subject. However,
 the following proposition (which is common for arbitrary inverse semigroups),
is one of reasons of
 easy proof of classification of representations.
 \hfill $\boxtimes$

\begin{proposition}
Let $X$ be an object of the category $\cK_\Lambda$, i.e., a finite
subset in $\U_\Lambda$.

\sm 

{\rm a)} An element $\frp\in\End(X)$ is self-adjoint {\rm($\frp=\frp^*$)} if and only if it is an idempotent
{\rm ($\frp^2=\frp$)}.

\sm 

{\rm b)} Any two idempotents in $\End(X)$ commute.
\end{proposition}

 This is an obvious general statement for partial bijections,
see Subsect. \ref{ss:partial}.

\sm

{\bf\punct Topology on $\boldsymbol{\End(X)}$.}
By the construction, the countable set $\End(X)$ is a quotient space of $\G_\Lambda$, so it is equipped with the natural
quotient topology. An element $\frp\in\End(X)$ is determined by the matrix
$$
\Delta[\frp](x,y)=d\bigl(p_+(x), p_-(y)\bigr), \qquad \text{where $x$, $y\in X$.}
$$

\begin{lemma}
\label{l:topology}
A sequence  $\frp_j\in \End(X)$ converges to $\frp$ if the following two conditions are satisfied:

\sm

{\rm 1)} if $\Delta[\frp](x,y)=0$, then $\Delta[\frp_j](x,y)$ converges to $\Delta[\frp](x,y)$;

\sm

{\rm 2)} if $\Delta[\frp](x,y)\ne 0$, then $\Delta[\frp_j](x,y)=\Delta[\frp](x,y)$ for sufficiently large $j$;
\end{lemma}

{\bf Proof.} This is a corollary of the following observation: Let $x$, $x_j$, $y$ be points of an ultrametric
space. Let $x_j\to x$. If $d(x,y)>0$, then $d(x_j,y)=d(x,y)$ for sufficiently large $j$.
\hfill $\square$

 \sm

{\bf \punct Semigroups $\boldsymbol{\End^{00}(X)}$ of near-unit elements.%
\label{ss:near-unit}}
For $x\in X$ denote 
$$
d_x:=\min_{y\in X,\, y\ne x} d(x,y).
$$
If $X$ is a one-point space, then $d=\infty$.

\begin{figure}
$$
{\mathrm a)}\quad \epsfbox{trees.8}\qquad \mathrm{b)}\quad \epsfbox{trees.9}
$$
\caption{a) To Subsect. \ref{ss:near-unit}. A tree $T(P)$ corresponding to 
$\frp\in \End^{00}(X)$. We mark points of $p_+(X)$ by circles, points of $p_-(X)$
by squares. The tree $T(p_+(X))$ is drawn by solid lines, the forest $T(P)\setminus T(p_+(X))\subset T(p_-(X))$
by dashed lines, $\dom \pi= T(p_+(X))\cap T(p_-(X))$ by fat lines.
\newline
b) To the proof of Lemma \ref{l:decomposition}. A tree $T(P)$ for an idempotent, whish is not a near-unit.
We obtain $Y$ removing from $X$ the point marked by the arrow.
}

\label{fig:near}

\end{figure}

We say that an element $\frp=(P,p_+,p_-)$ of $\End(X)$ is {\it near-unit} if
for all $x$ we have 
$$d(p_+(x),p_-(x))< d_x$$
 In other words, $p_-(x)$ is the nearest to $p_+(x)$ point of $P$. 
On the language of partial bijections of trees $\pi:T(X)\to T(X)$
an element is near-unit if the following two conditions are satisfied:

\sm

$1^*$. $\dom \pi$ is a subtree of $T(X)$ that contains all actual vertices of $T(X)$ of non-zero
height together with their neighbourhoods (see Fig. \ref{fig:near}.a).

\sm

$2^*$. $\pi$ is identical on $\dom\pi$.

\sm

 Near-unit elements
$\frz\{\lambda_x\}$ are enumerated by collections
$\{\lambda_x\}_{x\in X}$, where
$$
\lambda_x:=d(p_+(x),p_-(x))\in \Lambda\cup 0.
$$

The following statement is obvious.

\begin{lemma}
\label{l:near-unit}
All near-unit elements of $\End(X)$ form a commutative semigroup,
\begin{equation}
\frz\{\lambda_x\}\frz\{\mu_x\}=\frz\{\max(\lambda_x,\mu_x)\},
\end{equation}
all its elements are self-adjoint idempotents,
\begin{equation}
\frz\{\lambda_x\}^2=\frz\{\lambda_x\}, \qquad \frz\{\lambda_x\}^*=
\frz\{\lambda_x\}.
\end{equation}
\end{lemma}

We denote the semigroup of all near-unit elements by $\End^{00}(X)$.

\sm

Fix $d$, which is positive or $\infty$. Consider the semigroup $\cZ_d(\Lambda)$ consisting of elements
$\frz_\lambda$, where $\lambda$ ranges in $(\Lambda \cup 0)\cap[0,d)$,
satisfying the relations
$$
\frz_\lambda^*=\frz_\lambda,
\qquad
\frz_\lambda^2=\frz_\lambda,\qquad \frz_\lambda \frz_{\lambda'}
=\frz_{\max(\lambda,\lambda')}.
$$
We also assume that if $\lambda_j\to 0$, then the sequence $\frz_{\lambda_j}$ tends to $\frz_0$,
and assume there are no other convergent and not eventually constant sequences.
Then Lemma \ref{l:near-unit} can be reformulated in the form:

\newtheorem*{lemma1}{Lemma \ref{l:near-unit}$^{\circ}$}

\begin{lemma1}
$
\End^{00}(X)\simeq\prod_{x\in X}\cZ_{d_x}(\Lambda).
$
\end{lemma1}

\sm

{\bf \punct Semigroups $\boldsymbol{\End^{0}(X)}$ of  near-automorphisms.%
\label{ss:near-aut}}
We say that an element $P\in \End(X)$ is a {\it near-automorphism}, if there is an isometry $\kappa: X\to X$ such that for each point 
$x\in X$ the nearest point to $p_+(x)$ in $P$  is $p_-(\kappa(x))$. We denote the set of all near-automorphisms by $\End^{0}(X)$.

On the language of partial bijection of trees $T(X)$ this means that the condition
$1^*$ of the previous subsection  is fulfilled.

The following  statements is obvious.

\begin{lemma}
{\rm a)} The set of near-automorphisms of $X$ is a semigroup. The group of invertible elements is isomorphic to the group $\Isom(X)$.

\sm 

{\rm b)}
Any element of $\End^{0}(X)$ is a product of an automorphism and an element $\frz\{\lambda_x\}$.  

\sm

{\rm c)} For any  $\kappa\in\Isom(X)$, we have
$
\kappa^{-1}\frz\{\lambda_x\}\kappa=\frz\{\lambda_{\kappa(x)}\}.
$
 \end{lemma}

{\bf \punct Reduction of the problem of classification to semigroups $\boldsymbol{\End^0(X)}$.}

\begin{lemma}
\label{l:PP}
Let $\frp$ be an endomorphisms of $X$, which is not near-automorphism. Then the idempotent $\frp^*\diamond\frp$ is not a near-automorphism.  
\end{lemma}

{\sc Proof.} Let $\pi$ be the corresponding partial isomorphism $\pi:T(X)\to T(X)$.
 Then  $\dom \pi$ does not satisfy condition $1^*$ of Subsect. \ref{ss:near-unit}.
 But $\dom \pi^*\pi\subset \dom\pi$ and it also does not satisfy condition $1^*$.
\hfill $\square$

\begin{lemma}
\label{l:decomposition}
Let $\frq\in \End(X)$ be an idempotent, which is not a near-auto\-mor\-phism.
 Then there is a proper subspace $Y\subset X$  and a morphism
 $\frt:X\to Y$ such that
 $$
 \frt^*\diamond \frt=\frr.
 $$
\end{lemma}

{\sc Proof.} Let $\kappa$ be the partial isomorphism $T(X)\to T(X)$ corresponding to $\frq$.
 Since $\frq$ is not a near-automorphism, there
 is an actual vertex $v$ of $T(X)$ and an edge $[v,w]$ of $T(X)$ with upper vertex $v$ such
that $[v,w]\cap \dom \kappa$ is empty or consists of the point $v$.
Consider the set $Z$ of all points of $X$ under $w$ (if $w\in X$, we take $Z:=\{w\}$).
Set $Y:=X\setminus Z$, and define
a partial isomorphism $\tau:T(X)\to T(Y)$ such that $\dom\tau=\dom \kappa$ and 
$\tau$ is identical on $\dom \tau$, see Fig. \ref{fig:near}.b.
Then $\kappa=\tau^*\tau$. 
\hfill $\square$

\sm


\begin{lemma}
\label{l:reduction}
Consider a unitary representation $\rho$ of the group $\G_\Lambda$
in a Hilbert space $H$.
Let $X$ be a minimal element of the set of all $A$ such that $H[A]\ne0$.
Then for any $\frp\in \End(X)\setminus \End^{0}(X)$ we have $\wh\rho_{X,X}(\frp)=0$.
\end{lemma}

{\sc Proof.} By Lemma \ref{l:PP}, $\frp^*\diamond\frp\notin \End^{0}(X)$.
Applying Lemma \ref{l:decomposition}, we
represent $\frp^*\diamond\frp=\frt^*\diamond\frt$, where $\frt\in \End(X,Y)$
with
 $Y \subsetneqq X$. Since $H[Y]=0$, we have $\wh\rho_{X,Y}(\frq)=0$, therefore
$\wh\rho_{X,X}(\frp^*\diamond\frp)=0$ and $\wh\rho_{X,X}(\frp)=0$.
\hfill $\square$

\sm

{\bf \punct Description of  representations of semigroups $\boldsymbol{\End^{0}(X)}$.}
First, let us describe $*$-representations of the semigroup $\cZ_d(\Lambda)$.
Since it is commutative, irreducible representations are one-dimensional.
Clearly, elements $\frz_\lambda$ act in a one-dimensional
space by multiplications by $0$ or $1$. So the set
$(\Lambda \cup 0)\cap[0,d)$ splits into two subsets,
$L_0$, where $\frz_\lambda=0$, and $L_1$, where $\frz_\lambda=1$.
Any element  of the first set is larger than any element of the second set.
So we have a section of an  ordered set  $(\Lambda \cup 0)\cap[0,d)$ in the Dedekind sense. Let us formulate this more precisely.

Irreducible representations $ \chi_{\epsilon, r}$  of the semigroup $\cZ_d$ are enumerated by parameters
$(\epsilon,r)$, where $\epsilon$ is '$c$' or '$o$', and $r$ is real. They satisfy
the following restrictions:
 if $\epsilon$ is  $c$, then $r\in \Lambda^{+0}\cap [0,d)$; 
if $\epsilon$ is $o$, then $r$ is an element of $\Lambda^\circ\cap (0,d]$
(if $d=\infty$, then we allow $r=\infty$).

%

\sm

We define the function $\zeta(\lambda;\epsilon,r)$
by
$$
\zeta(\lambda;c,r)=\begin{cases}
1, \quad \text{if $r\ge \lambda$};
\\
0, \quad \text{otherwise},
\end{cases}
\qquad
\zeta(\lambda;o,r)=\begin{cases}
1, \quad \text{if $r> \lambda$};
\\
0, \quad \text{otherwise}.
\end{cases}
$$
Then a one-dimensional $*$-representation $\chi_{\epsilon, r}$ is given by
the formula
$$
\chi_{\epsilon, r}(\frz_\lambda)=\zeta(\lambda;\epsilon,r). 
$$

{\sc Remarks.}
This is clear, we only notice that:

\sm

---  Our conditions for pairs $(\epsilon,r)$ forbid the
representation $\rho(\frz_0)=1$, $\rho(\frz_\lambda)=0$ for $\lambda>0$ if
$0$ is not a limit point of $\Lambda$. This representation is discontinuous.

\sm

--- Our conditions forbid zero representation,
$\nu(\frz_\lambda)=0$ for all $\lambda$ since for a representations $\wh\rho_{X,X}$ of $\End(X)$ arising
from a representations $\rho$ of $\G_\Lambda$ we have  $\wh\rho_{X,X}(1)=1$ and the unit of 
$\End(X)$ is $\frz\{0,\dots,0\}$.
\hfill $\boxtimes$

\sm

\sm 

A {\it labeling} $\{\epsilon_x,r_x\}_{x\in X}$ in the sense of Subsect. \ref{ss:classification} means
that we assign an admissible pair $\{\epsilon_x,r_x\}$ to each point $x\in X$. 
So labelings enumerate  irreducible representations of $\End^{00}(X)$.

For a labeling $\{\epsilon_x^\centerdot,r_x^\centerdot\}_{x\in X}$,
we denote by $\Gamma\{\epsilon_x^\centerdot,r_x^\centerdot\}\subset \Isom(X)$
the subgroup preserving it.
Consider the homogeneous space 
$$\cX=\cX\{\epsilon_x^\centerdot,r_x^\centerdot\}:=
\Gamma\{\epsilon_x^\centerdot,r_x^\centerdot\}\backslash \Isom(X).$$
Points of this space are all labellings of $T(X)$, which are equivalent to $\{\epsilon_x^\centerdot,r_x^\centerdot\}$ under the action 
of $\Isom(X)$. 
Fix  an irreducible representation $\nu$ of 
$\Gamma\{\epsilon_x^\centerdot,r_x^\centerdot\}$
in a Euclidean space $V$.
These data, $\{\epsilon_x^\centerdot,r_x^\centerdot\}$ and $\nu$, determine an irreducible
representation of $\End^0(X)$ in the following way:

Consider the (finite-dimensional) Euclidean space 
$\ell^2(\cX\{\epsilon_x^\centerdot,r_x^\centerdot\},V)$.
Elements $\frz\{\lambda_x\}$ act in this space by multiplications:
\begin{equation}
\frz\{\lambda_x\} F\bigl(\bigr)=
\prod_{x\in X} \zeta(\lambda_x;\epsilon_x,r_x)\,
F\bigl(\{\epsilon_x,r_x\}\bigr).
\label{eq:zeta}
\end{equation}
A representation of 
$\Isom(X)$ in $\ell^2(\cX\{\epsilon_x^\centerdot,r_x^\centerdot\},V)$
is the representation induced from the representation
$\nu$ of $\Gamma\{\epsilon_x^\centerdot,r_x^\centerdot\}$. This means
 (see the definition of  induced representations
in Subsect. \ref{ss:induced}) that
the group $\Isom(X)$ acts by operators of the form
\begin{equation}
\tau(g)F\bigl(\{\epsilon_x,r_x\}\bigr)=
A\bigl(\{\epsilon_x,r_x\},g\bigr) F\bigl(\{\epsilon_{xg},r_{xg}\}\bigr),
\label{eqLtau} 
\end{equation}
where $A(\cdot,\cdot)$ is a function on $\cX\times \Isom(X)$ taking values
in the unitary group of $V$,
and for $h\in \Gamma\{\epsilon_x^\centerdot,r_x^\centerdot\}$
we have 
$$
A(\{\epsilon_x^\centerdot,r_x^\centerdot\},h)=\nu(h).
$$

\begin{lemma} 
\label{l:representations-semigroup}
{\rm a)} Any irreducible $*$-representation $\rho$ of the semigroup $\End^{0}(X)$
such that $\rho(1)\ne 0$
has this form.

\sm

{\rm b)} Any  $*$-representation of the semigroup $\End^0(X)$
is a direct integral of irreducible representations.

\sm

{\rm c)} The semigroups $\End^0(X)$ have type $I$.
\end{lemma}

{\sc Proof.} Let us show that our representation is irreducible. Consider
an intertwining operator $R$. The family of functions
$$
\Phi_{\{\lambda_x\}}\bigl(\{\epsilon_x,r_x\}\bigr)
:= \prod \zeta(\lambda_x;\epsilon_x,r_x)
$$
separate points of the $\Isom(X)$-orbit on the set of labelings,
 therefore the operator $R$ acts fiber-wise.
Since $\nu$ is irreducible, $R$ is scalar in the fiber over 
$\{\epsilon_x^\centerdot,r_x^\centerdot\}\in\cX$.
 By the  homogeneity of $\cX$, the operator $R$ is the same scalar in each fiber.

\sm

a) Consider a  semigroup algebra of $\End^0(X)$, i.e., the algebra of finite formal
linear combinations $\sum_{\fru\in \End^0(X)} c_\fru \fru$.
For each collection $\{\lambda_x\}_{x\in X}$ consider the sum
$$
\Delta_{\{\lambda_x\}}:=\sum_{g\in \Isom(X)} \frz_{\{\lambda_{xg}\}}.
$$
For $g\in\Isom(X)$ we have
\begin{equation}
g^{-1}\frz_{\{\lambda_x\}}g=\frz_{\{\lambda_{xg}\}}. 
\label{eq:conjugation}
\end{equation}
Therefore elements $\cZ_{\{\lambda_x\}}$ are central elements
of the semigroup algebra of $\End^0(X)$ and hence they act 
in irreducible representations by scalar operators.

Restrict an irreducible representation $\rho$  of    the semigroup
$\End^{0}$ to $\End^{00}$ and decompose the restriction
into an integral  of irreducible representations. The operators
$\rho(\Delta_{\{\lambda_x\}})$ act in this integral by multiplications
by certain $\Isom(X)$-invariant functions. It is easy to see that
such functions separate $\Isom(X)$-orbits on the set of labelings. So a representation
of $\End^{00}$ can be irreducible only if the direct integral is taken over
 one $\Isom(X)$-orbit on the set of labelings.
 
  So our representation
 is realized in a space of vector-valued functions supported by some orbit
 $\Gamma\{\epsilon_x^\centerdot,r_x^\centerdot\}\backslash \Isom(X)$ and the semigroup $\End^{00}(X)$
 acts by the operators \eqref{eq:zeta}. In other words
 our space is a direct sum of summands enumerated by points of the 
 orbit. By \eqref{eq:conjugation}, the group $\Isom(X)$ acts permuting
 summands, so we have an induced representation.

%

 \sm 
 
 b), c) Consider a $*$-representation $\rho$ of $\End(X)$
 and the corresponding $C^*$-algebra $\cC$, i.e., we consider
 the linear span of operators of the representation and its closure
 with respect  to norm operator topology. Since the semigroup $\End(X)$ is countable, we  get a separable
 $C^*$-algebra. Restricting any representation of $\cC$ to
 $\End(X)$ we get a $*$-representation of $\End(X)$. So all irreducible representations
 of $\cC$ are finite-dimensional. Hence (see Dixmier \cite{Dix}, 9.1) $\cC$
 is a GCR-algebra and it has type I. Since our algebra is separable,
any its representation can be decomposed to a direct integral
of irreducible representations
(see \cite{Dix}, 8.3). 
\hfill
$\square$

\sm 

{\bf \punct Induced representations of $\boldsymbol{\G_\Lambda}$.}
All irreducible representations of the semigroup
$\End^0(X)$ described in Lemma \ref{l:representations-semigroup} actually arise from unitary representations
of the group $\G_\Lambda$. Precisely, for a labeling 
$\{\epsilon_x,r_x\}_{x\in X}$ we assign a collection of balls
$\{\B^\epsilon(x,r_x)\}_{x\in X}$.
Consider the group $\G^\bullet\{\B^{\epsilon_x}(x,r_x)\}$ consisting of all $g\in\G_\Lambda$ sending the set
$\cup_x \B^{\epsilon_x}(x,r_x)$ to itself and its subgroup $\G\{\B^{\epsilon_x}(x,r_x)\}$ consisting
of transformations sending each ball $\B^{\epsilon_x}(x,r_x)$ to itself.
 Then the quotient
$\G^\bullet\{\B^{\epsilon_x}(x,r_x)\}/\G\{\B^{\epsilon_x}(x,r_x)\}$ is finite and is isomorphic to the subgroup  $\Gamma\{\epsilon_x,r_x\}$ of $\Isom(X)$
preserving the labeling.

The homogeneous space $\G\{\B^{\epsilon_x}(x,r_x)\}\backslash \G_\Lambda$ is countable, it is the space of all collections
of balls isometric to the initial collection $\{\B^{\epsilon_x}(x,r_x)\}_{x\in X}$. So we can apply
the construction of induction described in Subsect. \ref{ss:induced}. We consider an irreducible representation 
$\nu$ of $\Gamma\{\epsilon_x,r_x\}$, regard it as a representation of $\G^\bullet\{\B^{\epsilon_x}(x,r_x)\}$ and take the induced representation of $\G_\Lambda$,
\begin{equation}
\label{eq:ind-nu}
\rho=\Ind(\nu):=\Ind_{\G^\bullet\{\B^{\epsilon_x}(x,r_x)\}}^{\G_\Lambda}(\nu).
\end{equation}

 \begin{proposition}
 \label{pr:trivial}
{\rm a)} Any representation {\rm \eqref{eq:ind-nu}} of $\G_\Lambda$  is irreducible.

\sm 

{\rm b)} 
Representations corresponding to triples $(X, \{\epsilon_x,r_x\}, \nu)$, 
$(X', \{\epsilon_{x'},r_{x'}\}, \nu')$ can be equivalent only by a trivial reason, i.e., if there is
$g\in\G_\Lambda$ such that $X'=Xg$, 
$\{\epsilon_{x'},r_{x'}\}=\{\epsilon_{xg},r_{xg}\}$, $\nu'(h)=\nu(g^{-1}hg)$.

  \sm
  
 {\rm c)} The set $X$ is a  minimal subset in $\U_\Lambda$, for which the space of $\G[X]$-fixed vectors is nonzero.
 
 \sm
 
 {\rm d)}  The representation of $\End(X)$ in the space of $\G[X]$-fixed vectors coincides with the representation
 described above in Lemma {\rm \ref{l:representations-semigroup}}.
\end{proposition}
  
This statement must be considered as obvious for experts in representation theory. To be complete (and to explain  for non-experts the origin of the phenomenon), we present a proof.
We need the following lemma, which is a rephrasing of Mackey \cite{Mack0}, Theorem 3, or Corwin \cite{Cor}
in a convenient for us form. 

\begin{lemma}
\label{l:mackey}
 Consider a separable topological group $G$ and its open subgroups $L$, $M$.
 Let $\tau$ be a unitary finite-dimensional representation of $L$ in a Hilbert space $V$, $\sigma$
 be a unitary finite dimensional representation of $M$ in a Hilbert space $W$. 
 
 \sm
 
 {\rm a)} Let all orbits of $L$ on $M\backslash G$ be infinite. 
 Then the induced representations $\rho=\Ind_L^G(\tau)$ and $\tau=\Ind_M^G(\sigma)$
 are disjoint%
 \footnote{Two unitary representations are {\it disjoint} if they have no equivalent subrepresentations.  Equivalently, there are no intertwining 
 operators between them.}.
 
 \sm

 {\rm b)} Let $\nu$ be irreducible and all orbits of $L$ on $L\backslash G$ except the initial point be infinite.
 Then $\Ind_L^G(\nu)$ is irreducible.
 \end{lemma}
 
  {\sc Proof.} a) Representations $\rho$ and $\tau$ are realized 
  in the spaces $\ell^2(L\backslash G,V)$, $\ell^2(M\backslash G,W)$
 by formulas of the form
 \begin{align*}
 \rho(g) \phi(x)=A(x,g) \,\phi(xg), \qquad \tau(g)\psi(y)=B(y,g)\,\psi(yg),
 \end{align*}
 where $A(x,g)$, $B(y,g)$ are unitary cocycles.
  Any bounded
 operator 
 $$S:\ell^2(L\backslash G,V)\to \ell^2(M\backslash G,W)$$
 can be written in a form
 $$
 S \phi(y)=\sum_{x\in L\backslash G} \Theta(y,x) \phi(x),
 $$
 where $\Theta(y,x)$ is a function on $(M\backslash G)\times(L\backslash G)$
 taking values in the space of linear operators $V\to W$.
 Let $S$ be intertwining, i.e., $\tau(g)S=S\rho(g)$.
 The left hand side is
 $$
 \tau(g) S \phi(y)=\sum_{x\in L\backslash G} B(y,g)S(yg,x) f(x),
 $$
 the right hand side is
 $$
 S \rho(g) \phi(y)=\sum_{z\in L\backslash G} S(y,z) A(z,g)f(zg)=
 \sum_{x\in L\backslash G} S(y,xg^{-1}) A(x,g^{-1})^{-1}f(x).
 $$
 (we changed the summation index, $x=zg$).  Therefore the kernel $S$ satisfies the equation
 $$
 B(y,g)S(yg,x)=S(y,xg^{-1})A(x,g^{-1})^{-1}.
 $$
By $L^x\subset G$ we  denote the stabilizer of $x\in L\setminus G$.
 Let $h\in L^x$. Then 
 $$
 S(yh,x)=B(y,h)^{-1}S(y,x)A(x,h^{-1})^{-1}.
 $$
 
Denote by $p\delta_z\in \ell^2(L\backslash G,V)$ a function, which is zero outside a point $z$ and is $p$ at $z$.
Fix an orthonormal basis $v_1$, \dots, $v_n\in V$.
Our operator $S$ sends a function $v_j\delta_x$
to the function 
$$
F_j=\sum_{u\in y L^x} S(u,x)v_j \delta_u+ \sum_{u\notin y L^x} S(u,x)v_j \delta_u.
$$
We have 
\begin{equation}
\|F_j\|^2_{\ell^2}\ge \sum_{u\in y L^x} \|S(u,x)v_j\|^2_W.
\label{eq:norm-l}
\end{equation}
 Notice that an operator norm $\|S(u,x)\|$ does not depend
  on $u$ if $u$ lies in one 
$L^x$-orbit on $M\backslash G$.  Next, we notice that for any operator $T:V\to W$
there is a basis element $v_j$ such that $\|Tv_j\|\ge \frac 1n \|T\|$, where $n=\dim V$. 
Therefore if $S(y,x)\ne 0$  at least for one $j$, then the sum in the right hand side of \eqref{eq:norm-l} is infinite.
Therefore $S(y,x)=0$ for all $y$, $x$.

\sm
 
The statement  b) follows from the same considerations. 
\hfill $\square$

\sm

{\sc Proof of Proposition \ref{pr:trivial}.}
a), b) In Lemma \ref{l:mackey} we take $G=\G_\Lambda$, and 
$$L=\G^\bullet\bigl(\{\B^{\epsilon_x}(x,r_x)_{x\in X}\}\bigr),\qquad 
M=\G^\bullet\bigl(\{\B^{\epsilon_{x'}}(x',r_{x'})\}_{x'\in X'}\}\bigr).$$

Consider two balls $\B^\epsilon(x,r)$ and $\B_{\epsilon'}(x',r')$, and consider
their stabilizers $\Gamma$, $\Gamma'\in \G_\Lambda$. Then $\Gamma$-orbit on
$\Gamma'\backslash \G_\Lambda$ can be finite only if 
$\B_{\epsilon'}(x',r')\supset \B_{\epsilon}(xg,r)$ for some $g\in \G_\Lambda$.
Moreover, this finite orbit is unique and consists of one point.

Similarly, an $L$-orbit on $M\backslash \G_\Lambda$ can be finite,
only if there exists $g\in \G_\Lambda$ such that each ball 
$\B_{\epsilon'}(x',r')$ contains a ball $\B_{\epsilon}(xg,r)$.
This is possible, but we can interchange $L$ and $M$. So 
a nonzero intertwining operator can exists only if 
the collections $\{\B^{\epsilon_{x}}(x,r_{x}$, $\{\B^{\epsilon_{x'}}(x',r_{x'}$ are isometric.

\sm

c) Consider an induced representation of $\G_\Lambda$ in 
$\ell^2\bigl(\G^\bullet\{\B^{\epsilon_x}(x,r_x)\}\backslash \G_\Lambda,V\bigr).$
Let $\phi(\zeta)$ be an element of $\ell^2$ 
invariant with respect to $\G[Y]$. Then the function $p(\zeta):=\|\phi(\zeta)\|_V$ is constant 
on orbits of $\G[Y]$. So it can be nonzero only on finite orbits.
A union of balls $\cup_x \B^{\epsilon_x}(xg,r_x)$ can have a finite $\G[Y]$-orbit only
if each ball $\B^{\epsilon_x}(xg,r_x)$ contains 
 contains a point of $Y$. Therefore if the space of fixed functions is nontrivial, then
$Y$ contains a subset isometric to $X$. 

\sm

 d) Consider the  space  $H=\ell^2(\G^\bullet\{\B^{\epsilon^\centerdot_x}(x,r^\centerdot_x)\}\backslash \G_\Lambda,V)$ 
and the representation $\rho$ of $\G_\Lambda$ in this space. 
 The subspace $H[X]$ consists of $V$-valued functions, which are zero outside
 the $\Isom(X)$-orbit of $\cup_{x\in X}\B^{\epsilon^\centerdot_x}(x,r^\centerdot_x)$, denote this orbit by $\Omega$.
 Denote by $V\{\epsilon_x,r_x\}\subset H[X]$ the fiber over $\cup_{x\in X}\B^{\epsilon_x}(x,r_x)\in\Omega$,
 $$
 H[X]=\bigoplus_{\cup_x\B^{\epsilon_x}(x,r_x)\,\in\, \Omega} V\{\epsilon_x,r_x\}. 
 $$
Let us describe the action of the semigroup $\End^{00}(X)$ in $H[X]$. Consider elements
$\frz\{\lambda_x\}\in\G_\Lambda$ such that $\lambda_x=d(xh,x)<d_x$ for all $x\in X$.
The corresponding elements of $\End(X)$ are $\frz\{\lambda_x\}$. If for some ball
$\B^{\epsilon_x}(x,r_x)$ we have $h\{\lambda_x\}\bigl(\B^{\epsilon_x}(x,r_x)\bigr) \ne \B^{\epsilon_x}(x,r_x)$,
then $\Pi[X]\rho( h\{\lambda_x\})\bigl( V\{\epsilon_x,r_x\}\bigr)=0$, i.e., $\wh \rho_{X,X}(\frz\{\lambda_x\})$ sends  this fiber to 0. If for all $x\in X$ we have 
$h\{\lambda_x\}\bigl(\B^{\epsilon_x}(x,r_x)\bigr)=  \B^{\epsilon_x}(x,r_x)$,
then $\rho(h\{\lambda_x\})$ induces a unitary operator in the fiber $V\{\epsilon_x,r_x\}$, and therefore
$\wh \rho_{X,X}(\frz\{\lambda_x\})$ induces the same unitary operator in the fiber. Since $\wh \rho_{X,X}(\frz\{\lambda_x\})$ is a projector,
it acts in this fiber as 1.

We omit a watching of the action of $\Isom(X)\subset \End^0(X)$ in $H[X]$. 
\hfill $\square$

\sm

{\bf\punct End of proof the classification theorem.} 
It remains to note, that all necessary elements of the classification theorem are proved.

\sm

{\sc Proof of Theorem \ref{th:classification-main}.} By the Multiplicativity Theorem \ref{th:multiplicativity},
an irreducible representation $\rho$ of $\G_\Lambda$ in a Hilbert space $H$ generates a $*$-representation of the category $\cK_\Lambda$.
We take a minimal subset $X\subset \U_\Lambda$, for which $H[X]$ is nonzero. By Lemma \ref{l:reduction}, the representation
of $\End(X)$ in this space is reduced to a representation of $\End^0(X)$. By Proposition \ref{pr:red}, the representation
of $\End^0(X)$ is irreducible and determines uniquely the  representation $\rho$. 
Classification of irreducible irreducible representations of $\G_\Lambda$ now  follows from Lemma \ref{l:representations-semigroup}.a and Proposition
 \ref{pr:trivial}.d. 
 
For a decomposition of a reducible representation $\rho$ of $\G_\Lambda$ into a direct integral, we apply Corollary \ref{cor:splitting}.a
and get a decomposition \eqref{eq:bigoplus}. After this, for each summand we decompose to the integral  the representation of $\End(X^\alpha_{(n)})$
in $\G[X^\alpha_{(n)}]$-fixed vectors, this is possible by Lemma \ref{l:representations-semigroup}.b. By Proposition \ref{pr:red},
we come to a decomposition of $\rho$ into a direct integral.

It remains to verify that any unitary representations $\rho$ of $\G_\Lambda$ has type $I$. It is sufficient to verify
(see, e.g, \cite{Dix}, Addendum A.52)
that the von Neumann algebra of all intertwining operators has type $I$. Such operators preserve the decomposition
\eqref{eq:bigoplus}; by Proposition \ref{pr:red}, the algebra of intertwining operators for $\G_\Lambda$ in each $V[X_{(n)}^\alpha]$
coincides with the algebra of intertwining operators for $\End^0([X_{(n)}^\alpha)$ in the space of $\G[X_{(n)}^\alpha]$-fixed vectors.
   By  Lemma \ref{l:reduction}.c, this algebra has type $I$.
   \hfill $\square$
   
   \sm
   
 {\sc Proof of Theorem \ref{th:classification-main}$^{\circ\circ}$.}  We apply Lemma \ref{l:stages} to
 $G=\G_\Lambda$, $N=\G^\bullet\{\B^{\epsilon_x}(x,r_x)\}$ , $H=\G\{\B^{\epsilon_x}(x,r_x)\}$.
 \hfill $\square$
 
 \sm 

 {\sc Proof of Theorem \ref{th:classification-main}$^{\circ}$.} 
We fix $X$ and a collection $\{\B^{\epsilon_x}(x,r_x)\}$. Keeping the notation in this theorem, we write
$$
\bigotimes_{x\in X} \ell^2(\Xi_{\epsilon_x,r_x})=
\ell^2\Bigl(\prod_{x\in X} \Xi_{\epsilon_x,r_x} \Bigr).
$$
Let $\cO^\alpha$ be orbits of $\G_\Lambda$ on the countable set $\prod \Xi_{\epsilon_x,r_x}$.
Then our space $\ell^2$ splits into a direct sum
$
\oplus_{\alpha} \ell^2(\cO^\alpha).
$
At least one of such orbits is $\G\{\B^{\epsilon_x}(x,r_x)\}\backslash \G_\Lambda$,
and we apply Theorem \ref{th:classification-main}$^{\circ}$
\hfill $\square$

\sm

{\bf \punct An additional remark.} Apparently, all homomorphisms
$\G_\Lambda\to S_\infty$ admit a classification in the following sense%
\footnote{For several statements of this kind, see Tsankov \cite{Tsa}.}:

\begin{conjecture}
\label{conj:1}
Any open subgroup  in $\G_\Lambda$ is an intermediate subgroup between 
some
$\G^\bullet\{\B^{\epsilon_x}(x,r_x)\}$ and $\G\{\B^{\epsilon_x}(x,r_x)\}$.
\end{conjecture}

The argument for this conjecture is the following: an action
of $\G_\Lambda$ on a countable set must generate unitary representation
and they must appear somewhere in our classification. On the other hand open subgroups 
correspond to transitive actions on discrete spaces.

\section{Universal compactifications\label{s:university}}

\COUNTERS
 
 {\bf \punct The definition of the semigroups $\boldsymbol{\GammA_\Lambda}$.}
 Fix $\Lambda$. As in Subsect. \ref{ss:woolly}, denote
 $$
 \wt T(\U_\Lambda):=\begin{cases}
 T_{>0}(\U_\Lambda)\coprod \bigl(\,\coprod\limits_{v\in\ver T(\U_\Lambda)} \germ^\downarrow(v)\bigr),
 & \text{if 0 is limit point of $\Lambda$;}
 \\
  T(\U_\Lambda)\coprod \bigl(\,\coprod\limits_{v\in\ver T(\U_\Lambda)} \germ^\downarrow(v)\bigr),&
  \text{otherwise.}  
 \end{cases}
 $$
An element of the  $\GammA_\Lambda$ is a partial bijection 
$$
\pi: \wt T(\U_\Lambda)\to  \wt T(\U_\Lambda),
$$
whose domain is a woolly subtree $J$ in $T(\U_\Lambda)$ and 
which sends $J$ isomorphically (with a preservation of the height)
to a woolly subtree.

\sm

{\bf\punct Topology on $\boldsymbol{\GammA_\Lambda}$.}
Let $\Omega$ be a set (of arbitrary cardinality).
Consider the semigroup ${\End\vphantom{q}}_\PB\,(\Omega)$ of partial bijections
$\Omega\to\Omega$. Let $S$ be a finite set, consider  collections
of points $\{\omega_s\}_{s\in S}$, $\{\omega_s'\}_{s\in S}$ in $\Omega$.
Consider also a collection $\{\delta_s\}_{s\in S}$, where each $\delta_s$ is 0 or 1.
We define a subset $\cO(\{\omega_s,\omega'_s\}, \{\delta_s\})\subset {\End\vphantom{q}}_\PB\,(\Omega)$
as the set of all partial bijections $\psi$ satisfying the conditions:

\sm

--- $\omega_s'=\psi(\omega_s)$ if $\delta=1$;

\sm 

--- $\omega_s'$ is not $\psi$-image of $\omega$ if $\delta_j=0$ (i.e. $\omega_s\notin\dom\psi$, or 
$\psi(\omega_s)\ne \omega'_s$).

\sm

We define the topology on ${\End\vphantom{q}}_\PB\,(\Omega)$, whose base consists of  all possible sets 
$\cO(\{\omega_s,\omega'_s\}, \{\delta_s\})$. Then ${\End\vphantom{q}}_\PB\,(\Omega)$ becomes
a compact Hausdorff topological space, and the product in this semigroup
is separately continuous, see, e.g., \cite{mnogo}, Theorem 5.12.

\sm

{\sc Remark.} A graph of a partial bijection is a subset in $\Omega\times \Omega$. The set of all subsets
of $\Omega\times \Omega$ is in a one-to-one correspondence with the set of functions on $\Omega\times \Omega$
with values $0$, $1$. So  ${\End\vphantom{q}}_\PB\,(\Omega)$ can be identified with a subset
in the direct product of sets $\{0,1\}$, in which factors are enumerated by points of $\Omega\times\Omega$.
Our topology on ${\End\vphantom{q}}_\PB\,(\Omega)$ is induced from the topology of the direct product,
and compactness is provided by the Tychonoff theorem.
\hfill $\boxtimes$

\sm

We equip the semigroup $\GammA_\Lambda$ with the topology induced from the semigroup
${\End\vphantom{q}}_\PB\,\bigl(\wt T(\U_\Lambda)\bigr)$.

\begin{proposition}
{\rm a)} The semigroup $\GammA_\Lambda$ is compact.

\sm 

{\rm b)} The product in $\GammA_\Lambda$ is separately continuous.
\end{proposition}

{\sc Proof.} 
We must only verify that $\GammA_\Lambda$ is closed in the semigroup of partial bijections
$\wt T(\U_\Lambda))\to \wt T(\U_\Lambda))$.
 Let $\xi$, $\xi'\in \wt T(\U_\Lambda)$
correspond to perfect balls $\B^\epsilon(x,r)$, $\B^\epsilon(x',r)$, let 
$\zeta$, $\zeta'$ correspond to larger perfect balls $\B^\delta(x,\rho)$, $\B^\delta(x',\rho)$. Then  the condition
\begin{equation}
\text{'if $\pi$ sends $\xi$ to $\xi'$, then it sends $\zeta_1$ to $\zeta_2$'}
\label{eq:sentence}
\end{equation}
corresponds to the complement to the set 
$\cO\bigl(\{(\xi,\xi'),(\zeta,\zeta')\},\{1,0\}\bigr)$. This complement is closed.
We take intersection of all sets \eqref{eq:sentence}, and get 
the semigroup $\GammA_\Lambda$.
\hfill $\square$

\begin{proposition}
 The group $\G_\Lambda$ is dense in $\GammA_\Lambda$.
 \label{pr:dense}
\end{proposition}

Before the proof,  we  reduce the family of sets 
$\cO(\{\omega_j,\omega'_j)\},\{\delta_j\})$ defining a base of the topology.
 First of all, 
we can assume that $\omega_j$, $\omega_j'$ 
correspond to perfect balls of the same type $(\epsilon,r)$ (otherwise,
the condition 'a partial isomorphism $\pi$  does not send $\omega$ to $\omega'$' holds
automatically). Second, if $d(\omega_1,\omega_2)\ne d(\omega_1',\omega_2')$,
then
 $$\cO(\{(\omega_1,\omega_1'),(\omega_2,\omega_2')\};\{1,1\})=\varnothing.$$
On the other hand, under the same condition
$$
\cO(\{(\omega_1,\omega_1'),(\omega_2,\omega_2')\};\{1,0\})=
\cO(\{(\omega_1,\omega_1')\};\{1\}).
$$
Third, if $\omega_1=\B^\delta(x,\rho)$ is over $\omega_2= \B^\epsilon(x,r)$, $\omega_1'=\B^\delta(x',\rho)$ is over $\omega_2'=\B^\epsilon(x',r)$,
then 
$$
\cO(\{(\omega_1,\omega_1'), (\omega_2,\omega_2')\},\{0,1\})=\varnothing.
$$

For these reasons, the following family $\cO[\dots]$ form a base of the topology
in $\GammA_\Lambda$. We take two finite collections of points $\{\xi_\alpha\}$,
$\{\xi'_\alpha\}$ corresponding to two isometric systems of perfect balls.
Equivalently, $\xi_\alpha$, $\xi_\alpha'$ are not contained 
in naked intervals and woolly trees generated by these collections are
isomorphic.  Next, we take finite collections $\{\eta_\beta\}$, $\{\eta_\beta'\}$
such that $\eta_\beta$ and $\eta_\beta'$ correspond to perfect balls
of one type and for each $\alpha$, $\beta$ we have
$d(\xi_\alpha,\eta_\beta)=d(\xi_\alpha',\eta_\beta')$. Also there are no points
$\eta_\beta$ lying over $\xi_\alpha$.
Under these conditions we define sets 
$\cO\bigl[\{\xi_\alpha\},\{\xi'_\alpha\};\{\eta_\beta\},\{\eta_\beta\}\bigr]$
consisting of partial isomorphisms $\pi$ of woolly subtrees
such that $\pi$ sends $\xi_\alpha$ to $\xi_\alpha'$ for each $\alpha$ and does not
send $\eta_\beta$ to $\eta_\beta'$ for each $\beta$.

\sm

{\sc Proof of Proposition \ref{pr:dense}.}
Let $\pi\in\GammA_\Lambda$. Consider a 
finite collection of points $\{\xi_\alpha\}\in \dom \pi$, let
$\xi_\alpha':=\pi(\xi_\alpha)$. Choose arbitrary $\eta_\beta$, $\eta_\beta'$
satisfying the conditions of the previous paragraph.
We must show that the set $\cO\bigl[\{\xi_\alpha\},\{\xi'_\alpha\};\{\eta_\beta\},\{\eta'_\beta\}\bigr]$ contains a point of $\G_\Lambda$.

 Without loss of generality we can assume that $\xi_\alpha'=\xi_\alpha$,
 otherwise we multiply $\pi$ by an appropriate element of $\G_\Lambda$.
 Denote by $[[\{\xi_\alpha\}]]$  the minimal woolly subtree containing these points.
 Let $\G[[\{\xi_\alpha\}]]\subset \G_\Lambda$ be the stabilizer of this subtree.
 It is sufficient to show that there is an element $g$ of this stabilizer,
such that $\eta_\beta g\ne  \eta_\beta'$ for all $\beta$.

As in Subsect. \ref{ss:stabilizers}, our woolly subtree determines a splitting
of $\U_\Lambda$ into disjoint clopen sets isometric to perfect balls.
For each vertex $v$ of $[[\{\xi_\alpha\}]]$
we denote by $\germ^\downarrow_{\out}(v)$ the set all 
elements $\nu$ of  $\germ^\downarrow(v)$ that are contained in 
the woolly subtree $[[\{\xi_\alpha\}]]$ and look to outside the subtree
(i.e., there are no elements of $[[\{\xi_\alpha\}]]$ under $\nu$).
Each element  $\nu\in \germ^\downarrow_{\out}(v)$ determines an 
open ball  $D_\nu\subset\U_\Lambda$ lying under $\nu$.
For a vertex $v\in [[\{\xi_\alpha\}]]$ we define the set
$\wh D_v$ consisting of all $x\in\U_\Lambda$ such that the lowest vertex of 
$L_x \cap [[\{\xi_\alpha\}]]$ is $v$. Next, we define the set
$$
D_v:=\wh D_v\setminus \cup_{\nu\in \germ^\downarrow_{\out}(v)} D_\nu
$$
Then $\U_\Lambda$ splits into a disjoint union
$$
\coprod_{v\in \ver [[\{\xi_\alpha\}]]} \Bigl(D_v \coprod\amalg
 \Bigl(\coprod_{\nu\in \germ^\downarrow_{\out}(v)} D_\nu\Bigr)\Bigr).
$$  
The stabilizer $\G[[\{\xi_\alpha\}]]$ consists of
elements of $\G_\Lambda$ preserving this partition 
and is isomorphic to
$$
\prod_{v\in \ver [[\{\xi_\alpha\}]]} \Bigl(\Isom(D_v) \times
 \Bigl(\prod_{\nu\in \germ^\downarrow_{\out}(v)} \Isom(D_\nu)\Bigr)\Bigr).
$$ 

Under our conditions, each pair of balls $\eta_\beta$, $\eta_\beta'$ is contained
in one set $\wh D_v$. 
If balls $\eta_\gamma $, $\eta_\gamma '$ are contained in different elements
of the subpartion  $\wh D_v=D_v\amalg (\amalg_\nu D_\nu)$, then for any element $q$ of the stabilizer 
$q(\eta_\gamma)$ and $\eta_\gamma '$ remain to be in the same elements, and $q(\eta)\ne\eta'$.
 Let us   separate remaining pairs $(\eta_\beta,\eta_\beta')$.

Now we are in the following situation. Each set $D_v$ or $D_\nu$ is isometric to a certain perfect ball,
and it contains two finite families of subballs corresponding to $\eta_\beta$ and $\eta_\beta'$.
 So consider a perfect ball $\B^\sigma(z,R)$
and two collection of its proper perfect subballs $\{\B^{\epsilon_j}(x_j,r_j)\}$,
$\{\B^{\epsilon_j}(y_j,r_j)\}$. We find in this family a subball of the largest 
type, say $\B^\delta(z,\rho)$ and replace our collections by
 $\{\B^\delta(x_j,\rho)\}$, $\{\B^\delta(y_j,\rho)\}$. Since there is a countable
 set of subballs of this type  in $\B^\sigma(z,R)$, we can choose
 $q\in \Isom(\B^\sigma(z,R))$ such that $q(\B^\delta(x_j,\rho))\ne \B^\delta(y_j,\rho)$.
 \hfill $\square$

 \sm
 
 {\bf \punct Representations of semigroups $\boldsymbol{\GammA_\Lambda}$.}
 
 \begin{theorem}
 \label{th:extension}
 {\rm a)} Any unitary representation $\rho$ of the group $\G_\Lambda$ extends
 by continuity to a $*$-representation $\ov\rho$ of the semigroup
 $\GammA_\Lambda$. Moreover, the weak closure $\ov{\rho(\G_\Lambda)}$
 of the set of all operators $\rho(g)$, where $g$ ranges in $\G_\Lambda$,
 coincides with $\ov\rho(\GammA_\Lambda)$.
 \end{theorem}

Notice that a continuous extension is unique if it exists.

\sm 

{\sc Proof.} a)
 Denote by $\Xi_{\epsilon,r}$ the set of all perfect
balls in $\U_\Lambda$ of a type $(\epsilon,r)$. By the definition
the semigroup  $\GammA_\Lambda$ acts on each  space $\Xi_{\epsilon,r}$ by partial bijections.
Therefore $\GammA_\Lambda$ acts in $\ell^2(\Xi_{\epsilon,r})$
by 0-1-matrices. Consider an action of $\GammA_\Lambda$ in some tensor product
$\otimes_{j=1}^N\ell^2(\Xi_{\epsilon_j,r_j})$.
 Since the group $\G_\Lambda$ is dense in $\Gamma_\Lambda$, each $\G_\Lambda$-invariant subspace
 is $\GammA_\Lambda$-invariant. By
 Theorem \ref{th:classification-main}$^\circ$  $\GammA_\Lambda$ 
such tensor  products
 contain all irreducible
representations of $\G_\Lambda$.
So $\Gamma_\Lambda$ acts in any irreducible representation and therefore in any direct integral of irreducible representations.

\sm 

b) Indeed, a continuous image of a compact set is closed.
\hfill $\square$

\sm

{\sc Remark.}  Our claim about compactness 
of the semigroup $\GammA_\Lambda$ is based on the Tychonoff theorem 
(direct product of any family of compact spaces is compact, in our case
a family is continual). The general Tychonoff theorem is equivalent to the Axiom of choice; if factors are Hausdorff spaces, then the compactness of the product
follows from the Ultrafilter principle (any non-trivial filter can be embedded
to an ultrafilter), which is weaker than the Axiom of choice, see, e.g., \cite{Sche}, 17.16, 17.22. 
This dependence  in a context, which seems constructive,
seems strange to the author. It is interesting to find a proof using only
the Axiom of dependent choice (on discussion of this axiom, see, e.g., \cite{Sche}, 6.28, 14.77, 14.78.d, 17.33, 19.51, 20.16, 27). 
\hfill $\boxtimes$

 \section{Affine isometric actions and property (T)\label{s:affine}}

\COUNTERS

{\bf \punct Affine isometric actions.}
Let $H$ be a real Hilbert space. An affine isometric transformation
of $H$  has the form
$$
R(g)h= Uh+\gamma,
$$
where $h$ is an element of the Hilbert space, $U$ is an orthogonal operator ($UU^*=1=U^*U$), and $\gamma\in H$.
The condition
$$
R(g_1)R(g_2)=R(g_1g_2)
$$
is equivalent to a pair of conditions:

\sm

--- $U(g_1)U(g_2)=U(g_1g_2)$, i.e., $U(g)$ is a representation of $G$;

\sm

---  the '{\it cocycle}'
 $\gamma(g)$ satisfies the identity
 $$
 \gamma(g_1 g_2)=R(g_1)\gamma(g_2)+\gamma(g_1).
 $$
 
 \sm
 
It can happened that there exists a global fixed point $p$ for all such transformations.
Then we can shift  0 of the Hilbert space to $p$ and get a usual unitary linear representation.
In this case for any point $h$ we have 
$$\|R(g)h-h\|\le \|R(g)h-p\|+\|p-h\|=|R(g)h-R(g)p\|+\|p-h\|=2\|h-p\|$$
and so a $G$-orbit of each vector $h$ is bounded.

We say that an affine isometric action is {\it non-trivial}, if it has no fixed points.

\sm

{\bf \punct Affine isometric actions of the groups $\boldsymbol{\G_\Lambda}$.}

\begin{proposition}
Let a set $\Lambda$ have not a maximal element.
Then the group $\G_\Lambda$ has nontrivial affine isometric actions.
\end{proposition}

{\sc Proof.} First, let $\Lambda$ be unbounded.
We consider an isometric embedding $\psi$ of $\U_\Lambda$ to a real Hilbert space $H$,
which exists according Theorem \ref{th:embedding}.
We can assume that  there is no a proper affine subspace in $H$ containing the image $\psi(\U_\Lambda)$.
Otherwise, we pass to a minimal subspace, which contains $\psi(H)$. Then any isometry
of the embedded space can be extended in a unique way to an isometry of $H$.
So we get an affine isometric action of $\G_\Lambda$ on $H$.
An orbit of any point $\psi(x)$ is the whole set $\psi(\U_\Lambda)$. So orbits are unbounded
and our action has  not fixed points.

Second, let $\Lambda$ be bounded and $a$ be the supremum of $\Lambda$.
Then we consider a continuous monotone bijection $\kappa:[0,a)\to[0,\infty)$
and set $\wt\Lambda=\kappa(\Lambda)$. Then groups $\G_{\wt\Lambda}$ and $\G_\Lambda$
coincide and we refer to the previous paragraph.
\hfill $\square$

\sm 

{\sc Remarks.} a) The group $\Isom(H)$ has a series of unitary representations. For this reason, affine isometric actions is a tool of construction of unitary representations
of infinite-dimensional groups, consideration of such actions is necessary in 
investigation of many groups or classes of groups. Apparently, this way to construct unitary representations
was proposed by Araki \cite{Ara}, see some constrictions in this spirit exposed or refered in the book
\cite{Ner-book},  Sect. 9.5, 9.7, 10.1-10.4, see occasional appearance of this topic in our list of references 
\cite{Ner-bist}, \cite{Ner-field}, Sect. 3.14, \cite{Olsh-new}, Subsect. 5.5, \cite{Olsh-GB}, Sect. 19. Some representations
obtained in this way can be a point of a separate interest, see \cite{GGV}. 

\sm

b)
According the construction of Subsect. \ref{ss:ultra-hilbert}, we have  representations, say $R_s$, of the group
$\G_\Lambda$ in  Hilbert  spaces $\cE_s(\U_\Lambda)$.
It is easy to show, that these representations are reducible and are direct integrals of representations in the spaces
$\ell^2(\Xi_{\epsilon,r})$, such direct integrals can be written explicitly. We omit these details. See a continuation of this topic
in next section.
\hfill $\boxtimes$



\sm
   
{\bf \punct Property $\boldsymbol{(T)}$ for the groups $\boldsymbol{\G_\Lambda}$.}
The {\it Kazhdan property $(T)$} is the following condition: whenever a unitary representation 
of $G$ weakly contains%
\footnote{See, e.g., \cite{Kir}, 7.3, \cite{Dix}, 2.2.4, 13.1.3, \cite{BHV}, Appendix F.} the trivial representation, it contains the trivial representation as 
a subrepresentation (see different variants of the definition and further discussion in a book by Bekka, de la Harpe, and Valette \cite{BHV}).

\begin{proposition}
A group $\U_\Lambda$ has property $(T)$ if and only if $\Lambda$ has a maximal element.
\end{proposition}

{\sc Proof.}  If $\Lambda$ has not a maximal element, then it has a non-trivial affine isometric
action and therefore has not property $(T)$, see, e.g., \cite{BHV}, Theorem 2.12.4.

Let $\mu$ be the maximal element of $\Lambda$.  Then our space $\U_\Lambda$
is a closed ball $\B^c(q,\mu)$ of radius $\mu$. Let $v$ be the vertex of $T(\U_\Lambda)$ corresponding 
this ball. 
As in  Subsect. \ref{ss:stabilizers},
we observe that the group $\G_\Lambda$  is a wreath product
\begin{equation}
S_\infty(\Z)\ltimes \Bigl(\prod_{j\in \Z} \Isom(\B^o(p_j,\mu))\Bigr).
\label{eq:wreath}
\end{equation}

Let a unitary  representation $\rho$ of $\G_\Lambda$ weakly contain the trivial representation.
In particular, it weakly contains the trivial representation of $S_\infty(\Z)$.
The group $S_\infty$ has property $(T)$ (this easily follows from the Lieberman classification theorem,
or we can refer to a general theorem about Roelcke precompact groups \cite{Iba}).
Therefore $\rho$ contains $S_\infty(\Z)$-fixed vectors. We literally
repeat the proof of Lemma \ref{l:fixed-fixed} and get that such vectors are fixed
with respect to the whole wreath product \eqref{eq:wreath}.
\hfill $\square$


\section{Groups of isometries of Bruhat-Tits 
R-trees
 of non-Archimedean 
 non-locally finite fields \label{s:extension}}

\COUNTERS

Let as in Subsect. \ref{ss:urysohn-construction}, $\Sigma\subset \R$ be a countable dense subgroup, let $h$ be an element of the interval $(0,1)$,
and $\Lambda$ be the set of all numbers of the form $h^s$, where $s$ ranges in $\Sigma$.
Let $\Bbbk$ be a countable field,
let $\K=\K(\Sigma,\Bbbk)$ be the corresponding field of formal Puiseux series \eqref{eq:puizo} considered as model of the Urysohn space $\U_\Lambda$.
 Consider the $\R$-tree
$T(\K)$ and define a new metric $m(\cdot,\cdot)$ on the set $T_{>0}(\K)$. Let $\xi$, $\eta\in T(\U_\Lambda)$, let $\upsilon$ be the upper point of the segment $[\xi,\eta]$.
We set
$$
m(\xi,\eta)= -\log_h(\frh(\upsilon)/\frh(\xi))-\log_h(\frh(\upsilon)/\frh(\eta)).
$$
We also can say that an element of our new tree $\T_\Sigma$ are pairs $(S, f(t))$, where $S\in\R$ and
$f(t)$ is a Puiseux  polynomial
$$
f(t)=\sum_{j=1}^N a_{s_j} t^{s_j},\qquad \text{where $a_j\in \Bbbk$, $s_j\in \Sigma$, and $s_j<S$.}
$$
For two elements $\xi_1=(S_1, f_1(t))$, $\xi_2=(S_2, f_2(t))$ we take
the nonzero term $c_st^s$ of $(f_1-f_2)$ of lowest power $s$
and set
$$
m(\xi_1,\xi_2)=(S_1-s)+(S_2-s).
$$
Now points of $\K=\U_\Lambda$ do not contained in $\T_\Sigma$ (they become infinitely far), but the group of isometries
$\Isom(\T_\Sigma)$ contains $\G_\Lambda$ and is larger than $\G_\Lambda$. 
In particular, for each vertex $v$ the set $\germ(v)$ is homogeneous with respect to the stabilizer $\Isom^v(\T)$ of $v$ in $\Isom(\T)$.
Moreover, the natural map $\Isom^v(\T)\to S_\infty(\germ(v))$ is surjective.

\sm 

{\sc Remark.}
This $\R$-tree $\T_\Xi$ is well-known and it was one of reasons of introduction
of the notion of $\R$-trees, see e.g., \cite{Chi}. The group $\GL(2,\K)$ naturally acts on the projective
line $\mathbb{P}\K^1=\K\cup\infty$ by linear fractional transformations, and these transformations have canonical extension
to  the tree $\T_\Xi$. \hfill $\boxtimes$

\sm



\begin{conjecture}
An irreducible unitary representation of the group $\Isom(\T_\Sigma)$ is either a representation
in a Hilbert space $\cE_s(\T)$ {\rm(}see Corollary {\rm\ref{cor:olsh})}, or is induced from
a stabilizer of a finite woolly subtree%
\footnote{Cf. Olshanski's classification  theorems for groups of automorphisms of Bruhat-Tits trees \cite{Olsh-bruhat}, \cite{Olsh-new}.}
of $\T_\Sigma$.
\end{conjecture}

By  'finite woolly subtrees' we understand the following objects.
 We consider a finite collection of points $\xi_\alpha\in \T_\Sigma$ and its convex
hull $\cT$. We equip $\cT$ with a finite collection of germs $\in \coprod_{v\in(\ver \T_\Sigma)\cap \cT)}\germ(v)$
looking to outside of $\cT$.

\tt

Fakult\"at f\"ur Mathematik, Universit\"at Wien;

Institute for Information Transmission Problems;

Moscow State University

e-mail: yurii.neretin(frog)univie.ac.at

URL: mat.univie.ac.at/$\sim$neretin

\end{document}